\documentclass[11pt]{article}
\input epsf
\usepackage{amsmath}
\usepackage{amssymb}
\usepackage{amscd}
\newcommand{\qed}{\hbox{\unskip\nobreak\hfil
        \penalty50\hskip1em\hbox{}\nobreak\hfil
        $\square$\parfillskip=0pt\finalhyphendemerits=0 \par}}

\newtheorem{dfn}{Definition}[section]
\newtheorem{rem}[dfn]{Remark}
\newtheorem{Not}[dfn]{Notation}
\newtheorem{conv}[dfn]{Convention}
\newtheorem{thm}[dfn]{Theorem}
\newtheorem{lem}[dfn]{Lemma}
\newtheorem{prop}[dfn]{Proposition}
\newtheorem{cor}[dfn]{Corollary}

\newtheorem{ex}[dfn]{Example}

\newtheorem{defn}[dfn]{Definition}

\def\no{\noindent}
\def\nbds{neighborhoods~}
\def\qh{quasi-homogeneous\ }
\def\={=\joinrel =}
\def\i{\iota}
\def\0{\emptyset}
\def\half{\frac{1}{2}}
\def\proof{\par\medskip\noindent{\it Proof: }}
\def\lra{\longrightarrow}

\def\>{\rangle}
\def\<{\langle}

\def\C{\mathbb C}
\def\t{\tilde}
\def\R{\mathbb R}
\def\P{\mathbb P}
\def\Z{\mathbb Z}
\def\E{{\cal E}}

\def\a{{\cal A}}
\def\Q{\mathbb Q}

\def\V{{\cal V}}
\def\A{{\mathbb A}}
\def\M{{\cal M}}
\def\eps{\epsilon}
\def\H{\mbox{\rm H}}

\def\al{\alpha}
\def\be{\beta}
\def\ga{\gamma}
\def\Ga{\Gamma}
\def\del{\delta}
\def\Del{\Delta}
\def\Si{\Sigma}
\def\si{\sigma}
\def\L{{\cal L}}
\def\r{{\cal R}}
\def\la{\lambda}
\def\La{\Lambda}
\def\8{\infty}

\def\om{\omega}
\def\D{\partial}
\def\M{{\cal M}}
\def\N{{\cal N}}

\def\E{{\cal E}}
\def\J{{\cal J}}
\def\gg{{\mathfrak g}}

\def\nbd{neighborhood~}
\def\hook{\hookrightarrow}
\def\k{\hbox{\bf k}}
\def\G{\hbox{\bf G}}
\def\1{\hbox{\bf 1}}
\def\m{{\mathfrak m}}

\def\u{\sqcup}
\def\h{/\!/}
\def\ov{\overrightarrow}
\def\BE{\begin{equation}}
\def\EE{\end{equation}}
\def\bul{\bullet}
\def\dga{differential graded algebra~}
\def\dgla{differential graded Lie algebra~}
\def\dglas{differential graded Lie algebras~}
\def\RA{\joinrel\relbar\joinrel\relbar\joinrel\relbar\joinrel\relbar\joinrel\relbar\joinrel\rightarrow}
\def\Hom{\mbox{\rm Hom}}
\def\BHom{\mbox{\rm BHom}}
\def\p{{\cal P}}
\def\LA{\leftarrow\joinrel\relbar\joinrel\relbar\joinrel\relbar\joinrel\relbar\joinrel\relbar\joinrel}

\begin{document}
\title{On representation varieties of Artin  groups, projective arrangements 
and the fundamental groups of smooth complex algebraic varieties\footnote{Preliminary version.}}
\author{Michael Kapovich\thanks{This research was partially supported by 
NSF grant DMS-96-26633 at University of Utah.} \ 
and John J.\ Millson\thanks{This research was partially supported by 
NSF grant DMS-95-04193 at  University of Maryland.}
}
\date{February 8, 1997}
\maketitle

\begin{abstract}
\no We prove that for any affine variety $S$ defined over $\Q$ 
there exist Shephard and Artin groups $G$ such that a Zariski open 
subset $U$ of $S$ is biregular 
isomorphic to a Zariski open subset  of the character  
variety $X(G, PO(3))= \Hom(G,  PO(3))\h PO(3)$. The subset $U$ contains 
all real points of $S$. As an application we construct new examples 
of finitely-presented groups which are not fundamental groups of 
smooth complex algebraic varieties. 
\end{abstract}

\section{Introduction} 

The goal of this paper is to understand representation varieties 
of Artin and Shephard groups and thereby obtain information on Serre's 
problem of determining which finitely-presented groups are fundamental 
groups of smooth complex (not necessarily compact) 
algebraic varieties. The first  
examples of finitely-presented  groups which are not fundamental 
groups of smooth complex algebraic varieties were given by 
J.~Morgan \cite{Morgan1}, \cite{Morgan2}. 
We find a new class of such examples which consists of certain Artin and 
Shephard groups. Since all Artin and Shephard groups have quadratically 
presented Malcev algebras, Morgan's test does not suffice 
to distinguish Artin groups from fundamental groups of smooth complex 
algebraic  varieties or even from fundamental groups of 
compact K\"ahler manifolds, see \S \ref{completion}.

Our main results are the following theorems (Artin and Shephard 
groups are defined in \S \ref{groups} below): 

\begin{thm}
\label{moremorgan} 
There are infinitely many distinct Artin groups that  are not isomorphic 
to fundamental groups of smooth complex algebraic varieties.
\end{thm}

\begin{thm}
\label{morereps}
For any affine variety $S$ defined over $\Q$ 
there are Shephard and Artin groups $G$ such that a Zariski open 
subset $U$ of $S$ is biregular 
isomorphic to  Zariski open subsets of the character  
varieties $X(G, PO(3))= \Hom(G,  PO(3))\h PO(3)$. 
The subset $U$ contains all real points of $S$. 
\end{thm}

The surprising thing about Theorem \ref{moremorgan} is that Artin groups 
look very similar to the fundamental groups of smooth complex 
quasi-projective varieties. For example the free group on $n$ letters is the fundamental group of $\C$ with $n$ points removed and it is the Artin group 
associated with the graph with $n$ vertices and no edges. On the other 
extreme, take a finite 
complete  graph where each edge has the label $2$. The 
corresponding Artin group is free Abelian, hence it is the 
fundamental group of the quasi-projective variety $(\C^{\times})^n$. 
Yet another example is the braid group which is the Artin group 
associated with the permutation group $S_n$. Theorem \ref{moremorgan} 
is a consequence of Theorems \ref{moremnev}, \ref{t3}, \ref{t4}, 
\ref{t6}, \ref{t7} and Corollary \ref{c5}
below. The main body of this paper is concerned with a study of the following diagram\footnote{Here and in what 
follows we do not assume our varieties are reduced or irreducible, i.e. they are schemes locally of finite type.}: 

\[
\begin{array}{ccc}
\framebox{\parbox{2in}
{Moduli spaces of arrangements in the projective plane}}
& \begin{array}{c}
\stackrel{\tau}{\RA}\\  
\stackrel{\displaystyle \LA}{\scriptstyle geo}
\end{array}
& \framebox{\parbox{.6in}{Varieties over $\Q$}} \\
  {\ }_{alg\Bigg\downarrow} & ~  & ~ \\ 
~ & ~ & ~ \\
\framebox{\parbox{2in}{Character varieties of representations
of  Shephard groups 
into $PO(3)$}} & \stackrel{\mu}{\RA} 
&  \framebox{\parbox{2in}
{Character varieties of representations of Artin groups into $PO(3)$}}
\end{array}
\]

The arrow $\tau$ is tautological, the arrow $\mu$ is pull-back. 
The arrows $geo$ and $alg$ are defined below. In 
\S\S \ref{Abstract arrangements}, \ref{Projective  arrangements} 
we define {\em abstract arrangements} (essentially bipartite graphs) 
and their projective realizations. The space $R(A)$ of projective 
realizations of a given abstract arrangement has a canonical 
structure of a projective variety (i.e. projective scheme, see 
the preceeding footnote) over $\Q$. In fact we refine the notion 
of projective arrangement to obtain the affine variety $BR_0(A)$ 
of {\em finite based realizations}. The variety $BR_0(A)$ injects 
as an open subvariety into 
the {\em moduli space} ${\cal M}(A)=BR(A)\h PGL(3)$ of the arrangement $A$. 
Our version  of Mnev's theorem \cite{Mnev} is then 

\begin{thm}
\label{moremnev}
For any affine algebraic variety $S$ defined over $\Q$ there is a 
marked based abstract arrangement $A$ such that the varieties 
$BR_0(A)$, $S$ are isomorphic. 
\end{thm}

\begin{rem}
It appears that Mnev's theorem \cite{Mnev} implies only that there is a 
stable homeomorphism between the sets of real points of $BR_0(A)$ and $S$. 
In addition Mnev gives only an outline of the proof. For our application 
to Serre's problem it is critical to prove an isomorphism on the scheme 
level. 
\end{rem}

The key idea  in proving Theorem \ref{moremnev} is to construct a 
cross-section 
$geo$ to $\tau$ (over the category of affine varieties) 
by showing that one can do ``algebra via geometry'', 
that is one can describe elementary algebraic operations over any 
commutative ring using projective arrangements 
(see \S \ref{Algebraic operations}). This idea actually goes back to the 
work of von Staudt \cite{St} (the ``fundamental theorem of the 
projective geometry''). The abstract arrangement $A$ corresponding 
to $S$ under $geo$ depends upon a choice of 
affine embedding $\t{S}$ (i.e. defining equations) for $S$ 
and upon a choice of particular formulae 
describing these equations.  Moreover we obtain isomorphism $\t{S} \to BR_0(A)$ of affine schemes over $\Q$. We will abuse notation and use $geo$ to denote this isomorphism as well. Thus the symbol $geo$ has two meanings: $A= geo(S)$ 
  is the arrangement determined by the above data, but if $x\in \t{S}$ then $\psi= geo(x)$ is a point in $BR_0(A)$. The meaning to be assigned to $geo(\cdot)$ will be clear from its argument. We make analogous abuses of notation with mappings $\tau$ and $alg$ below. 
 
The arrow {\em alg} encodes the points and lines of an {\em anisotropic  
projective realization} (see  \S \ref{S14} for the definition)  
$\psi$ of the abstract arrangement $A$ into a representation 
$$
\rho= \rho_{\psi}: G^s_A \lra PO(3, \C)
$$ 
of the Shephard group $G^s_A$ associated to the abstract arrangement 
$A$. A choice of a nondegenerate bilinear form on $\C^3$ 
determines anisotropic points and lines in $\P^2$ (we choose the 
bilinear form  so that all real points of $\P^2$ are anisotropic). 
Each anisotropic point in $\P^2$  determines the Cartan involution 
in $PO(3, \C)$ around this point or the rotation of order $3$ having 
this point as neutral fixed  point (i.e. a fixed  point where the 
differential of the rotation has determinant 1). There are 
two such rotations of order $3$, we choose one of them. Similarly 
every anisotropic  line $L$ 
uniquely determines the reflection in $PO(3, \C)$ which keeps 
$L$ pointwise fixed. Finally one can encode the incidence relation 
between points and lines in $\P^2$ using algebra: 
two involutions generate the subgroup $\Z/2\times \Z/2$ in $PO(3,\C)$ 
iff the neutral fixed point of one belongs to the fixed line of the 
other, rotations of orders $2$ and $3$  anticommute iff the 
neutral fixed point of the rotation of order $3$ belongs to the
fixed line of the involution, etc. We get a morphism  
$$
alg: \hbox{~~based anisotropic arrangements~~} \longrightarrow 
\hbox{~~representations}
$$
In the following theorem we shall identify $alg(\psi)$ with its 
projection to the character variety 
$$
X(G^s_A, PO(3)):= \Hom(G_A^s, PO(3))\h PO(3)\ \quad . $$ 

\begin{thm}
\label{t3}
The mapping $alg: BR(A, \P^2_0)\to X(G_A^s, PO(3))$ 
is a biregular  isomorphism onto a Zariski open (and closed) 
subvariety $\Hom_f^+ (G^s_A, PO(3))\h PO(3)$.
\end{thm}

\smallskip
The mapping $alg$ has the following important property: Let $S$ be an 
affine variety defined over $\Q$ and $O\in S$ be a rational point. 
Then we can choose an arrangement $A$ so that $O$ corresponds to a 
realization $\psi_0$ under the mapping $geo: S\to BR_0(A)$ 
such that the image of the representation 
$alg (\psi_0)$ is a finite subgroup of $PO(3,\C)$ with trivial 
centralizer. 

There is an Artin group $G^a_A$ and a canonical epimorphism  
$G^a_A\to G^s_A$ associated with the Shephard group $G^s_A$. 
It remains to examine the  morphism 
$\mu:  \Hom^+_f(G^s_A, PO(3)) \to \Hom (G^a_{A}, PO(3))$ 
given by pull-back of homomorphisms. 

\begin{thm}
\label{t4}
Suppose that $A$ is an admissible based arrangement. Then the morphism 
$\mu$ is an isomorphism onto a union of Zariski connected components. 
\end{thm}

\begin{cor}
\label{c5} 
The variety $\Hom (G^a_A, PO (3))$ inherits all the
singularities of the representation variety $\Hom (G^s_A, PO (3))$
corresponding to points of $BR(A, \P^2_0)$, whence (since all 
real points of $BR(A)$ are anisotropic) to all singularities 
of $BR(A)$ at real points.
\end{cor}

Combining Corollary \ref{c5} with Theorem \ref{moremnev} we obtain

\begin{thm}
\label{t6}
Let $S$ be an affine algebraic variety defined over $\Q$ and 
$O\in S$ be a rational point. Then there exists an admissible 
based arrangement $A$ and a 
representation $\rho_0: G^a_A \to PO(3,\R)$ with finite image such that 
the (analytical) 
germ $(X(G^a_A, PO(3, \C)), [\rho_0])$  is isomorphic to $(S, O)$. 
\end{thm}

Thus the singularities of representation varieties of Artin groups 
at representations with finite image are at least as bad as germs 
of affine varieties defined over $\Q$ at rational points. 

On the other hand, if $M$ is a (connected) smooth complex algebraic  
variety and $G$ is an algebraic Lie group, the singularities of 
$\Hom(\pi_1(M), G)$ at representations with finite image are 
severely limited by Theorem \ref{t7} below. We will need the following

\begin{defn}
Let $X$ be a real or complex analytic space, $x\in X$ and $G$ a Lie group acting on $X$. We say that there is a local {\em cross-section} through $x$ to the $G$-orbits if there is a $G$-invariant open \nbd $U$ of $x$ and a closed analytic subspace $S\subset U$ such that the natural map $G \times S \to U$ is an isomorphicm of analytic spaces.    
\end{defn}  

\begin{thm}
\label{t7}
 Suppose $M$ is a smooth connected complex algebraic variety, $G$ is an 
algebraic  Lie group  and $\rho: \pi_1(M)\to G$ is a representation 
with finite image. Then the germ  
$$(\Hom(\pi_1(M), G), \rho)$$
 is analytically isomorphic to a 
quasi-homogeneous cone with generators of weights $1$ and $2$ and 
relations of weights $2, 3$ and $4$. 
In the case there is a local cross-section through $\rho$ to $Ad(G)$-orbits, then 
the same conclusion is valid for the quotient germ 
$$(\Hom(\pi_1(M), G)\h G, [\rho])\quad . $$  
\end{thm}
  
We present two proofs of this result: in  \S \ref{Hain's theorem} 
we deduce it  from a theorem of  R.~Hain \cite{Hain}
%(under the extra condition that $\H^0(\pi_1(M), ad\circ\rho)=0$)
 and, since Hain's paper is still in preparation, in \S \ref{direct} 
we also deduce Theorem \ref{t7} from results of J.~Morgan \cite{Morgan2} 
on Sullivan's minimal models of smooth complex algebraic varieties. 

Note that in the case when $M$ is a {\em compact} smooth 
complex algebraic variety (or more generally a compact K\"ahler manifold) 
then a more general theorem holds:  

\begin{thm}
\label{gms}
(W.~Goldman, J.~Millson \cite{GM}, C.~Simpson \cite{Si}). Suppose that $M$ is a compact K\"ahler manifold $G$ is an algebraic  Lie group  and 
$\rho: \pi_1(M)\to G$ is a representation such that the Zariski closure of $\rho(\pi_1(M))$ is a reductive subgroup of $G$. Then the germ  
$$(\Hom(\pi_1(M), G), \rho)$$
 is analytically isomorphic to a (quasi)-homogeneous cone with 
generators of weight $1$ and relations of weight $2$ (i.e. a 
homogeneous quadratic cone). In the case there is a local cross-section 
through $\rho$ to $Ad(G)$-orbits, then 
the same conclusion is valid for the quotient germ 
$$(\Hom(\pi_1(M), G)\h G, [\rho])\quad . $$  
\end{thm}

Our proof of Theorem \ref{t7} is in a sense analogous to the proof of Theorem 
\ref{gms} in \cite{GM}, \cite{Si}: we construct a differential graded Lie 
algebra ${\cal Q}^{\bul}$ which is weakly equivalent to the algebra of bundle-valued differential forms $\a^{\bul}(M, adP)$ on $M$ so that  
${\cal Q}^{\bul}$ {\em controls} a germ which is manifestly a 
quasi-homogeneous cone with the required weights. 

In Figure \ref{Fig9} we describe the graph of an Artin group $G^a_A$ which 
admits a representation with finite image $alg(\psi_0)= \rho_0: G^a_A\to PO(3,\C)$ such that the germ $(X(G^a_A, PO(3,\C)), [\rho_0])$ 
is isomorphic to the germ at $0$ defined by $x^5=0$. Thus Theorem 
\ref{t7} implies that the group $G^a_A$ is not the fundamental 
group of a smooth complex algebraic variety. 

\begin{rem}
Our convention for Coxeter graphs is different from the standard convention 
for Dynkin diagrams. Namely, if two vertices are not connected by an edge 
{\bf it does not mean that corresponding generators commute}. If on our 
diagram an edge has no label, we assume that the edge has the label $2$. 
On the diagram for a {\bf Shephard group} if a vertex has no label this 
means that the corresponding generator has {\bf infinite order}. 
\end{rem}

There is a local cross-section to the $PO(3)$-orbit through the 
representation $\rho_0$ (that appears is Theorem \ref{t6}), 
hence we apply Theorem \ref{t7} and conclude that $G^a_A$ is not the  
fundamental group of a smooth complex algebraic variety. To see that 
there are infinitely many distinct examples we may proceed as follows. 

Take  the varieties $V_p:= \{x^p=0\}$, $p\ge 2$ are prime numbers.   
Clearly $V_p$ is not analytically isomorphic to $V_q$ for $q\ne p$. 
Thus for all 
$p\ge 5$ the varieties $V_p$ are not analytically  isomorphic to 
quasi-homogeneous varieties described in the Theorem \ref{t7}. 
Hence the Artin groups  $G^a_{A_p}$ 
corresponding to $V_p$ are not fundamental groups of smooth 
complex quasi-projective varieties.  Note 
that among the groups $G^a_{A_p}$ we have infinitely many ones 
which are not  mutually isomorphic. The reason is that for any 
finitely-generated group $\Ga$ the character
 variety $X(\Ga, PO(3, \C))$ has only finite number of isolated 
singular points   whose germs are isomorphic to one of $V_p$.   
This proves Theorem \ref{moremorgan}. 

In \S \ref{triv} we use the results of \S \ref{completion} to show that 
for every Artin group $\Ga$ and Lie group $G$ the germ 
$$
(\Hom(\Ga, G), \rho)$$
is quadratic where $\rho$ is the trivial representation. 

\medskip
{\bf Acknowledgements.} The first author is grateful to A.\ Vershik 
for a lecture on Mn\"ev's result 
in 1989.  The authors are grateful to E.~Bierstone, J.~Carlson, R.~Hain, 
J.~Kollar, P.~Millman, C.~Simpson  and D.~Toledo for helpful conversations. 

\tableofcontents

\section{Singularities}
\label{singu}

Let $\k$ be a field of zero characteristic.  Let $V$ be a variety 
defined over the field {\bf k}, $o \in V$ be a point and 
$\widehat{{\cal O}_{V, o}}$ the complete local ring. We denote by 
$$
T^m _o (V)= \Hom_{ {\bf k}- alg}(\widehat{{\cal O}_{V, o}}~,  
{\bf k}[t]/t^{m+1})$$
the  $m$-th order Zariski tangent space at $o \in V$ and by 
$\pi_{mn} : T^m _o \to T_o^n (V)$ the natural projection ($m\ge n\ge 1$). 
We will think of elements in $T^m_o$ as formal curves of 
degree $m$ which are tangent to $V$ at $o$ up to the order $m$. 
If $(V, o)$, $(W, p)$ are two germs as above and $f: (V, o)\to (W, p)$ 
is an analytic isomorphism of germs then it induces an 
isomorphism between the corresponding 
towers of Zariski tangent spaces of finite order 

\[
\begin{array}{ccc}
T^m _o (V) & \lra &  T^m _p (W) \\
\pi_{mn}\Big\downarrow & ~ &  \pi_{mn}\Big\downarrow\\
T^n_o (V)  & \lra &  T^n_p (W) \\
\end{array}
\]

\begin{lem}
\label{LL3.1}
Suppose that $f: (V, o)\to (W,p)$ is a morphism of analytic germs such that
the derivatives $D_o ^{(m)}f: T^m _o (V) \to T^m _p (W)$ are isomorphisms 
for all $m$. Then $f$ is an isomorphism of germs.
\end{lem}
\proof Let $R:= \widehat{O_{V,o}}$ and $S:= \widehat{O_{W,p}}$ be the 
complete local rings and let $\m_R$ and $\m_S$ be their maximal ideals. 
By Artin's theorem (see \cite[Theorem 3.1]{GM}) it suffices to prove that 
the induced map $f^*: S\to R$ is an isomorphism. We observe that 
$D_o ^{(n+1)}f$ induces an isomorphism from $\ker(\pi_{n+1,n}: 
T^{n+1}_o(V)\to T^n_o(V))$ to
 $\ker(\pi_{n+1,n}: T^{n+1}_p(W)\to T^n_p(W))$. It is clear that the above 
kernels are canonically isomorphic to the dual vector spaces 
$(\m^n_R/\m_R^{n+1})^*$ and $(\m^n_S/\m_S^{n+1})^*$. Hence $f^*$ induces an 
isomorphism
$$
Gr_n(f^*): \m^n_R/\m_R^{n+1}\lra \m^n_S/\m_S^{n+1}
$$
and consequently the induced map $Gr(f^*)$ of the associated graded rings is 
an isomorphism. The Lemma follows from \cite[Lemma 10.23]{AM}. $\qed$

\begin{rem}
Suppose that $(Z, 0), (W, 0)$ are {\bf minimal} germs of  
varieties in $\A^n$ (i.e. $\k^n$ equals to the both Zariski tangent spaces 
$T_0(Z)$ and $T_0(W)$), and these germs are analytically 
isomorphic. Then there is an analytic diffeomorphism $f: \A^n\to 
\A^n$ whose restriction to $Z$ induces an isomorphism of germs 
$(Z, 0)\to (W, 0)$. See for instance \cite[Proposition 3.16]{Dimca}. 
\end{rem}

\medskip
Suppose that we have a 
collection of polynomials $F= (f_1, ..., f_m)$ in $\k^n$, we assume that 
all these polynomials have trivial linear parts. The  
polynomial $f_j$ is said to be {\em weighted homogeneous} if there 
is a collection of positive integers (weights) $w_1> 0,..., w_n >0$ and 
a number $u_j\ge 0$ so that 
$$
f_j((x_1 t^{w_1}),..., (x_n t^{w_n}))= t^{u_j} f_j(x_1,..., x_n)
$$
for all $t\in\k$. We will call the numbers $w_i$ the {\em weights of 
generators} and the numbers $u_j$  the {\em weigths of relations}. 
Let $Y$ denote the variety given by the system of equations
$$
\{ f_1 =0, ...., f_m =0\}
$$ 
(Note that the germ $(Y,0)$  is necessarily minimal.) 
We say that $(Y,0)$ is a {\em quasi-homogeneous} if we can choose 
generators $f_1,..., f_m$ for its defining ideal such that all the 
polynomials $f_j$ are weighted homogeneous with the same weights 
$w_1,..., w_n$ (we do not require $u_j$ to be equal for distinct 
$j=1,..., m$). In particular, if $(Y,0)$ is a {\em quasi-homogeneous} 
then $(Y,0)$ is invariant under the $\k^{\times}$-action on $\k^n$ 
given by the weights $w_1,..., w_n$. 

\begin{rem}
The variety $Y$ given by a system of quais-homogeneous equations 
$\{ f_j=0\}$ is also called a {\bf quasi-homogeneous} 
(or {\bf weighted homogeneous}) cone. 
\end{rem}

We now give an intrinsic characterization of \qh germs. Suppose that 
$(Y,0)$ is \qh. Let $S_m\subset \k[x_1,..., x_n]$ be the subspace of 
polynomials which are homogeneous of degree $m$ (in the usual sense). 
We may decompose the subspace  $S_m$ into one dimensional eigenspaces 
under $\k^{\times}$ (since the multiplicative group of a field is a 
reductive algebraic group). We obtain a bigrading
\[
\k[x_1,..., x_n]= \oplus_{m,n} S_{m,n}
\]
where $m$ is the degree and $n$ is the weight of polynomial under 
$\k^{\times}$ ($f$ transforms to $t^n f$). We obtain  a new grading of 
$\k[x_1,..., x_n]$ by weight 
\[
\k[x_1,..., x_n]= \oplus_{n=1}^{\8} S_n'
\]
where $S_n'$ is the subspace of polynomials of weight $n$. We let $I$ be 
the ideal of $Y$. Then $I$ is invariant under the action of 
$\k^{\times}$ (since its generators are). We claim
$$
I= \oplus_{n=1}^{\8} I\cap S_n' 
$$
This follows by decomposing the action of $\k^{\times}$ in the 
finite dimensional subspaces $I \cap \oplus_{m=1}^N S_m$. Thus if 
$f\in I$ we may write 
$$
f= \sum_{n=1}^{\8} f_n\ , \quad \hbox{with~~} f_n\in I\cap S_n'
$$
(the sum is of course finite). Let $R= \k[Y]= \k[x_1,.... x_n]/I$. 
Then $R$ is a graded ring, $R= \oplus_{n=1}^{\8} R_n$ with $R_0=\k$, 
$R_n= S_n'/(I\cap S_n')$. 

We let $\hat{R}$ be the  completion of $R$ at $\m$ where $\m$ is the 
ideal of zero, i.e. the ideal generated by $\{x_1,..., x_n\}$. Hence 
$$
\hat{R}\cong \k[[x_1,..., x_n]]/\hat{I}
$$ 
where $\hat{I}$ is the ideal generated by $I$ in $\k[[x_1,..., x_n]]$. 
Hence $\hat{R}= \widehat{O_{Y,o}}$. The ring $\hat{R}$ is not graded but 
it has a decreasing filtration $W^{\bul}$ such that $W^N \hat{R}$ is the 
closure of  $\oplus_{n=N}^{\8} S_n'$. Define $Gr^W_{n}( \hat{R}) = W^{n}(\hat{R})/W^{n+1}(\hat{R})$.  
The filtration $W^{\bul}$ satisfies: 

(i) $Gr^W_0 = \k$.  

(ii) $\cap_{n=0}^{\8}W^n = 0$. 

(iii) $\dim_{\k} Gr^W_n (\hat{R})  < \8$ for all $n$. 

\noindent The inclusion $R\hook \widehat{R}$ induces an 
isomorphism $R\cong Gr^W(\hat{R})$ and we obtain: 

\begin{lem}
\label{crit-qh}
(a) $\widehat{O_{Y,o}}$ admits a decreasing filtration 
$W^{\bul}$ satisfying the  properties (i)---(iii) above. 

(b) There is a monomorphism of filtered rings 
$Gr^W(\widehat{O_{Y,o}} )\to 
\widehat{O_{Y,o}}$ with dense image, so 
$\widehat{O_{Y,o}}$ is the 
completion of $Gr^W(\widehat{O_{Y,o}} )$. 

(c) Conversely if  $\widehat{O_{Y,o}}$ satisfies (a) and 
(b) then $(Y,o)$ is \qh. 
\end{lem}
\proof It remains to prove (c). Define a $\k^{\times}$-action on 
$Gr^W\widehat{O_{Y,o}}$ so that the elements in the $n$-graded 
summand have weight $n$ for the action of $\k^{\times}$. Let $\m$ 
be the ideal of $o$. Choose a basis of eigenvectors under $\k^{\times}$ 
action on $\m/\m^2$. Lift these vectors to eigenvectors $f_1,..., f_n$ of 
$\m$. Then by a standard argument if we set $\pi(x_i)=f_i$  we obtain a 
surjection 
$$
\pi: \k[[x_1,..., x_n]] \lra \widehat{O_{Y,o}}
$$ 
which is $\k^{\times}$-equivariant. Hence the induced map of graded rings
$$
\pi': \k[x_1,..., x_n] \lra Gr^W \widehat{O_{Y,o}}
$$
is also surjective. Let $I$ be the kernel of $\pi'$ and 
let $Y$ be the affine variety corresponding to $I$. $\qed$ 

\begin{defn}
We will say that a complete local ring $R$ is {\em \qh } if it 
satisfies (a) and (b) as in Lemma \ref{crit-qh}. We will say a germ 
$(Y,o)$ is \qh if the complete local ring $\widehat{O_{Y,o}}$ is \qh. 
\end{defn}

\medskip
Here are several  examples. The polynomial $f(x, y, z)= 
x^2 + y^5 + z^3$ 
is quasi-homogeneous with the weights of generators $15$, 
$6$ and $10$ respectively. The weight of the relation is $30$. 
 Let $g(x)= x^n$, then 
$g$ is quasi-homogeneous for any weight $w$ of the  generator 
and the weight $nw$ of the relation. 

\medskip
Another example is the  germ $(V_p, 0)= (\{x^p=0\}, 0)$, $p\ge 2$ is 
prime. Let's prove that this germ is not quasi-homogeneous for any 
weigths of relations $< p$. Indeed, suppose $Y$ is a quasi-homogeneous 
cone whose germ at zero is isomorphic to $(V_p, 0)$. Since we assume 
that $Y$ is minimal, hence $Y\subset \k$. Polynomials defining 
$Y$ must be monomials (since $Y$ is quasi-homogeneous). Then $Y$ 
clearly  can be defined by a single monomial equation $x^m=0$. 
Isomorphism of germs $(Y,0)\to (V_p, 0)$ induces isomorphisms of 
finite order tangent spaces, hence $m=p$. 

\medskip
In a certain sense {\em generic} germs are not 
quasi-homogeneous. This is discussed in details in 
\cite{Arnold1}, \cite{Arnold2}. Here is one explanation, 
in the case of germs in the affine plane $\A^2$, one which doesn't require 
knowledge of singularity theory but is based on 3-dimensional 
topology.  
Suppose that $Y\subset \A^2$ is a minimal affine curve (defined over $\C$),  
which is invariant under weighted action of $\C^{\times}$ on $\A^2$ 
(with the weights $w_1$, $w_2$). Then the set of complex points 
$$
Y(\C) \subset \C^2
$$ 
is invariant under the $\C^{\times}$-action with the weights 
$w_1$, $w_2$. The corresponding weighted  action of $S^1$ 
preserves a small sphere $S^3$ around  
zero and the link $Y_{\C} \cap S^3= L$. Thus $S^3- L$ admits a 
free $S^1$-action, therefore $S^3- L$ is a Seifert manifold. 
{\em Generic} singularities do not have such property, see 
\cite{Neu}.  For convenience of the reader we describe a way 
to produce examples of singularities which do not admit 
$\C^{\times}$-action. Our discussion follows \cite{Neu}. 
Start with a finite ``Piuseaux series'' 
$$
y= x^{q_1/p_1} (a_1 + x^{q_2/p_1p_2}( a_2 + x^{q_3/p_1p_2 p_3}
(... (a_{s-1} + a_s x^{q_s/p_1... p_s})\ldots )
$$
where  $(p_i, q_i)$ are pairs of positive coprime integers. 
The numbers $a_j$ are nonzero integers. Then $y, x$ satisfy 
some polynomial equation $f(x, y)=0$ with integer 
coefficients, the link $L$ of the singularity at zero 
is an {\em iterated torus knot}, the number $s$ is 
the depth of iteration, numbers $p_i, q_i$ describe 
cabling invariants. The complement $S^3- L$ is not a 
Seifert manifold provided that $s \ge 2$. The simplest 
example is when $s=2$, 
$$
y= x^{q_1/p_1} (a_1 + a_2 x^{q_2/p_1p_2}) 
$$ 
For instance take $a_1= a_2= 1$, $p_1= p_2= 2$, $q_1= q_2=3$ 
(the iterated trefoil knot),  then
$$
y^{2}= x^3 + x^9 + 2 x^6 \quad.
$$
Another example of a singularity which is not \qh is
$$
x^2 y^2 + x^5 + y^5 = 0
$$
see \cite[Page 122]{Dimca}. 

\section{Analytic isomorphisms of algebraic varieties}

In this section we discuss the following question:

Let $\k$ be a field of zero characteristic. 
Suppose that $X'\subset X, Y'\subset Y$ are subvarieties 
in smooth quasi-projective varieties $X, Y$ over $\k$, 
$\eta: X\to Y$ is a biregular isomorphism which carries 
$X'$ bijectively to $Y'$. Does $\eta$ induce a biregular  
isomorphism $X'\to Y'$ ?

Clearly the answer is ``yes'' if both subvarieties $X', Y'$ 
are reduced. The simple example:
$$
X= Y= \C, \quad X'= \{z^2=0\}, \quad Y'= \{z^3= 0\} ,
\quad \eta= id
$$
shows that in the nonreduced case we need some extra 
assumptions to get the positive answer. Our goal is to 
prove that the answer is again positive if we assume 
that $\eta$ induces an {\em analytic} isomorphism 
between $X', Y'$ (Theorem \ref{cbasic}).  

\bigskip
Let $R$ be a ring, $m$ is a maximal ideal in $R$, $R_m$ 
is the localization of $R$ at $m$ and $\hat{R_m}$ is 
the completion of $R$ at $m$. 

\begin{lem}
\label{basic}
Suppose $R$ is a Noetherian ring and $f\in R$ has the 
property that its image in $\widehat{R_m}$ is zero for all 
maximal ideals $m$. Then $f= 0$. 
\end{lem}
\proof By Krull's theorem \cite[Corollary 10.19]{AM}, 
the induced map $R_m \to 
\widehat{R_m}$ is an injection. Thus the image of $f$ in 
$R_m$ is zero for all maximal ideals $m$. Hence for 
every such $m$ there exists $s\notin m$ with $sf=0$. 
Therefore $Ann(f)$ is contained in no maximal ideal. 
This implies that $Ann(f)= R$ and $f=0$. $\qed$

\begin{lem}
\label{lbas}
Let $\phi: R\to S$ be a ring homomorphism and 
$I\subset R, J\subset S$ be ideals. Suppose that 
for every maximal ideal $m$ in $S$ with $m\supset J$ 
we have
$$
\phi(I)\otimes \widehat{S_m} \subset J\otimes \widehat{S_m}
$$
Then $\phi(I)\subset J$.  
\end{lem}
\proof It suffices to prove this when $\phi(I)$ is 
replaced by an element $f$ of $S$. Thus we assume 
that the image of $f$ in $S_m$ is contained in 
$J\otimes S_m$ for all maximal ideals $m\subset S$. 
We want to conclude that $f\in J$. We further simplify 
the situation by dividing by $J$. We have an exact sequence
$$
0\to J \to S\to S/J \to 0
$$
We use \cite[Proposition 10.15]{AM}  to conclude that
$$
\widehat{S/J_{m/J}} \cong S/J \otimes \widehat{S_m}
$$
and \cite[Proposition 10.14]{AM} to conclude
$$
S/J \otimes \widehat{S_m}\cong \widehat{S_m}/(J\otimes \widehat{S_m})
$$
Replace $f$ by its image in $S/J$. We find that the image 
of $f$ in all completions of $S/J$ at maximal ideals of 
$S/J$ is zero. Then $f$ is zero by the Lemma \ref{basic}. 
$\qed$ 

\begin{thm}
\label{cbasic}
Suppose that $X, Y$ are nonsingular (quasi-) 
projective varieties over $\k$ and $\eta: X\to Y$ is an isomorphism. 
Let $X'\subset X, Y'\subset Y$ are subvarieties so that: 
$\eta':= \eta|_{X'}: X'\to Y'$ is a bijection which is an 
analytic isomorphism.  Then $\eta': X'\to Y'$ is a 
biregular isomorphism. 
\end{thm}
\proof It is enough to check the assertion on  open subsets, 
so we may as well assume that $X, Y$ are affine with coordinate 
rings $S, R$, affine subvarieties $Y', X'$ are given by ideals 
$I\subset R, J\subset S$. Coordinate rings of $Y', X'$ are 
$R/I, S/J$. Let $m$ be a maximal ideal in $R/I$, then there 
is a maximal ideal $M\subset R$ such that $m= M/I$. Thus
$$
\widehat{R/I_m} \cong \frac{\widehat{R_{M}}}{I\otimes  
\widehat{R_{M}} } 
$$  
Let $\phi: R\to S$ be the  isomorphism induced by 
$\eta: X\to Y$. 
Since $\eta'$ is an analytic isomorphism it induces 
isomorphisms of all completions 
$$
\frac{\widehat{R_{M}}}{I\otimes  \widehat{R_{M}} } \lra 
\frac{\widehat{R_{\phi(M)}}}{J\otimes \widehat{R_{M}} } 
$$
Thus the assertion of Theorem follows from Lemma \ref{lbas}.  $\qed$ 

%\medskip
%The above Theorem has the following obvious corollary:

%\begin{cor}
%\label{corbasic}
%Suppose that $X, Y, X', Y', \eta, \eta'$ are as in Theorem 
%\ref{cbasic} and $X, Y, \eta$  are defined over $\R$.  Then $\eta'$ 
%induces a biregular isomorphism of the corresponding real varieties. 
%\end{cor}

\section{Coxeter, Shephard and Artin groups}
\label{groups}

Let $\La$ be a finite graph where two vertices are connected by at most one
edge, there are no loops (i.e. no vertex is connected by an edge to itself) 
and each edge $e$ is assigned an integer 
$\eps(e)\ge 2$. We call $\La$ a {\em labelled} graph, let $\V(\La)$ 
and $\E(\La)$ denote the sets of vertices and edges of $\La$. 
When drawing $\La$ we will omit labels $2$ 
from the edges (since in our examples most of the labels are $2$). 
Given $\La$ we construct two finitely-presented 
 groups corresponding to it. The first group $G_{\La}^c$ is called the
 {\em Coxeter group} with the {\em Coxeter graph} $\La$, the second is the 
{\em Artin group} $G_{\La}^a$.   
The sets of generators for the both groups are $\{ g_v, v\in \V(\La)\}$. 
  Relations in $G_{\La}^c$ are:
$$
g_v^2= 1, v\in \V(\La), \ (g_{v} g_{w})^{\eps(e)}= \1, 
\hbox{~~over all edges~~} e= [v, w]\in \E(\La)
$$
 Relations in $G_{\La}^a$ are:
$$
\underbrace{g_{v} g_{w} g_v g_w ... }_{\eps\hbox{~multiples~}}= 
\underbrace{g_{w} g_{v} g_{w} g_v ... }_{\eps\hbox{~multiples~}}\ , \quad 
\eps= \eps(e) \hbox{,~~over all edges~~} e= [v, w]\in \E(\La)
$$
We let $\eps(v, w)= \eps([v,w])$ if $v, w$ are connected by the edge
$[v, w]$ and $\eps(v, w)= \8$ if $v, w$ are not connected by any edge. 
For instance, if we have an edge $[v, w]$ with the label $4$, then the Artin 
relation is
$$
g_{v} g_{w} g_v g_w = g_{w} g_{v} g_{w} g_v
$$
Note that there is an obvious epimorphism
$G_{\La}^a \longrightarrow G_{\La}^c$. 
We call the groups $G_{\La}^c$ and $G_{\La}^a$ {\em associated} 
with each other. 
Artin groups above appear as generalizations of the Artin braid group. 
 Each Coxeter group $G_{\La}^c$ admits a canonical discrete faithful linear
 representation
$$h: G_{\La}^c\lra GL(n, \R)\subset GL(n, \C)$$
where $n$ is the number of vertices in $\La$. Suppose that the Coxeter
 group $G_{\La}^c$ is finite, then remove from $\C^n$ the collection
 of fixed points of elements of $h( G_{\La}^c- \{\1\})$ and denote the 
resulting 
complement $X_{\La}$. The group $G_{\La}^c$ acts freely on $X_{\La}$ and 
the quotient $X_{\La}/G_{\La}^c$ is a smooth complex quasi-projective variety 
with the fundamental group $G_{\La}^a$, see \cite{Bri} for details. Thus the 
Artin group associated to  a finite Coxeter group is  the fundamental 
group of a smooth complex quasi-projective variety. 

\medskip
The construction of Coxeter and Artin groups can be generalized as follows. 
Suppose that not only edges of $\La$, but also  its vertices $v_j$ 
have labels $\del_j= \del(v_j)\in \{0, 2, 3,...\}$. 
Then take the presentation of the 
Artin group $G_{\La}^a$ and add the relations:
$$
g_{v}^{\del(v)}= \1 , \quad v\in \V(\La)
$$
If $\del(v)=2$ for all vertices $v$ then we get the Coxeter group, in 
general the resulting group is called the {\em Shephard group}, they were 
introduced by Shephard in \cite{Sh}. Again there is a canonical epimorphism 
$G^a_{\La}\to G^s_{\La}$. 

Given a Coxeter, Artin or Shephard group $G$ associated with 
the graph $\La$ we define {\em vertex} subgroups 
$G_v, v\in \V(\La)$ and {\em edge} subgroups  
$G_e= G_{vw}, e= [v,w]\in \E(\La)$ as subgroups generated 
by the elements: $g_v$ (in the case of the vertex subgroup) 
and $g_v, g_w$ (in the case of the edge subgroup). 

\medskip
We will use the fact that several of Shephard groups are 
finite. Here is the list of the finite groups that we will 
use (see \cite{Coxeter}):

\begin{figure}[h]
\leavevmode
\centerline{\epsfxsize=3.5in \epsfbox{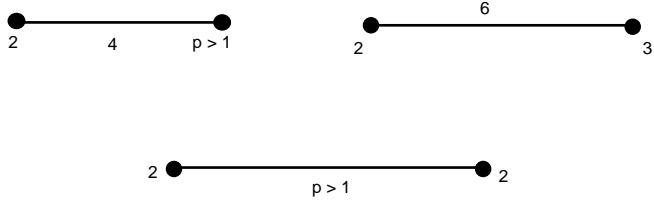}}
\caption{\sl Graphs for certain finite Shephard groups.}
\label{Fig28}
\end{figure}

Some of the Coxeter groups above admit extensions to {\em extended Artin} 
groups, which we describe below. Suppose that $\eps(e)\in \{2, 4,  \8 \}$. 
Enumerate vertices of $\La$. 
Then define a generalized Cartan matrix $N$ as follows. Diagonal 
entries of $N$ are equal to   $2$. Now consider off-diagonal 
entries $n_{ij}$ of $N$ assuming $i< j$. 
If $\eps(v_i, v_j)= 2$ we let $n_{ij}= n_{ji}= 0$. 
If $\eps(v_i, v_j)= 4$ we let $n_{ij}= -1, n_{ji}= -2$. Finally, if
$\eps(v_i, v_j)= \8$ we let $n_{ij}= n_{ji}= -2$. 

\begin{rem}
More generally one is also allowed to consider the case $\eps(e)= 3, 6$.
%,but we will not need this in our paper. 
\end{rem}

Under these conditions there is a {\em generalized root system} associated 
to $N$ and one constructs the {\em extended Artin group} $G_{\La}^e$  
as follows (\cite{Lo}, \cite{Le}):
 
Let $m$ denote the number of vertices in $\La$. 
Let $\tau_j, 1\le j\le m,$ denote generators of $\Z^m$ and 
$g_i, 1\le i\le m,$  denote generators of $G^a$. In the 
free product $G^a * \Z^m$ introduce the following relations: 

(a) Suppose that $-n_{ji}= 2r$, then introduce the relations:  
$$
g_i \tau_j = \tau_j \tau_i^r g_i \tau_i^{-r}
$$ 

(b) Suppose that $-n_{ji}= 2r +1$, then introduce the relations:  
$$
g_i \tau_j = \tau_j \tau_i^{r+1} g_i^{-1} \tau_i^{-r}
$$ 
Note that for edges $[v_i, v_j]$ labelled $2$ we get $n_{ij}= n_{ji}=0$ and 
$$
[g_i, \tau_j] =\1, \quad [g_j, \tau_i] =\1
$$
The resulting group is $G^e_{\La}$.

The extended Artin groups arise naturally as fundamental groups of 
certain Stein manifolds $M(G^e)$,  which are  complements of discriminants of 
moduli  spaces of versal deformations of certain 
singularities ({\em simply elliptic} and {\em cusp} singularities), see
 \cite{Lo}, \cite{Le}. 
Below is an outline of the construction. Let $G^c$ be the 
Coxeter group associated with $G^a$, $h: G^c\to GL(n, \R)$ 
be its canonical (discrete and faithful) representation. 
The group $h(G^c)$ acts properly discontinuously on a 
convex cone $C\subset \R^n$, which is the ``maximal'' domain 
of discontinuity of $h(G^c)$. Complexify this action. The group $G^c$ acts 
on $\Z^n\subset \R^n$. Identify the $G^c$-module $\Z^n \otimes \R$ with the   
``imaginary part'' of $\C^n$. Thus we get a properly 
discontinuous action of $G:= \Z^n \ltimes G^c$ on 
$W:= V \times (\Z^n \otimes \R) \subset \C^n$. Take the 
quotient $W/G$ and remove from it projections of fixed-points 
of elements of $G- \{\1\}$. The result is the manifold 
$M(G^e)$. 

\section{Local deformation theory of representations}
\label{defof}
Let $\k$ be a field of zero characteristic, define the Artin local 
$\k$-algebra  $\Pi_m$ as $\k[t]/t^{m+1}$. Take a finitely-generated 
group $\Ga$ and 
an algebraic  Lie group $G$ over $\k$, let $\G:= G(\k)$ be the 
set of $\k$-points 
of $G$. Then the Lie group $G(\Pi_m)$ splits as the semidirect product 
$G_0(\Pi_m)\rtimes \G$, where $G_0(\Pi_m)$ is the kernel of the natural 
projection $p_m: G(\Pi_m)\to \G$ induced by $\Pi_m \to \k$, which is given 
by $t\mapsto 0$. The group $G_0(\Pi_m)$ is $\k$-unipotent. We consider the 
germ
$$
(Hom(\Ga, \G), \rho)
$$
for a certain representation $\rho: \Ga \to \G$. Then $m$-th order 
Zariski tangent space \newline 
$T^m_{\rho}( \Hom(\Ga, \G))$ is naturally isomorphic 
to the space of homomorphisms
$$
Z^1_{(m)}(\Ga, \G; \rho):= 
\{ \xi: \Ga \to G(\Pi_m) | p_m \circ\xi = \rho\}
$$
(see \cite{GM}). If $m=1$ then $Z^1_{(m)}(\Ga, \G; \rho)\cong 
Z^1(\Ga, ad\circ\rho)$, where $ad$ is the adjoint action of $\G$ on 
its Lie algebra.  

We call $\xi\in T^m_{\rho}( \Hom(\Ga, \G))$ a {\em trivial 
deformation} of $\rho$ if there is an element 
$h= h(\xi)\in G(\Pi_m)$ such that 
$\xi(g)= h \rho(g) h^{-1}$ for all $g\in \Ga$. If $m=1$ then 
the infinitesimal deformation $\xi$ is trivial iff the 
corresponding cocycle is a coboundary. 

Suppose that $G^s$ is a Shephard group, $G^a$ is the corresponding 
Artin group, $q : G^a \to G^s$ is the canonical projection. 
Let $V_L$ denote the set of vertices with nonzero 
labels  in the graph of $G^s$. Let $\G$ be a group of $\k$-points 
of an algebraic Lie group.  Consider 
a homomorphism $\rho: G^s\to \G$, let $\t\rho= \rho \circ q$. The 
projection  $q$ induces an injective
 morphism of the representation varieties
$$
q^* : \Hom(G^s, \G)\lra \Hom(G^a, \G)
$$
and injective morphisms of the corresponding finite order Zariski 
tangent spaces
$$
D^{(n)} q^*: T^n_{\rho} \Hom(G^s, \G)\lra T^n_{\t\rho} \Hom(G^a, \G)
$$

\begin{lem}
\label{L7.1}
Suppose that $\xi\in T^n_{\t\rho} \Hom(G^a, \G)$ is an element whose 
restriction to each cyclic vertex subgroup $G_v^a, v\in V_L$ 
is a {\bf trivial deformation} of $\t\rho|_{G_v^a}$. Then $\xi$ 
belongs to the image of  $D^{(n)} q^*$. 
\end{lem}
\proof Recall that we identify $T^n_{\t\rho} \Hom(G^a, \G)$ and 
$T^n_{\rho} \Hom(G^s, \G)$ with spaces of certain representations of 
$G^a, G^s$ into $G(\Pi_n)$. We have the exact 
sequence
$$
\1\lra \<\< \{ g_v^{\del(v)}, v\in V_L\}\>\> \lra G^a \lra G^s \lra \1 
$$
Then $\xi$ belongs 
to the image of $D^{(n)} q^*$ if and only if $\xi(g_v^{\del(v)})=\1$ for 
all $v\in V_L$. We assume that the restriction of $\xi$ to $G_v^a$ is
a trivial deformation, thus there is an element $h:= h_{v}\in \G(\Pi_n)$ 
such that $\xi(g_v)= h \t\rho(g_v) h^{-1}$. Since $\t\rho(g_v^{\del(v)})=\1$ 
we conclude that  $\xi(g_v^{\del(v)})=\1$ as well. $\qed$

\begin{cor}
\label{CC}
Suppose that  $H^1(G^s, ad \circ\rho)=0$ (i.e. the representation $\rho$ is 
infinitesimally rigid) and for each $v\in V_L$ the restriction 
homomorphism 
$$
Res_v: \H^ 1(G^a, ad\circ\rho)\to \H^ 1(G_v^a, ad\circ \rho)
$$ 
is zero. Then 
$$
q^* : (\Hom(G^s, \G), \rho) \lra (\Hom(G^a, \G), \t\rho)
$$
is an analytic isomorphism of germs. 
\end{cor}
\proof Let $Z(\rho)$ denote the centralizer of $\rho(G^s)$ if $\G$. Then 
the representation variety $\Hom(G^s, \G)$ is smooth near $\rho$ and is 
naturally isomorphic to the quotient $\G/Z(\rho)$, see \cite[Theorem 2.6]{LM}. Now the assertion follows from Lemmas \ref{LL3.1}, \ref{L7.1}.  $\qed$ 

One example when the first condition of the corollary are satisfied is 
the case when the Shephard group $G^s$ is finite. 

\section{Projective reflections}
\label{S3}

Fix the bilinear form $\flat= x_1 y_1 + x_2 y_2 + x_3 y_3= 
(x_1, x_2, x_3) \cdot (y_1, y_2, y_3)$ on the vector 
space $V= \C^3$, we shall also use the notation 
$\< \cdot , \cdot \>$ for $\flat$.  
Let $\varphi$ denote  the quadratic form corresponding to $\flat$, 
let $O(3,\C)$ be the group of automorphisms 
of $\flat$. Let $\pi: V \to \P(V)$ denote the 
quotient map and let 
$CO(3, \C):= \C^{\times}\cdot O(3,\C) \subset End(V)$. The inclusion 
$CO(3,\C) \to End(V)$ induces an embedding $PO(3,\C) \hook \P(End~V)$. Let 
$p: End(V)  \to \P(End~V)$ be the quotient map. Note that 
$PSL(2, \C)\cong PO(3, \C)$ and $SO(3,\R)\cong PO(3, \R)$. 

\subsection{The correspondence between  projective reflections and their  
fixed points}

In this section we study projective properties of elements of order $2$ in 
the group $PO(3, \C)$. Consider an element $A\in O(3, \C)$ such that the 
projectivization $p(A)$ is an involution acting on $\P^2(\C)$. The 
fixed-point set of $p(A)$ consists of two components: an isolated point 
$a$ and a projective line  $l$ dual to $a$ (with respect to $\flat$).  
Our goal is to describe the 
correspondence $p(A)\leftrightarrow a$ in algebraic terms. 

Let $R\subset SO(3,\C)$ be the affine subvariety of involutions. Note that 
$-\1\in O(3,\C)$ doesn't belong to $R$. 
We leave the proof of the following lemma to the reader. 

\begin{lem}
$PO(3, \C)$ acts transitively by conjugations 
on $PR(\C)$ (the image of $R(\C)$ in $PO(3, \C)$). 
\end{lem}

We now determine $\overline{p(R)}$, the Zariski closure of $p(R)$ in 
$\P(End~V)$. We define a morphism $\eta: V \to End(V)$ by $\eta(v)(x)= 
\varphi(v)x - 2\<v, x\> v$. 
If $v$ is an anisotropic vector then $\eta(v)$ is a multiple 
of the reflection through the hyperplane in $V$ orthogonal to $v$. 
The reader will verify that $\eta$ is an $O(3)$-equivariant morphism, i.e. 
$\eta(g v)= g \eta(v) g^{-1}$, $g\in O(3), v\in V$, and that $\eta$ induces 
an equivariant embedding  
$$
\eta: \P(V) \lra \P(End~V)
$$  
which is an immersion of one smooth complex manifold into another. 

\begin{lem}
The image of $\eta$ is contained in $\overline{p(R)}$. 
\end{lem}
\proof Let $V_0$ be the complement in $V$ of $X= \{ v \in V: \varphi(v)=0\}$.  
Then $\P_0(V):= \pi(V_0)$ is Zariski dense in $\P(V)$. But also 
$\eta(V_0)\subset \C^{\times} R$. Hence 
$\overline{\P_0(V)} \subset \overline{p(R)}$. $\qed$ 

We now consider $\P(V)$ and $\P(End(V)$ as varieties over $\Q$. 

\begin{lem}
\label{fixed point}
The morphism $\eta$ induces an isomorphism of varieties $\P(V) \cong \overline{p(R)}$. 
\end{lem} 
\proof Since $\eta: \pi(V_0)\to p(R)$ is equivariant it is easy to verify 
that it is onto (in fact a bijection). Hence $\eta(\pi(V_0))= p(R)$ and 
accordingly 
$\overline{\eta\pi(V_0)} = \overline{p(R)}= \eta(\P(V))$. But we have seen 
that the morphism $\eta: \P(V)\to \eta(\P(V))$ is an analytic isomorphism of 
smooth compact complex manifolds. Hence (by the GAGA-principle) 
it is an isomorphism of projective varieties. $\qed$ 

\medskip
Let $N:= p^{-1} (\overline{p(R)})- p^{-1} (p(R))$. Then $N$ may be 
described as follows. The bilinear form $\flat$ induces an 
isomorphism $\t\flat: V\otimes V \to V^* \otimes V = End(V)$. Then 
$N= \t\flat(X \otimes X)$. We note that $n_v := \t\flat(v \otimes v)$ 
is then given by $n_v(x)= \<v, x\> v$. Hence the set of real points 
$N(\R)$ is empty. 

Let $Q\subset End(V)$ be the affine cone defined by 
$Q:= p^{-1} (\overline{p(R)})$. Hence $N \subset Q$. We define 
$Q_0:= Q- N$. Then $\eta$ induces a commutative diagram
$$
\begin{array}{ccc}
V & \longrightarrow & Q  \\
\Big\uparrow & ~ & \Big\uparrow\\
V_0 & \longrightarrow & Q_0\\
\end{array}
$$

\begin{rem}
It can be shown that the cone $Q\subset End(V)$ is defined by the equations:
\begin{enumerate}
\item $X X^{\top} = X^{\top} X$, 
\item $X X^{\top} E_{ij} = E_{ij}  X X^{\top}$, $1\le i, j\le 3$,
\item $X^2 E_{ij} = E_{ij} X^2$, $1\le i, j\le 3$.
\end{enumerate}
Here $E_{ij}$ is the matrix with $1$ in the $(ij)$-th position and $0$ 
elsewhere. The equations (1) and (2) define the closure 
$\overline{PO(3)} \subset \P(End~V)$. We will not need the explicit 
equations for $\overline{p(R)}$ in what follows. 
\end{rem}

We let $\zeta: PQ \to \P^2$ be the inverse of $\eta$ and abbreviate 
$\P(V_0)$ to $\P^2_0$. Note that \newline 
$\zeta: PR = PQ_0 \to \P^2_0$\  assigns 
to each projective reflection its neutral (isolated) 
fixed-point. Thus we have described the correspondence 
$p(A) \leftrightarrow a$ algebraically. Note that we have 
$\P^2(\R)= \P^2_0(\R)$ and $PQ(\R)= PQ_0(\R)$. Let   
$\hat\flat: V \to V^*$ be the isomorphism induced by $\flat$. 
Define $(\P^2_0)^{\vee}$ by 
$(\P^2_0)^{\vee}= \hat\flat(\P^2_0)$. 
Hence the space of {\em anisotropic lines} 
$(\P^2_0)^{\vee}$ is the space of lines dual to the set of 
{\em anisotropic points}  $\P^2_0$.

\subsection{Fixed points}

Suppose that $g\in PO(3, \C)$ is a nontrivial element. 
In this section we discuss the fixed-point set for the 
action of $g$ on $\P^2(\C)$. 

\begin{defn}
A fixed point $x$ for the action of $g$ on $\P^2(\C)$ is called {\bf neutral} 
if the determinant of the differential of $g$ at $x$ is equal to 1. 
\end{defn}

There are two classes of nontrivial elements $g\in PO(3, \C)$: 
(a) $g$ is unipotent, (b) $g$ is semi-simple.

\medskip
{\bf Case (a).} In this case $g$ has a single fixed 
point $a\in \P^2$, the point $a$ belongs to an 
invariant projective line $L\subset \P^2$. If we 
choose coordinates on $L$ such that $a= \8$, then 
$g$ acts on $L- \{a\}$ as a translation.  The flag 
$(L, a)$ is determined by the element $g$ uniquely. 
On the other hand, the flag $(L, a)$ uniquely 
determines the 1-parameter maximal unipotent 
subgroup in $PO(3)$ which contains $g$. Finally, 
the fixed point $a$ of $g$  uniquely determines 
the line $L$. It is easy to see that $a$ is the {\em neutral} 
fixed point of $g$.  

The collection $P\subset \P^2$ of fixed points of all 
unipotent elements in $PO(3,\C)$ is the projectivization of the 
cone $\{ \ov{x}\in \C^3 : \varphi(\ov{x})= 0\}$. 

\medskip
{\bf Case (b).} Each semisimple element of $O(3)$ 
is conjugate (in $GL(3, \C)$) to 
$$
A= \left( \matrix  \pm \la & 0 & 0\\
0 & \la^{-1} & 0\\
0 & 0 & 1\endmatrix \right).
$$
(up to $\pm$). We have two possible cases depending on whether 
or not $A^2=\1$. If $A^2\ne \1$ then $A$ has three distinct 
complex eigenvectors, one of them $\ov{f}$  
 is fixed by $A$. Thus, in the case the 
transformation $p(A)$ has three fixed points on $\P^2(\C)$, 
one of them $f= p(\ov{f})$ is the {\em neutral} fixed point. 
The maximal torus in $PO(3, \C)$ 
containing $p(A)$ is uniquely determined by the (neutral) 
fixed point $f$ in the both real and complex cases. 

Finally, consider the case $A^2= \1$. Then  $p(A)$ has one 
isolated fixed point $f$ on $\P^2(\C)$ (which is  
{\em neutral}) and the fixed projective line disjoint from $f$.

\medskip
We will realize the Lie algebra $so(3, \C)$ of $PO(3, \C)$ 
as the subalgebra in the space of vector-fields on $\P^2$. 
Again, each nonzero element $\zeta\in so(3,\C)$ has either one 
or $3$ zeroes on $\P^2(\C)$. One of them is the {\em neutral zero}, 
it corresponds to the neutral fixed point of the  subgroup  
$\exp(\C \cdot \zeta)\subset PO(3,\C)$. The maximal abelian subalgebra 
$\C \cdot \zeta$ of $so(3,\C)$ containing $\zeta$  is determined 
by the neutral zero of $\zeta$.  

\subsection{Commuting and anticommuting elements}
\label{Commuting}
 
Suppose that $\al, \al'$ are  involutions in $PO(3,\C)$, thus 
they have isolated fixed points $\la, \la'$ and fixed 
projective lines $\La, \La'\subset \P^2(\C)$. 

\begin{lem}
The elements $\al, \al'$ commute if and only if either:

 (1) $\al= \al'$ and $\la= \la', \La= \La'$, or (2) $\la\in \La'$, 
$\la' \in \La$ and $\La$ intersects $\La'$ orthogonally 
(with respect to the quadratic form $\varphi$). 
\end{lem}  
\proof Proof is obvious and is left to the reader. $\qed$

\medskip
Suppose that $\al$ is an involution in $PO(3,\C)$ and $\be\in 
PO(3,\C)- \{\1\}$.  
We say that the elements $\al, \be$ {\em anticommute} if 
$\al \be \al \be = \1$ (i.e. $\al \be \al^{-1}=  \be ^{- 1}$) and 
$\al\ne \be$. 

\begin{lem}
\label{anticom}
The elements $\al, \be$ anticommute iff the neutral fixed point 
$\mu$ of $\be$  belongs to the fixed projective line $\La$ of 
$\al$. 
\end{lem}  
\proof If $\al, \be$ anticommute then $\mu$ must belong to the 
fixed-point set of $\al$. If $\mu$ is the neutral point of $\al$ 
then $\al$ and $\be$ commute (and $\be^2\ne \1$) 
or $\al= \be$ which contradicts our assumptions.  
$\qed$ 

\begin{rem}
Suppose that $\al, \be$ are anticommuting elements as above. 
Then they satisfy the relation
$$
\al \be \al \be = \be \al \be \al
$$
(since both right- and left-hand side are equal to $\1$). 
\end{rem}

\section{Cohomology computations}
\label{elem}

In this section we consider cohomologies of certain {\em elementary} 
Artin and Shephard groups (their graphs have only two vertices). The section 
is rather technical, its material will be needed in Section \ref{S14}. 

Let $G$ be a finitely-generated group, $\rho: G\to PO(3, \C)$ be  
a representation. We will denote the action of elements $g$ of $G$ 
on vectors $\xi$ in the Lie algebra $so(3,\C)$ by 
$$
g\xi:= ad\circ\rho(g) \xi$$

\subsection{2-generated Abelian group}

\begin{figure}[tbh]
\leavevmode
\centerline{\epsfxsize=3.5in \epsfbox{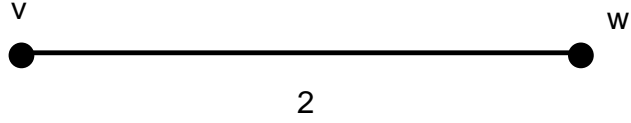}}
\caption{\sl Graph $\Ga$ for 2-generated free Abelian group.}
\label{Fig21}
\end{figure}

Let $\Ga$ be the graph with two vertices $v, w$ connected by the edge 
$e$ with the label $2$ (Figure \ref{Fig21}). 
The corresponding Artin group $G^a$ is free 
Abelian group on two generators $a= g_v, b= g_w$. Take two distinct   
commuting involutions $\al, \be \in PO(3)$ and the homomorphism 
$\rho_0: G^a\to PO(3,\C)$ given 
by $\rho_0: a\mapsto \al$, $\rho_0: b\mapsto \be$.  

\begin{lem}
\label{L5.3}
The representation $\rho_0$ is infinitesimally rigid, i.e. 
$\H^ 1(G^a, ad\circ \rho_0)= 0$. The point 
$\rho_0\in \Hom(G^a, PO(3,\C))$ is nonsingular. 
\end{lem}
\proof $G^a= \pi_1(T^2)$ (where $T^2$ is the 2-dimensional torus) 
and the group $PO(3)$ is reductive, hence Poincare duality gives 
an isomorphism 
 $\H^ 0(G^a, ad \circ \rho_0)\cong \H^ 2( G^a, ad \circ \rho_0 )$. 
The centralizer of $\rho_0( G^a)$ in $PO(3,\C)$ is trivial. 
Therefore, 
$$
0= \H^ 0(G^a, ad \circ \rho_0)\cong \H^ 2( G^a, ad \circ \rho_0 )
$$
Since the Euler characteristic of $T^2$ (and hence of $G^a$) 
equals zero, we conclude that $0= \H^ 1(G^a, ad \circ \rho_0 )$. 
The second assertion of Lemma follows from \cite[Theorem 2.6]{LM}. 
$\qed$

\medskip
Note that the associated Coxeter group $G^c$ is the 
finite group $\Z/2 \times \Z/2$. Let $\rho^c_0$ 
denote the homomorphism $G^c\to PO(3,\C)$ corresponding to $\rho_0$. 

\begin{cor}
\label{C1}
We have natural isomorphism between germs of representation
varieties:
$$
(\Hom(G^c, PO(3,\C)), \rho_0^c)\cong (\Hom(G^a, PO(3,\C)), \rho_0)
$$
given by composing homomorphisms $G^c\to PO(3,\C)$ with 
the projection $G^a \to G^c$. 
\end{cor}
\proof Follows from the Lemma \ref{L5.3} and Corollary 
\ref{CC}. $\qed$ 

\medskip
Now we consider the global structure of $\Hom(G^a, PO(3))$. 

\begin{lem}
\label{L5.8}
The variety $\Hom(G^a, PO(3, \C))$ is the disjoint union 
of two Zariski closed subsets:

(1) Reducible representations 
$S_1:= \{\rho : \dim \H^ 0(G^a, ad\circ \rho)\ge 1\}$;

(2) The orbit  $S_2:= Ad(PO(3,\C)) \rho_0$, where $\rho_0$ 
is the  representation  whose image is the group 
$\Z/2 \times \Z/2$.  
\end{lem}
\proof First we verify that $\Hom(G^a, PO(3, \C))$ is the 
union $S_1\cup S_2$. Let $A\subset PO(3,\C)$ be an abelian 
subgroup, $\bar{A}$ is the Zariski
closure of $A$. If $A$ is infinite then the abelian group 
$\bar{A}$ is the finite extension of a 1-dimensional 
connected abelian Lie 
subgroup $\bar{A}^0$ of $PO(3,\C)$. Hence $\bar{A}^0$ 
is either a maximal torus 
or a maximal unipotent subgroup of $PO(3,\C)$: both are 
maximal abelian subgroups in $PO(3,\C)$, thus 
$A\subset \bar{A}^0$. We apply this to the 
group $A= \rho(G^a)$ and conclude that in this case
  $\dim \H^ 0(G^a, ad\circ \rho)=1$. 

Now we consider the case when $A$ is finite, hence, after 
conjugation, we get: $A\subset SO(3,\R)$. It follows from 
classification of finite subgroups in $SO(3,\R)$ that $A$ 
either has an invariant vector in $\R^3$ (which again means 
that $\dim \H^ 0(G^a, ad\circ \rho)\ge 1$ ) or 
$A\cong \Z/2 \times \Z/2$ is generated by two involutions 
with orthogonal axes. Hence $S_1\cup S_2=\Hom(G^a, PO(3, \C))$.  

\medskip
It's clear that $S_1$ and $S_2$ are disjoint. 
Since the representation $\rho_0$ is locally rigid it 
follows that the orbit 
$Ad(PO(3,\C) )\rho_0$ is open, \cite[Theorem 2.6]{LM}. 
The representation $\rho_0$ is {\em stable} (see \cite[Theorem 1.1]{JM}) 
since $\rho_0(G^a)= \Z/2 \times \Z/2$ is 
not contained in a proper parabolic subgroup of $PO(3, \C)$. 
Hence the orbit  $Ad(PO(3,\C) )\rho_0$ is closed. $\qed$

\subsection{Finite elementary Shephard groups}

Take $\Ga$ be any of the labelled graphs from the Figure 
\ref{Fig28}. The corresponding Shephard group 
$G:= G^s_{\Ga}$ is finite, hence we have  

\begin{prop}
\label{L5.4}
All  representations of $G$ to $PO(3)$ are
infinitesimally rigid and \newline
$\Hom(G, PO(3))$ is smooth. 
\end{prop}

We will need a slight modification of the above proposition. Let 
$L$ be an algebraic Lie group with the Lie algebra ${\cal L}$, 
$G$ is a finitely-generated group such that all representations 
$\rho: G\to L(\C)$ are infinitesimally rigid, pick elements 
$a\in G$ and $\al\in L(\C)$ and consider the subvariety 
$$
F= F_{a, \al}(G, L)= 
\{ \rho: G\to L| \rho(a)= \al \} \subset \Hom(G, L)
$$

\begin{prop}
\label{PF}
The subvariety $F$ is smooth. 
\end{prop}
\proof 
The space $\Hom(G, L(\C))$ is the union of $L(\C)$-orbits 
of representations $\rho_j, 1\le j\le m$. If $\al\ne \rho(a)$ 
for all $\rho\in \Hom(G, L(\C))$ 
then there is nothing to prove. Otherwise we can assume that 
$\al= \rho_j(a)$, $j\in J\subset\{1,..., m\}$.  
Since the representations $\rho_j$ are locally rigid we get
$$
F(\C)\cong \bigcup_{j\in J} Z_{L(\C)}(\<\al\>)/ Z_{L(\C)}(\rho_j(G))
$$
where $Z_{L(\C)}(H)$ denotes the centralizer of a subgroup 
$H\subset L(\C)$. 
It is enough to verify smoothness of $F$ at the representations 
$\rho_j, j\in J$.  Consider the Zariski tangent space 
$T= T_{\rho_j}F(\C)$. It is naturally isomorphic to
$$
\{ \xi\in Z^1(G, {\cal L}_{ad\rho_j}) : \xi(a)=0\}
$$
However infinitesimal rigidity of $\rho_j$ implies 
that\footnote{Recall  that ${\cal L}^H$ denotes the subspace of 
$H$-invariant vectors for a group $H$ acting on ${\cal L}$.} 
$T\cong {\cal L}^{\<a\>}/{\cal L}^G$, 
where $G$ and $\<a\>$ act on ${\cal L}$ via the adjoint 
representation $ad\rho_j$. Hence the dimension of $F(\C)$ 
(as a complex manifold) at $\rho_j$ 
is equal to the  dimension of its Zariski tangent space at 
$\rho_j$, which implies that $F$ is smooth at $\rho_j$. $\qed$ 

As a particular case we let $a$ be one of the generators of 
$G= G^s_{\Ga}$, pick an element $\al\in PO(3,\C)$ and 
consider the subvariety 
$$
F_{a, \al}(G, PO(3))
= \{ \rho: G\to PO(3,\C)| \rho(a)= \al \} \subset \Hom(G, PO(3,\C))
$$

\begin{cor}
\label{PFC}
Suppose that $G= G^s_{\Ga}$ is a Shephard group as above, $a$ is 
one of the generators of $G$, $\al\in PO(3, \C)$. Then the 
$F_{a, \al}(G, PO(3))$ is smooth. 
\end{cor}

\subsection{The infinite cyclic group}
\label{cy}

Consider the infinite cyclic group $G= \<b\>$ and a representation 
$\rho: \<b\>\to PO(3,\C)$ so that $\rho(b)=\be$ is a nontrivial  
semisimple element. 

Note that the Lie algebra $so(3,\C)$ has $ad\circ \rho(b)$-invariant splitting $L\oplus L^{\perp}$ where $L$ consists of vectors vanishing at the neutral 
fixed point point $B\in \P^2(\C)$ of the element $\be$. The action 
of $b$ on $L^{\perp}$ has no nonzero invariant subspaces. Thus 
$\H^ 1(G, L^{\perp})= 0$ and 
$$
\H^ 1(G, L)\oplus \H^ 1(G, L^{\perp})  = \H^ 1(G, L)\cong  \H^ 1(G, so(3,\C))
$$

This proves the following

\begin{prop}
\label{cycl}
Any cocycle $\si \in Z^1(G, ad\circ \rho)$ has the form $\si(b)= \tau + 
b\xi - \xi$ where $b\tau= \tau$. The element $\tau$ depends only on the cohomology class of $\si$. 
\end{prop}

\subsection{Elementary Shephard groups with the edge-label $4$.}

\begin{figure}[tbh]
\leavevmode
\centerline{\epsfxsize=3.5in \epsfbox{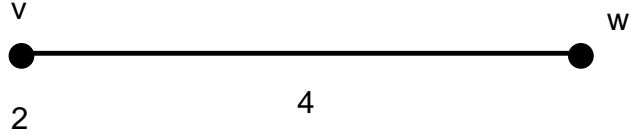}}
\caption{\sl Graph for a Shephard group}
\label{Fig23}
\end{figure}

Consider the graph with two vertices $v, w$ connected by the edge with 
the label $4$, we put the label $2$ on the vertex $v$, see 
Figure \ref{Fig23}. 
The corresponding Shephard group $G$ has the presentation
$\< a= g_v , b= g_w | a^2= \1, (ab)^2= (ba)^2\>$.

 Consider a representation $\rho: G \to PO(3, \C), \  \rho: a\mapsto 
\al, \  \rho: b \mapsto \be$\ , 
where we choose $\al$ to be an involution, $\be\ne \1$ is an  
element such that $\al,\be$ anticommute, i.e. $\al \be \al \be =\1$ and 
$\al \ne \be$. We define another group $\bar{G}$ (which is 
 a Coxeter group) by
$$
\bar{G}= \< a, b|  a^2=\1, (ab)^2=\1 \>\cong \Z/2 * \Z/2
$$
Let $q: G\to \bar{G}$ denote the natural projection. 
It's clear that we have a homomorphism 
$\bar{\rho}: \bar{G}\to PO(3, \C)$ such 
that $\bar\rho\circ q=\rho$. Thus we have the induced monomorphism 
$$
q^*: \H^ 1(\bar G, ad\circ \bar\rho)\lra \H^ 1(G, ad\circ \rho)
$$

\begin{prop}
\label{L5.5}
The mapping $q^*$ is an isomorphism. Any cocycle 
is cohomologous to $\si\in Z^1(G, ad\circ\rho)$ of the form
$$
\si(a)= 0, \quad \si(b)= \tau, \hbox{~~so that~~} b\tau= \tau
$$
If $\si|_{\< b\>}$ is a coboundary then $\si=0$. 
\end{prop}
\proof Let $\si\in Z^1(G, ad\circ\rho)$,
 by normalizing $\si$ via coboundary we can assume that $\si(a)= 0$,
 hence there is an element
$\xi\in so(3, \C)$ such that $\si(b)= b\xi - \xi + \tau$, where 
$b\tau= \tau$ (see Proposition \ref{cycl}).  Now we need 

\begin{lem}
The cocycle $\si$ can be chosen within its cohomology class in 
$\H^ 1(G, ad\circ\rho)$ so that $\si(a)=0, a\xi= - \xi$. 
\end{lem}
\proof Let $\theta:= (\xi + a \xi)/2$, then $a\theta= \theta$. Let 
$\xi':= \xi - \theta$, let $\si'(a):= 0, \si'(b):= b\xi' - \xi' + \tau$. 
Then its clear that $\si'$ is cohomologous to $\si$ and a direct 
computation shows that  $a\xi'= -\xi'$. $\qed$ 

Now we suppose that $\si$ is chosen as in the above lemma. 
Since $\al$ anticommutes with $\be$ it also anticommutes with all elements of 
the 1-parameter subgroup $B\subset PO(3, \C)$ containing $\be$. Suppose for 
a moment that $\tau \ne 0$, since $b\tau = \tau$ we conclude that 
$B= \exp(\C \cdot \tau)$; hence $a\tau = -\tau$. The same of course is true if  
$\tau=0$. Thus $a\tau= -\tau$ in any case. The Artin relation in $G$
is $ab ab  = ba ba $. Therefore we must have:
$$
\si(ab ab )= \si( ba ba)
$$
Direct calculation shows that
$$
\si(ab ab)= - 2 b^{-1} \xi +2\xi ~~, \quad \si(ba ba)= 2b\xi - 2\xi
$$
Thus we must have: $b\xi + b^{-1}\xi - 2\xi= 0$, equivalently 
$b\xi - \xi= \xi- b^{-1}\xi$. Let $\eta:= b\xi - \xi$, then the above 
equation implies: 
$$
b^{-1} \eta= \eta , \quad b  \eta= \eta
$$
On the other hand $a \eta= \eta$ because 
$$
a \eta = ab \xi - a \xi = b^{-1} a \xi + \xi = 
\xi- b^{-1}\xi = b\xi - \xi= 
\eta $$
Since $\al, \be$ anticommute we conclude that $\eta=0$, which implies that 
$b \xi= \xi$.
 Thus $\si(b)= \tau$. Let's check that the cocycle $\si$ lies in 
the image of the mapping $q^*$, for this it's enough to verify that 
$\si(ab ab)=0$. Indeed, $\si(ab ab)= -\tau + \tau=0$. 
It is left to prove that the assmption that $\si$ is a coboundary 
on $\<b\>$ implies that $\tau=0$. 
Let $L\subset so(3, \C)$ be the subspace fixed by the adjoint action of 
$\be$ (as in Section \ref{cy}). By Proposition \ref{cycl}, 
if cocycle $\si$ is a coboundary on $\<b\>$ then there exists 
$\eta\in L$ such that $\si(b)= \tau= b\eta - \eta= \eta -\eta=0$. $\qed$ 

\begin{rem}
\label{smo}
Note that the cocycles $\si$ described in the above Lemma are 
integrable, they
correspond to deformations $\rho_t$  of the representation $\rho$ 
which fix $\rho(a)= \al$ and 
change the element $\rho(b)= \be$ within the corresponding 
1-parameter subgroup 
in $PO(3, \C)$. For such  representations the elements 
$\rho(a), \rho_t(b)$ anticommute. 
\end{rem}

\subsection{Elementary  Artin group with the edge-label $6$.}

\begin{figure}[tbh]
\leavevmode
\centerline{\epsfxsize=3.5in \epsfbox{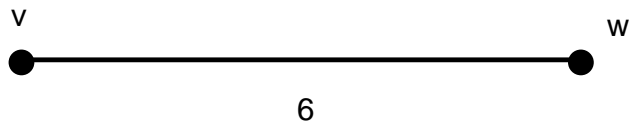}}
\caption{\sl Graph for the Artin group $G^a$}
\label{Fig25}
\end{figure}

Now consider 2-generated Artin group $G^a$ given by  the relation 
$(ab)^3= (ba)^3 $, see Figure \ref{Fig25}. 
Take a representation $\rho: G^a\to PO(3,\C)$ which maps 
$a$ and $b$ to elements $\al, \be$ so that $\al^2=1$, 
$[\al ,\be^2] \ne \1$ and the product $\al\be$ has the order $3$. 
(In particular $\be \ne \1$.) 

\begin{lem}
\label{L5.6}
Let $\si\in Z^1(G, ad\circ\rho)$ be a cocycle  such that
$\si(a)=0, \si(b)= \tau$, where $b\tau= \tau$. Then $\tau= 0$. 
\end{lem}
\proof Direct calculation shows that
$$
\si(ab ab ab)= a \tau + \tau + aba \tau= a\tau + \tau + b^{-1} a\tau$$
$$
 \si(ba ba ba)= \tau + ba \tau + a\tau
$$
Therefore $ad(\rho(b^{-1} a))\tau = ad(\rho( ba)) \tau$. 
We suppose that $\tau\ne 0$, then it has unique zero in $\P^2(\C)$  
which is the same as the neutral fixed point of $\be$. This implies that 
$$
\be^{-1}\al \ \ga\ (\be^{-1}\al )^{-1} = \be \al \ \ga\ (\be \al)^{-1}  
$$
for all elements $\ga$ of the 1-parameter subgroup 
$\exp(\C \cdot \tau)$, in 
particular it is true for $\be= \ga$. Since $\al^2=1$, 
$(\al\be)^3= (\be \al)^3 =1$ we conclude that 
$\be^{-2} \al = \al \be^{-2}$, 
which contradicts our restrictions on $\al, \be$. $\qed$ 

\subsection{A non-elementary Shephard group.}

\begin{figure}[tbh]
\leavevmode
\centerline{\epsfxsize=3.5in \epsfbox{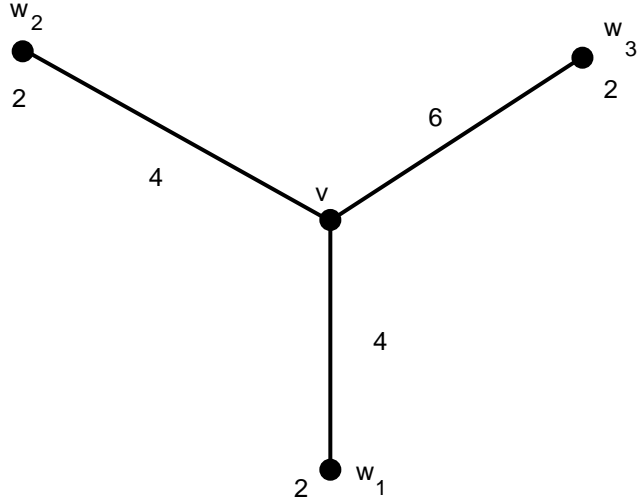}}
\caption{\sl Graph for the nonelementary Shephard group $G^s$}
\label{Fig26}
\end{figure}

Now suppose that we have a group $G^s$ with the presentation:
$$
\< a_1, a_2, a_3,   b | a_j^2= \1, \quad j=1, 2, 3\ ; \quad 
(a_i b)^2=  (b a_i)^2 , i=1,2 ; ~~~ (a_3 b)^3= (b a_3)^3\>
$$
See the graph on the Figure \ref{Fig26}.

Consider a representation $\rho_0: G^s \to PO(3, \C)$ so that

\begin{enumerate}
\item $\rho_0(a_j)\ne \1$, $1\le j \le 3$;  
\item $\rho_0(b)^3= \1$, $\rho_0(b_3)\ne \1$;
\item the neutral fixed points of the elements 
$\rho(a_1), \rho(a_2), \rho(b)$ do not belong to a common projective 
line in $\P^2$; 
\item $\rho_0(a_3 b)^3 =\1$.  
\end{enumerate}

Take the  subvariety $V\subset \Hom(G, PO(3,\C))$ which contains 
$\rho_0$ and consists of homomorphisms that are constant on the 
generators $a_j, 1\le j\le 3$. 

\begin{lem}
\label{L5.10}
The point $\rho_0$ is an isolated reduced point in $V$.  
\end{lem}
\proof Note that $\rho_0(b^2)= \rho_0(b^{-1})$. Thus, if $\rho_0(b^2)$ commutes 
with $\rho_0(a_3)$ then $\rho_0(a_3)^3=1$ which is impossible. 
Take a cocycle 
$\si\in Z^1(G^s, ad\circ\rho)$ which is tangent to the variety $V$. Hence
$\si(a_j)=0, 1\le j\le 3$. According to Proposition \ref{cycl} the value of $\si$ on the generator $b$ equals $\tau + b\xi -\xi$. On the other hand, by 
Proposition \ref{L5.5} we can find coboundaries $\del_{\theta_j} \in 
B^1(\< a_j, b\>, so(3,\C))$, $j=1,2$ so that:
$$
\si_j:= \si - \del_{\theta_j}, \si_j(a_j)=0, \si_j(b)= \tau
$$
The coboundary $\del_{\theta_j}$ is given by
$$
\del_{\theta_j}(x)=  x\theta_j -\theta_j, x\in \< a_j, b\>, 
\theta_j\in so(3,\C) 
$$
Thus 
$$
a_j \theta_j = \theta_j, b(\theta_1- \theta_2)= \theta_1- \theta_2
$$
Note however that the condition (3) on the representation $\rho_0$ implies that 
the (1-dimen\-si\-onal) fixed-point sets for the adjoint actions of 
$\rho_0(a_1), \rho_0(a_2), \rho_0(b)$ on $so(3,\C)$ are linearly independent. Therefore we conclude that $\theta_j=0$, $j=1,2$, thus 
$$
\si(b)= \tau
$$ 
and by Lemma \ref{L5.6} we have $\si(b)=0$. Hence $\si=0$ and the Zariski tangent space to $\rho_0$ in $V$ is zero. $\qed$

\subsection{Nondegenerate representations}
\label{appli}

Let $\Ga$ be a labelled graph where all vertices and edges have nonzero even 
labels. Let $G^s$ denote the Shephard group corresponding to the graph $\Ga$. 
The following technical definition will be used in 
Section \ref{Representations of Shephard groups}. 

\begin{defn}
A representation $\rho: G^s \to PO(3,\C)$ will be called 
{\bf nondegenerate on the edge} $e= [v, w]\subset \Ga$ if the element 
$\rho(g_v g_w)$ has the order 
$$
\left\{ \begin{array}{cc}
\eps(e) \ \ , & \hbox{ if $\del(v)= \del(w)=2$} \\
\eps(e)/2 \ \  , & \hbox{otherwise} \end{array} \right. 
$$
A representation 
$\rho$ will be called 
{\bf nondegenerate  on the vertex} $v\in V(\Ga)$ if $\rho(g_v)\ne \1$. A 
representation $\rho$ will be called {\bf nondegenerate} if it is 
nondegenerate on each edge and each vertex of $\Ga$. Let 
$\Hom_f(G^s, PO(3,\C))$ denote the space of all nondegenerate
 representations. 
\end{defn}

\begin{prop}
\label{P6}
Suppose that for each edge $e\subset \Ga$ the corresponding edge subgroup  
$G_e\subset G^s$ is finite. Then 
$\Hom_f(G^s, PO(3,\C))$  is Zariski open and closed in  $\Hom(G^s, PO(3,\C))$.
\end{prop}
\proof Since each $G_v, G_e \subset G^s$ is finite, 
$\Hom(G_v, PO(3,\C))$, $\Hom(G_e, PO(3,\C))$ are disjoint unions  
of finite numbers of $PO(3, \C)$-orbits of rigid representations. 
Since each orbit is Zariski open the proposition follows in the 
case $G^s= G_v, G_e$. 
Let 
$$Res_e: \Hom(G^s, PO(3,\C)) \to \Hom(G_e, PO(3,\C)),$$ 
$$Res_v: \Hom(G^s, PO(3,\C)) \to \Hom(G_v, PO(3,\C))
$$
 be the restriction morphisms. Then $\Hom_f(G^s, PO(3,\C)) = $ 
$$
\bigcap_{v\in \V(\Ga)} Res_v^{-1} 
\Hom_f(G_v, PO(3,\C)) \cap \bigcap_{e\in \E(\Ga)} Res_e^{-1} 
\Hom_f(G_e, PO(3,\C))
$$
and the proposition follows.  $\qed$

\section{Arrangements}
\subsection{Abstract arrangements}
\label{Abstract arrangements}

An {\it abstract arrangement} $A$ is a disjoint 
union  of two finite sets $A = \p \u \L$, with the  set of 
``points'' $\p = \{v_1\ , \, v_2\ , \cdots\}$ 
and the set of ``lines'' $\L = (l_1\ , l_2\ , \cdots)$ 
together with the {\em incidence  relation}  
$\iota = \iota_A \subset \p \times \L;
\iota (v, l)$ is interpreted to mean ``the point $v$ lies on the 
line $l$~". We may represent the arrangement $A$ by a 
bipartite graph, $\Gamma =\Gamma_A$: 
vertices of $\Ga$ are elements of $A$, two vertices are connected 
by an edge if and only if the corresponding elements of $A$ are 
incident ($\Ga$ is also called the {\em Hasse diagram} of 
the arrangement $A$).

\begin{figure}[tbh]
\leavevmode
\centerline{\epsfxsize=3.5in \epsfbox{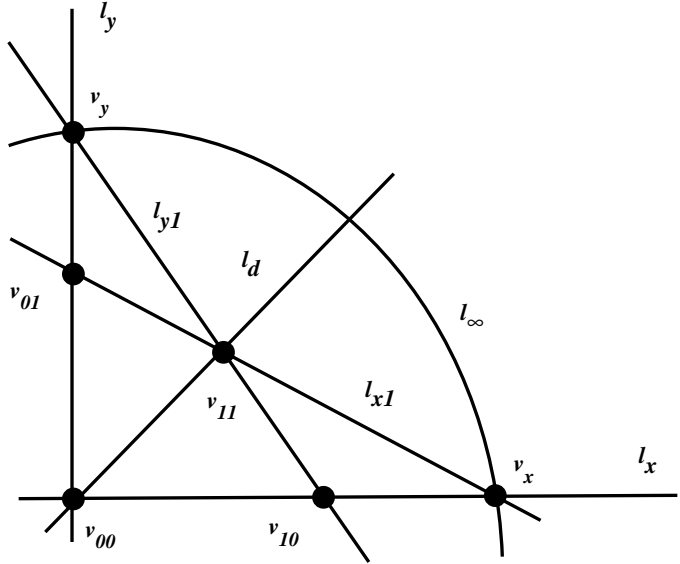}}
\caption{\sl Standard triangle.}
\label{Fig1}
\end{figure}

\begin{conv}
When drawing an arrangement we shall draw {\em points} as {\em solid points} and {\em lines} as {\em lines}. If $\i(v, l)$ then we shall draw the point 
$v$ on the line $l$. 
\end{conv}

An example of an abstract arrangement is the {\em standard triangle} 
$T$ described on the Figure \ref{Fig1} (this is a triangle with the complete set of bisectors): 
$$
T = \{v_{00}, v_x, v_y , v_{1,0}, v_{01}; l_x, l_y, l_{\infty}, 
l_d, l_{y1}, l_{x1}: \i(v_{00}, l_x), \i(v_{01}, l_y),   
\i(v_x, l_x), \i(v_x, l_{\infty}),$$
$$
  \i(v_y, l_y), \i(v_y, l_{\infty}), 
 \i(v_{11}, l_d), \i(v_{00}, l_d), \i(v_{00}, l_y), \i (v_{10}, 
l_{y1}), \i(v_{01}, l_{x1}), $$
$$
\i(v_{10}, l_x), 
%%%%%\i(v_{dd}, l_{\8}), \i(v_{dd}, l_d), 
\i(v_y, l_{y1}), \i (v_x, l_{x1}), \i(v_{11}, l_{y1}), \i(v_{11}, l_{x1}) \} 
$$

\begin{figure}[tbh]
\leavevmode
\centerline{\epsfxsize=3.5in \epsfbox{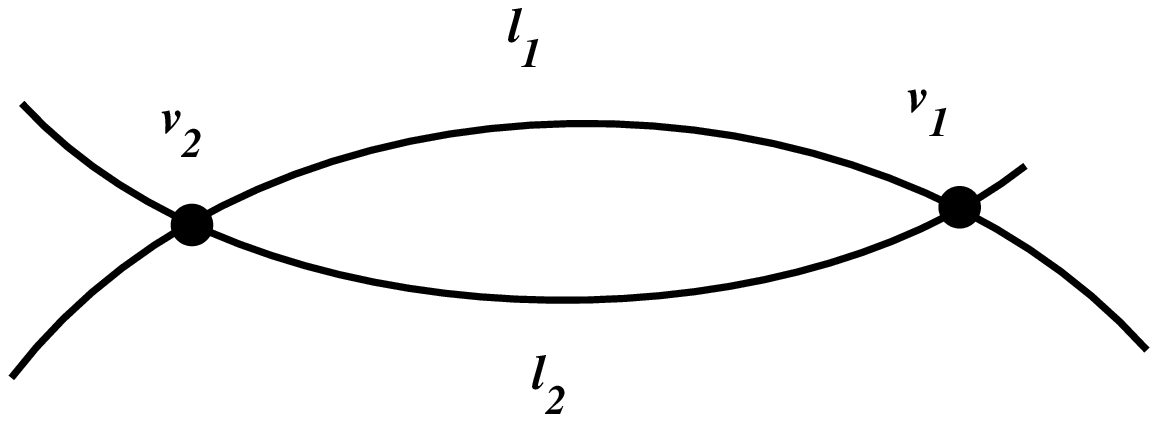}}
\caption{\sl Bigon.}
\label{Fig2}
\end{figure}

\medskip
Here is another example of arrangement (Figure \ref{Fig2}), we take 
$A= \{ v_1, v_2; l_1, l_2\}$ with the incidence relation:
$$
\i(v_1, l_1), \i(v_1, l_2), \i(v_2, l_1), \i(v_2, l_2) \quad .
$$
 
Suppose that $(A, \i_A), (B, \i_B)$ are abstract arrangements, 
$\phi: B\to A$ is a map which sends points to points and lines to lines. 
We say that $\phi$ is a {\em morphism} of arrangements if $\i_B(x, y)$ 
implies $\i_A(\phi(x), \phi(y))$. A {\em monomorphism} of arrangements 
is an injective morphism. 
An {\em isomorphism} of arrangements is an invertible morphism. Suppose 
that $\phi: B \to A$ is a monomorphism of arrangements, we call the 
image $\phi(B)$ a {\em subarrangement} in $A$. Note that if 
 we work with the corresponding bipartite graphs $\Ga_A, \Ga_B$, 
then morphisms $A\to B$ are morphisms of these bipartite graphs 
which send points to points, lines to lines. 

\begin{defn}
 An abstract {\bf based arrangement} $A$ is an arrangement together with\newline  a  
monomorphism of the standard triangle 
$T\hook A$ that we call the {\bf canonical} {\bf embedding}.
\end{defn}

\medskip
\no An arrangement $A$ is called {\em admissible} if it satisfies the axiom:   

\medskip 
{\bf (A1)} {\em Every element of $A$ is incident to at least two 
distinct elements.}  (I.e. every \newline point belongs to at least two 
lines and every line contains at least two 
points.)\footnote{This axiom will be needed only in the Section 
\ref{repart}.} 

\medskip
Suppose that $A$, $B$, $C$ is a triple of arrangements and 
$\phi: C\to A, \psi: C\to B$ are monomorphisms. We define the 
 {\em fiber sum} $A\times_C B$ as follows. First we take 
the disjoint union of the arrangements $A$ and $B$. Then
identify in $A\cup B$ the elements $\phi(c), \psi(c)$ for 
all $c\in C$. If $A, B$ are based arrangements, $C$ is as above, 
then their join $A*_C B$ is defined as $A\times_{T\u C} B$, where 
$T$ is the standard triangle with canonical embedding into $A, B$. 
If $C$ is an  arrangement which consists of a single {\em point} $c$ and 
$\phi(c)= a\in A, \psi(c)=b\in B$, then we denote
$$
A*_{a\equiv b} B := A *_C B 
$$

\subsection{Fiber products}
\label{3.2}

We remind the reader of the  definition of the  {\em fiber product} 
of varieties (recall that our varieties are neither reduced nor 
ireducible). Let $f: X\to Z, g: Y\to Z$ be morphisms. Then the 
{\em fiber product} $X\times_Z Y$ of $X$ and 
$Y$ with respect to $Z$ is a variety $X\times_Z Y$ together 
with canonical morphisms  $\Pi_X: X\times_Z Y \to X$ and 
$\Pi_Y: X\times_Z Y \to Y$
 such that the fol\-lo\-wing diagram is commutative 
\[
\begin{array}{ccc}
X\times_Z Y & \lra & X\\
 \downarrow & ~ & \downarrow \\
Y & \longrightarrow & Z\\
\end{array}
\]
These data satisfy the %universal 
property that given a 
variety $W$ and a commutative diagram
\[
\begin{array}{ccc}
W & \lra & X\\
 \downarrow & ~ & \downarrow \\
Y & \longrightarrow & Z\\
\end{array}
\]
we obtain a commutative diagram
\[
\begin{array}{ccccc}
W &  \RA &  X\\
~& \parbox[l]{.2in}{{\large $\searrow$}} 
\quad \quad\quad\quad \quad \parbox[r]{.3in}
{{\large $\nearrow$}} & ~ \\
\Big\downarrow &  X\times_Z Y &  \Big\downarrow\\
~ & \hbox{{\large $\swarrow$ \quad\quad \quad\quad  \quad }} &  ~ \\
Y &  \RA &  Z\\
\end{array}
\]
The fiber product of quasi-projective (resp. projective) varieties 
is again a quasi-projective (resp. projective). 

We recall  how to describe graphs of morphisms via fiber 
products. Suppose that \newline $f: X= \A^n \to \A^m= Z$ is a morphism. 
Let $g: Y= \A^m \to \A^m$ be the identity morphism. Define 
$\Ga_f$, graph of $f$, to be the subvariety 
$$
\Ga_f= X\times_Y Y = \{ (x, y)\in X\times Y: f(x)= g(y)= y\} 
$$
Clearly $\pi_X : \Ga_f \to X$ is an isomorphism of affine varieties. 

Let $f: X= \A^n \to Y= \A^m$ be a morphism, $\Ga_f$ its graph and  
$\pi_X: \Ga_f\cong \A^n$ be the canonical projection. We split 
$\A^n$ as $\A^{n-k}\times \A^k$, so $x\in \A^n$ is written as \newline  
$x= (x' , x'')$, $x'\in \A^{n-k}, x''\in \A^k$. We obtain projections
$$
\Pi'_X: \Ga_f \to \A^{n-k} \hbox{~~and~~} \Pi''_X: \Ga_f \to \A^{k}$$
 defined by $\Pi'_X(x)= x'$ and $\Pi''_X(x)= x''$. Now let 
$g: \A^s \to \A^{n-k}$ be a morphism with graph 
$\Ga_g$, let $y\in \A^s$ denote the variable in this space. 
Using the morphism\newline 
$\Ga_g\to \A^{n-k}$ (the second projection) and  
$\Pi'_X: \Ga_f \to \A^{n-k}$ 
we form the fiber product 
$$\Ga_g \times_{\A^{n-k}} \Ga_f\quad .$$ 
Now let $h: \A^s \times \A^k \to \A^m$ be the morphism given by 
$$
h(y, x'')= f(g(y), x'')
$$

\begin{lem}
\label{fib}
The projection map 
$$
p: (\A^s \times \A^{n-k})\times (\A^{n-k}\times \A^k) \times \A^m \to 
\A^s \times \A^k \times \A^m = T
$$
given by
$$
p((y, x'), (u', u''), z)= (y, u'', z)
$$
induces an isomorphism 
$$
\Ga_g \times_{\A^{n-k}} \Ga_f\lra \Ga_h \quad .$$
\end{lem}
\proof Obvious. $\qed$

\begin{cor}
The morphism 
$$
q: \Ga_g \times_{\A^{n-k}} \Ga_f\lra \A^s \times \A^k$$
 (given 
by the restriction of $p$ to $\Ga_g \times_{\A^{n-k}} \Ga_f$ and the 
the projection on the first and second factor of $T$) is an isomorphism. 
\end{cor}

\begin{cor}
\label{C3.2}
Suppose that  $k=0$, then the composition of $q$  with the 
projection on the first factor 
$$
r: \Ga_g \times_{\A^{n-k}} \Ga_f\lra \A^s
$$
is an isomorphism. 
\end{cor}

\subsection{Intersection operations in the projective plane}
\label{Intersection}

Let $\r$ be a commutative ring. We recall that the projective space $\P(M)$ 
for a projective $\r$-module $M$  of rank $n$ is defined by 
$$
\P(M):= \{ V\subset M: V ~~~ \hbox{is a projective submodule of rank 
1 such that~} M/V \hbox{~is projective}\}
$$
We then define $\P^n(\r)$ and $\P^n(\r)^{\vee}$ by
$$
\P^n(\r):= \P(\r^{n+1}), \quad \P^n(\r)^{\vee}:= \P(\Hom_{\r} 
(\r^{n+1}, \r))
$$
We refer to \cite[\S 1.3.4, \S 1.3.9]{DG}, to see that this is consistent 
with the usual definition of $\P^n$. 

Note that an element 
$\al\in \P^2(\r)^{\vee}$ gives rise to a projective  
$\r$-submodule $L\subset \r^3$ of rank $2$, $L:= \ker(\al)$, such that 
$\r^3/L$ is projective. We will call both $\al$ and $L$ {\em lines} in $\P^2(\r)$. We say that a point $V\in \P(M)$ belongs to a line $L$ (corresponding to $\al\in \P^2(\r)^{\vee}$ if and only if $V\subset L$; equivalently $\al(V)=0$. Suppose $V\in \P^2(\r)$ correspond to rank one 
{\em free} submodule of $\r^3$ with the basis $u= (x, y, z)$, then  
 we will write $V:= [x: y: z]$ (these are the homogeneous 
coordinates of $V$).

We now show how to do projective geometry over $\r$. We define two 
elements $Span(u_1)$, $Span(u_2)\in \P^2(\r)$ 
to be {\em independent} if the submodule 
$$
L= Span\{u_1, u_2\} := \r u_1 + \r u_2$$
 is a projective summand of $\r^3$ of rank $2$. In this case we will also 
say that $u_1, u_2$ are independent.

\begin{lem}
If $Span(u_1)$ and $Span(u_2)$ are independent over $\r$ then 
$Span\{u_1,u_2\}$ is the unique projective summand of $\r^3$ containing 
$u_1$ and $u_2$.
\end{lem}
\proof Let $M$ denote $Span\{u_1,u_2\}$. Suppose that $N$ is a projective 
summand containing $u_1$ and $u_2$. Then $N$ contains $M$. We want to 
prove that $M=N$. We may assume that $\r$ is local whence $M$ and $N$ are 
free. Let $\{v_1, v_2\}$ be a basis for $N$. Note that $\{u_1,u_2\}$ 
is a basis for $M$. Write
$$              
u_1 = av_1 + cv_2 , \quad u_2 = bv_1 + dv_2
$$
Let $\k$ be the residue field of $\r$. The images of $u_1$ and $u_2$ in 
$\k^3$  are independent so the image of $ad-bc$ in $k$ is nonzero. 
Hence $ad-bc$ is a unit in $\r$. $\qed$ 

\medskip
Thus two independent points $U_1= Span(u_1), U_2= Span(u_2) \in \P^2(\r)$ 
belong to the unique line $L= Span\{u_1,u_2\}$ in $\P^2(\r)$. We shall use 
the notation 
$$
L:= U_1 \bullet U_2
$$
for the line $L$ through the points $U_1, U_2$. 
If $u\in \r^3$ we let $u^{\vee}$ denote the 
element of $(\r^3)^{\vee}$ given by $u^{\vee}(v):= u\cdot v$.  
We have the following sufficient condition  for  independence:

\begin{lem}
 Suppose that there exists 
$u_3= (x_3, y_3, z_3)\in\r^3$ such that 
$(u_1 \times u_2)\cdot u_3= 1$. Then $u_1, u_2$ are independent, moreover
$$
Span\{u_1, u_2\} = \ker (u_1 \times u_2)^{\vee}
$$ 
\end{lem}
\proof The determinant of the matrix with the columns $u_1, u_2, u_3$ equals $1$, whence \newline 
$\{u_1, u_2, u_3\}$ is a basis for $\r^3$. Furthermore, suppose 
$v= a u_1 + bu_2 + cu_3$. Then $(u_1 \times u_2)\cdot v= c$, so $c=0$ if
and only if $v\in \ker (u_1 \times u_2)^{\vee}$. $\qed$

\begin{rem}
We observe that $u_3$ as above always exists (and hence $u_1, u_2$ 
are independent) if one of the coordinates of $u_1\times u_2$ is a 
unit in $\r$. In this case we will say that $u_1$ and $u_2$ satisfy 
 the {\bf cross-product test} for independence. 
\end{rem}

Dual to the correspondence $\bullet: \P^2 \times \P^2 \to (\P^2)^{\vee}$ 
there is  an operation of intersection of lines in $\P^2$. Namely, if 
$\la, \mu$ are  lines in 
$\P^2$ such that $\la^{\vee}, \mu^{\vee}$ are independent points in $\P^2$,  then we let 
$(\la\bullet \mu)^{\vee}= \la^{\vee}\bullet \mu^{\vee}$. 
Clearly $\la\bullet \mu= \ker(\la) \cap \ker(\mu)$. We have

\begin{lem}
If $u_1, u_2\in \R^3$ above satisfy  the {\bf cross-product test}, 
then $V_1= Span(u_1)$, $V_2= Span(u_2)$ can be joined by the 
unique projective line $V_1\bul V_2$ in $\P^2(\r)$ 
corresponding to $(u_1\times u_2)^{\vee}$. 
\end{lem} 

Dually let $L_j \in (\P^2)^{\vee}$ are dual to rank one free modules 
with bases $\si_j= (\al_j, \be_j, \ga_j)$, ($j=1,2$). We will write 
$L_j= [ \al_j: \be_j: \ga_j]$. We have 

\begin{lem}
If $\si_1, \si_2$ above are independent then $L_1, L_2$ intersect in the 
unique point $L_1\bul L_2$ with the homogeneous coordinates 
$[\si_1\times \si_2]$. 
\end{lem}

The {\em incidence variety} 
${\cal I}\subset \P^2 \times (\P^2)^{\vee}$ is given by the equations:
$$
\{ (p, l)\in \P^2 \times (\P^2)^{\vee} | l(p)=0\}
$$
Let $x, y, z$ be the coordinate functions on 
$\C^3$ relative to the standard basis and $\al,\be, \ga$ be the 
coordinate functions on $(\C^3)^*$ relative to the basis dual 
to $\ov{e_1}, \ov{e_2}, \ov{e_3}$. 
Then the homogeneous coordinate ring of ${\cal I}$ is isomorphic to  
$$
\frac{\C[x, y, z; \al, \be, \ga]}{(x \al + y\be + z\ga)}
$$
For a general commutative ring $\r$ the set of $\r$-points 
$I(\r)\subset \P^2(\r)\times \P^2(\r)^{\vee}$ 
consists of pairs $(V, \al)$ such that $V\subset L= \ker(\al)$. 

Pick a point $t\in \P^2 (\C)$. The {\em relative} incidence variety 
${\cal I}(t)\subset  (\P^2)^{\vee}(\C)$ is given by the equation:
$$
\{ l\in  (\P^2)^{\vee}(\C) | l(t) =0\}
$$
By dualizing we define the {\em relative} incidence variety 
${\cal I}(l)$ for any element $l \in (\P^2)^{\vee}(\C)$. 
We define {\em anisotropic} incidence varieties 
${\cal I}_0$, ${\cal I}_0(t)$ and ${\cal I}_0(l)$ 
by intersecting with $\P^2_0 \times (\P^2_0)^{\vee}$. 
The proof of the following lemma is a straightforward calculation 
and we leave it to the reader.

\begin{lem}
\label{is smooth}
For any $t$ and $l$ the varieties ${\cal I}$, ${\cal I}(t)$, ${\cal I}(l)$ 
are smooth. 
\end{lem}

\begin{Not}
We will make the following convention about inhomogeneous 
coordinates of points in $\P^2$: if $q= [x: y: 1]$ then we let $q:= (x, y)$, 
if $q= [0: 1: 0]$ we let $q:= (0, \8)$, if  $q= [1: 0: 0]$ we 
let $q:= (\8, 0)$ and if $q= [1:1:0]$ we let $q:= (\8, \8)$.  
\end{Not}

\subsection{Projective  arrangements}
\label{Projective  arrangements}
A {\em geometric realization} of the abstract arrangement 
$A= \p \u \L$ is a map
$$
\phi : \p \u \L \to \P^2(\C) \u (\P^2(\C))^{\vee} $$
which sends {\em points} to points and {\em lines} to lines. This map must 
satisfy the following condition:
\begin{equation}
%\label{inc}
\i(v, l) \Rightarrow \quad [\phi(v) \in \phi(l) \iff \quad 
\phi(l)^{\vee}\cdot  \phi(v) = 0\ ]
\end{equation}
The image of $\phi$ is a {\em projective  arrangement} in $\P^2$. Usually 
we shall denote {\em lines} of arrangements $A$ by uncapitalized letters 
($l, m, $ etc.), and their images under geometric realizations by 
the corresponding capital letters ($L, M, $ etc.). 

\begin{figure}[tbh]
\leavevmode
\centerline{\epsfxsize=3.5in \epsfbox{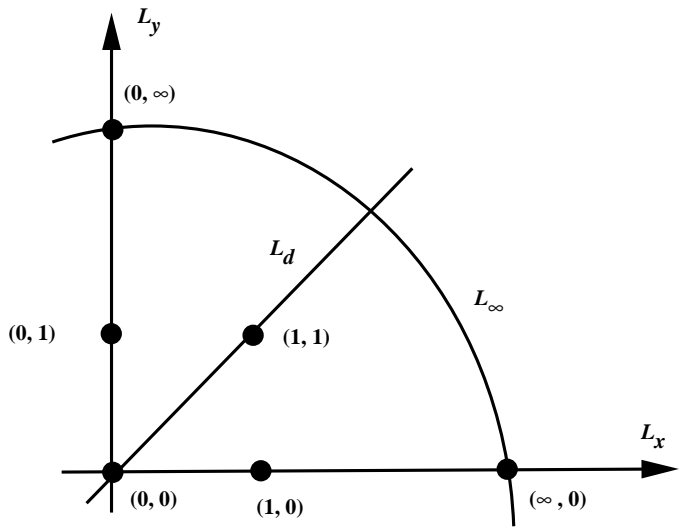}}
\caption{\sl Standard realization of the standard triangle.}
\label{Fig11}
\end{figure}

\begin{ex}
Consider the  standard triangle $T$. Define a geometric realization 
$\phi_T$ of $T$ in $\P^2$ so that:
$$
\phi_T(v_{00})= (0,0), \phi_T(v_{x})= (\8,0), \phi_T(v_{y})= (0,\8),
\phi_T(v_{11})= (1,1)
$$
This realization uniquely extends to the rest of $T$. 
See Figure \ref{Fig11}. We call 
$\phi_T$ the {\em standard} realization of $T$. 
\end{ex}

\medskip
The {\em configuration space} of an abstract arrangement $A$ is the 
space $R(A,\P^2(\C))$ of all geometric realizations of $C$. 
The space $R(A,\P^2(\C))$ 
is the set of $\C$-points of  a projective variety $R(A)$ 
defined over $\Z$  with equations determined by the condition (1). 

\begin{rem}
We will consider $R(A)$ as a variety over $\Q$. 
\end{rem}

We now give a concrete description of the homogeneous coordinate ring 
of $R(A)$ and the functor of points $R(A)(\bullet)$. Recall that $x, y, z$ 
be the coordinate functions on 
$\C^3$ relative to the standard basis and $\al,\be, \ga$ are the 
coordinate functions on $(\C^3)^*$ relative to the basis dual 
to $\ov{e_1}, \ov{e_2}, \ov{e_3}$. 
Now let $\Ga$ be the (bipartite) graph corresonding to the abstract 
arrangement $A$. We let $\{v_1,..., v_m\}$ denote 
the vertices of $\Ga$ corresponding to the {\em points} of $A$ and let 
$\ell_1,..., \ell_n$ denote the vertices of $\Ga$ corresponding to 
the {\em lines} of $A$. Let $P$ be the polynomial ring defined by
$$
P:= \C[x_1, y_1, z_1, ..., x_m, y_m, z_m; \al_1, \be_1, \ga_1,..., \al_n, 
\be_n,\ga_n]
$$
If $v_i$ and $\ell_j$ are incident we impose the relation $x_i \al_j + 
y_i \be_j + z_i \ga_j=0$. The resulting quotient of $P$ is the 
homogeneous coordinate ring of $R(A)$ which we denote $\C\, [R(A)]$. 
We next describe the functor of points $\r\to R(A)(\r)$ where $\r$ is a 
commutative ring. It follows from the description of $\P^2(\r)$ in 
\S \ref{Intersection} that the set $R(A)(\r)$ is described by 

\begin{lem}
\label{descr of points}
Let $\r$ be a commutative ring. Then an $\r$-point $\psi\in R(A)(\r)$ 
consists of an assignment of a {\em point} $V_i\in \P^2(\r)$ for each 
{\em point}-vertex $v_i\in \Ga$ and a {\em line} $L_j \in \P^2(\r)^{\vee}$  
for each {\em line}-vertex $\ell_j\in \Ga$ such that $\i(v_i, \ell_j)$ 
implies $V_i\subset L_j$. 
\end{lem}

\medskip
A {\em based realization} is a realization $\phi$ of a based 
abstract arrangement such that the restriction of $\phi$ to 
the canonically embedded triangle  $T$ is the {\em standard} 
realization $\phi_T$. 

The space of based realizations of an arrangement $C$ will be denoted 
by  $BR(A, \P^2(\C))$. It also is the set of  $\C$-points of a projective 
variety $BR(A)$ defined over $\Q$. Suppose that $A$ is a based arrangement 
which has distinguished set of {\em points} $\nu=\{v_1,..., v_n\}$ which 
lie on the {\em line} $l_x\in T$. We call $A_{\nu}$ a {\em marked} arrangement. 
A {\em morphism of marked arrangements} is a mapping 
$h: A_{\mu}\to B_{\nu}$ 
such that $h(\mu)\subset \nu$. Note that fiber sum of 
two marked arrangements has the natural structure of a marked 
arrangement.  If  $\nu=\{v_1,..., v_n\}$  is a marking 
we define the space of {\em finite realizations} as 
$$
BR_0(A_{\nu})= \{\psi\in BR(A): \quad \psi(v_j)\notin L_{\8}\ , 
\quad v_j\in \nu \}
$$
Clearly $BR_0(A)$ is a quasi-projective subvariety in $BR(A)$. 

\begin{lem}
The maps $R$, $BR$ and $BR_0$ define  contravariant functors from the 
category of abstract arrangements, based abstract arrangements and 
marked arrangements to the category of projective 
and quasi-projective varieties defined over $\Q$. 
\end{lem}
\proof Obvious. $\qed$

We leave the description of the coordinate rings and functors of 
points of $BR(C)$ and $BR_0(C)$ to the reader. 

\begin{thm}
\label{ftof}
Under the functors $R$, $BR$ and $BR_0$ the operation of fiber 
sum of arrangements corresponds to the operation of fiber 
product of projective varieties and quasi-projective varieties. 
\end{thm}
\proof Note that fiber sum of arrangements fits into the 
following commutative diagram of monomorphisms  of arrangements: 
\[
\begin{array}{ccc}
~ & A & ~ \\
~ & \nearrow \quad  \quad \searrow & ~  \\
~& C \longrightarrow A\times_C B & ~\\
 ~ & \searrow \quad  \quad \nearrow & ~  \\
~ & B & ~ \\
\end{array}
\] 
We shall consider the case of the functor $BR$, the other two 
cases are similar. It  follows from the previous lemma that we 
have commutative diagram of morphisms of varieties: 
\[
\begin{array}{ccc}
~ & BR(A) &  ~\\
 ~ & \nearrow \quad \quad \quad \quad \searrow & ~ \\
~ & BR(A\times_C B) \RA  BR(C) & ~ \\
 ~ & \searrow \quad \quad \quad \quad \nearrow & ~  \\
~ & BR(B) & ~ \\
\end{array}
\]
Thus, by the universal property we get a  morphism
$$
\la: BR(A\times_C B) \to BR(A)\times_{BR(C)} BR(B)
$$ 
To see that $\la$ is an isomorphism we have only to check that it 
induces bijections of $\r$-points, for each commutative ring $\r$ 
(see \cite[Proposition IV-2]{EH}). 
This is clear by Lemma \ref{descr of points}. $\qed$

\subsection{The moduli space of a projective arrangement}

In this section we will construct a distinguished Mumford quotient 
$R(A, \P^2(\C))\h PGL(3,\C)$ which we will refer to as the 
{\em moduli space} $\M(A, \P^2(\C))$ for a based arrangement 
$A$. Since the equations defining 
$R(A, \P^2(\C))$ are invariant under $PGL(3,\C)$ it suffices to construct a 
Mumford quotient of 
$$
(\P^2)^m \times ((\P^2)^{\vee})^n
$$ 
where $m$ is the number of {\em points} in $A$ and $n$ is the number of 
{\em lines}. 
The {\em Mumford quotient} $R(A, \P^2(\C))\h PGL(3,\C)$ of 
$R(A, \P^2(\C))$ will be then the subvariety of 
the quotient variaty cut out by the incidence equations. In the 
next section we will identify $\M(A, \P^2(\C))$ with $BR(A, \P^2(\C))$ 
which will be seen 
to be a cross-section to the action of $PGL(3,\C)$ on an open subvariety 
of $R(A, \P^2(\C))$. 

To construct a (weighted) Mumford quotient of  $(\P^2(\C))^m \times 
((\P^2(\C))^{\vee})^n$ we 
need a projective embedding. Such an embedding corresponds to a choice 
of polarizing line bundle over each factor of the product. 
Since the group of isomorphism classes of line bundles on $\P^2(\C)$  is 
infinite cyclic this amounts to assigning a positive integer weigth to 
each factor (i.e. to each vertex of the graph $\Ga$ of $A$). It will be 
more convenient to assign positive rational weigths to each vertex, 
then the integer weights are obtained by clearing the denominators. 
We choose a small positive rational number $\eps$ and assign the 
weight $\frac{1}{4}- \eps$  to each of the four {\em point vertices} 
$v_{00}, v_{x}, v_{y}, v_{11}$ of $T\subset A$ and the weight $\eps$ 
to all other vertices. Let $W:= \{v_{00}, v_{x}, v_{y}, v_{11}\}$. 

Note that all semistable configurations for 
the four-point case are stable, see below. 
Thus it is clear that if $\eps$ is small enough 
the calculation of stable and semistable points will reduce to the 
corresponding calculation for $PGL(3,\C)$ acting on the product of 
four copies of $\P^2(\C)$ corresponding to the four {\em point}-vertices 
described above. This calculation is well-known (see \cite{Newstead}). 
A configuration is stable (resp. semistable) iff less than (resp. 
no more than) $1/3$ of the total weight is concentrated on any point 
and less than $2/3$ (resp. no more than $2/3$) of the total weight 
is concentrated on any line. We obtain

\begin{lem}
\label{L8.6}
All semistable configurations in $R(A, \P^2(\C))$ are stable. 
A realization 
$\psi: A\to \P^2(\C) \u (\P^2(\C))^{\vee}$ is stable if and only if 
no three points of $\psi(W)$ lie on the same projective line in $\P^2(\C)$.  
\end{lem}

\subsection{The moduli space of the standard triangle}

As an important example we consider the configuration space 
and moduli space of the standard triangle $T$. First of all 
we note that $R(T, \P^2(\C) )$ is not irreducible. Here is the reason. 
Let $\phi_T\in R(T, \P^2(\C))$ be the standard realization. Then all 
$\phi\in R(T, \P^2(\C))$ nearby are equivalent to $\phi_T$ under some 
projective transformation. However there are other (degenerate) 
realizations $\psi_d\in R(T, \P^2(\C))$.  Namely, send all the {\em points} 
of $T$ to the origin $(0,0)\in \A^2\subset \P^2$. The triangle 
$T$ has 6 {\em lines}, which can be sent to any 6 lines in $\P^2$ 
passing through $(0,0)$. This gives us a 6-parameter family 
$F$ of degenerate realizations (which is the product of six  
copies of $\P^1$). After we mod out by the stabilizer of 
$(0,0)$ in $PGL(3, \C)$ we get  3-dimensional quotient. There 
are several other components which interpolate between $\phi_T$ 
and $\psi_d$, namely when three of the {\em points} 
$v_{00}, v_x, v_y, v_{11}$ belong to a common projective line. To 
remedy the problem we consider the {\em Mumford quotient} 
$R(T, \P^2(\C))\h PGL(3, \C)$ where we assign  weights as in the 
previous section.  Let $R_{ss}(T, \P^2(\C))$ be 
the set of semi-stable points with respect to these weights. 
Our definition of weights  prevents images  of any three of the 
{\em points}  $v_{00}, v_x, v_y, v_{11}$ from belonging to a 
common projective line in $\P^2(\C)$. It is clear that    
$R_{ss}(T, \P^2(\C))= \ PGL(3, \C) \{\phi_T\}$. Thus we get

\begin{lem}
The weighted quotient ${\cal M}(T)(\C) = R(T, \P^2(\C))\h PGL(3, \C)$ 
consists of a 
single point which we can identify with the cross-section $\{\phi_T\}$ 
for the action of  $PGL(3, \C)$ on $R_{ss}(T, \P^2(\C))$.
\end{lem}

\noindent Suppose now that $A$ is a general based arrangement, we assign 
weights as above. Then $BR(A, \P^2(\C))\subset 
R_s(A, \P^2(\C))= R_{ss}(A, \P^2(\C))$. By Lemma \ref{L8.6}, 
$BR(A, \P^2(\C))$ is a cross-section to the action of $PGL(3, \C)$ on 
$R_s(A, \P^2(\C))$. We obtain

\begin{thm}
\label{T8.7}
The inclusion $BR(A, \P^2)\hook R(A)$ induces an isomorphism 
$BR(A, \P^2(\C))$ $\cong \M(A, \P^2(\C))$ between projective varieties.  
\end{thm}

\subsection{Functional arrangements} 

Suppose that $A$ is a based marked arrangement 
with the {\em marking} $\nu=\{v_1,..., v_n\}$ (see \S 
\ref{Projective  arrangements}). 
 We call the points in $\nu$ the {\em input} points and we shall assume 
that $\nu \cap T= \0$. We also suppose that $A$ has the second marking 
$\mu= \{w_1,..., w_s\}$, which consists of distinct 
{\em output} points  $\{w_1,..., w_s\}$. (The sets $\mu$ and 
$\nu$ can intersect and we allow $\mu\cap T\ne \0$.) Recall 
that $\iota(v_i, l_x), \iota(w_j, l_x)$ for all $i, j$.  

Define the {\em projection maps} from the spaces of finite 
realizations of  $\Pi: BR_0(A_{\nu}) \to \A^n = \C^n$, 
$\Del: BR_0(A_{\nu}) \to \P^s$ by 
$$
\Pi: \phi \mapsto (\phi(v_1),..., \phi(v_n))= (z_1,..., z_n)\in \A^n
$$
$$
\Del: \phi \mapsto (\phi(w_1),..., \phi(w_s))= (y_1,..., y_s)\in \P^s
$$
(here we identify $L_x - \{(\8, 0)\}$ with the affine line $\A$). 

\begin{defn}
Suppose that arrangement $A$ above satisfies  the following axioms:  

\begin{itemize}
\item {\bf (A2)} $BR_0(A_{\nu})\subset BR_0(A_{\mu})$, i.e. 
$$ 
\psi(v_j)\notin L_{\8} , 1\le j\le n\quad \Rightarrow \quad 
\psi(w_i)\notin L_{\8} , 1\le i\le s~.
$$
\item {\bf (A3)}  The projection $\Pi$ is a biregular isomorphism of 
the variety $BR_0(A_{\nu})$ onto $\A^n$. 
\end{itemize}
Such arrangement $A$ is called a {\bf functional} arrangement on 
$n$ variables. 
\end{defn}

Each functional arrangement  {\em defines}
 a vector-function $f: \C^n \to \C^s$ by
$$
f(z_1,..., z_n) := \Del(\phi) = (\phi(w_1),..., \phi(w_s))$$
where $\phi\in BR_0(A_{\nu})$ corresponds to $(z_1,..., z_n)$ under 
the map $\Pi$. We shall record this by writing $A= A^f$. 
It's easy to see that the  vector-function $f$ must be polynomial. 
Later on we will give some examples of functional arrangements 
and we will prove that {\bf any} $m$-tuple of polynomials in 
$\Z[x_1,..., x_n]$ can be defined by a functional arrangement. 

\begin{lem}
The space $BR_0(A^f)$ is biregular isomorphic to the graph 
$\Ga_f$ of the function $f: \A^n \to \A^s$. 
\end{lem}
\proof Indeed, the natural projection $\pi: \Ga_f \to \A^n$ is 
an isomorphism. Compose it with the isomorphism $\Pi^{-1}$. The 
result is the required isomorphism $\Ga_f \to BR_0(A^f)$. $\qed$

\medskip
Now suppose that we are given two functional arrangements 
$A^f , A^g$ 
which define the functions $f: \C^n \to \C^s$ and $g: \C^t \to \C$. 
We denote the variables for $f$ by  $(z_1,..., z_n)$ and the variables for 
$g$ by $(x_1,..., x_t)$. We assume that they correspond to the 
{\em input points} $p_1,..., p_n$ and $q_1,..., q_t$ respectively. 
Denote by $w_0\in A^g$ the {\em output 
point}. We would like to find an arrangement which {\em defines} the 
function $h= f(g(x_1,..., x_t), z_2,..., z_n)$. To do this we 
let $A^h$ be the join $A^f *_{p_1\equiv w_0} A^g$.

\begin{lem}
\label{comp}
The arrangement $A^h$ is functional and it {\em defines} the polynomial  
$h$ which is the above composition of the functions $f, g$.
\end{lem}
\proof The Axiom (A3) for $A_h$ follows from the Axiom (A2) 
for the arrangements $A_{f}, A_g$. The Axiom (A4) follows from 
the fact that $BR_0(A^h)$ is the fiber product of the varieties 
$BR_0(A^f)$, $BR_0(A^g)$, see Lemma \ref{fib}.  $\qed$ 

\begin{figure}[tbh]
\leavevmode
\centerline{\epsfxsize=3.5in \epsfbox{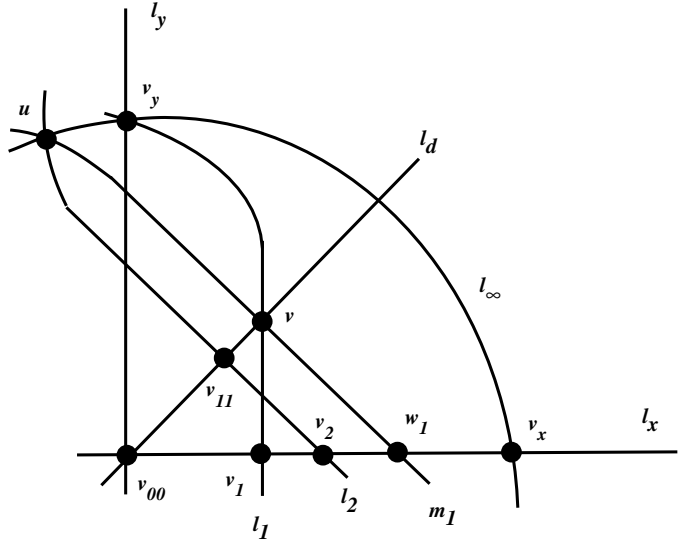}}
\caption{\sl The abstract arrangement $C_M$ for multiplication.}
\label{Fig3}
\end{figure}

\section{Algebraic operations via   arrangements}
\label{Algebraic operations}

\begin{figure}[tbh]
\leavevmode
\centerline{\epsfxsize=3.5in \epsfbox{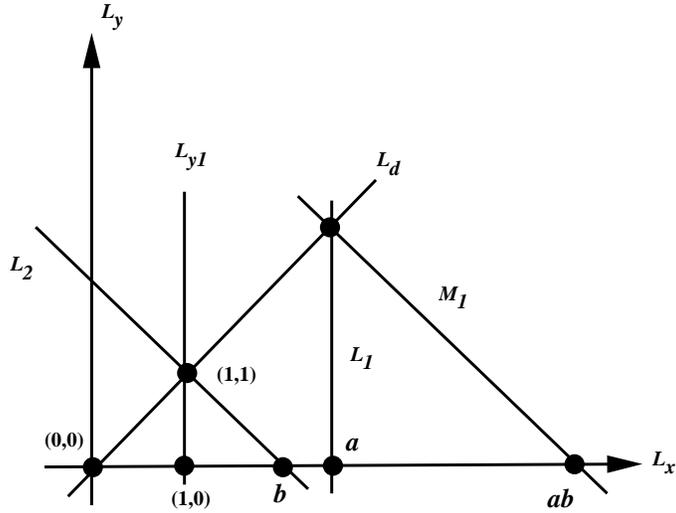}}
\caption{\sl Projective arrangement for multiplication.}
\label{Fig4}
\end{figure}

\begin{figure}[tbh]
\leavevmode
\centerline{\epsfxsize=3.5in \epsfbox{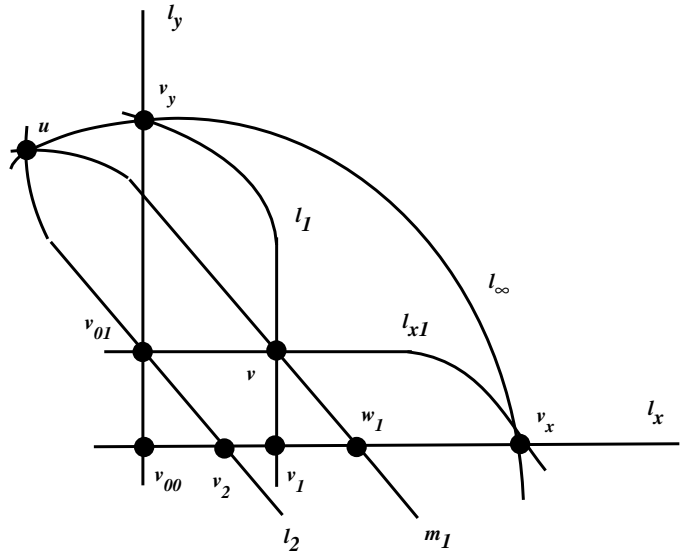}}
\caption{\sl The abstract arrangement $C_A$ for  addition.}
\label{Fig5}
\end{figure}

\begin{figure}[tbh]
\leavevmode
\centerline{\epsfxsize=3.5in \epsfbox{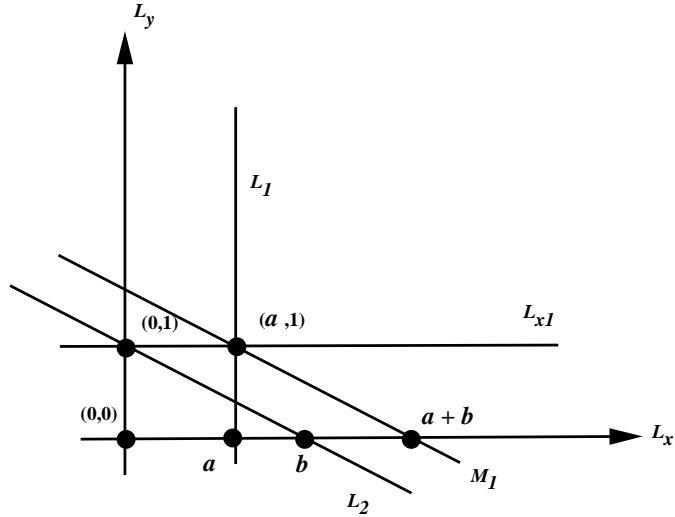}}
\caption{\sl Projective arrangement for addition.}
\label{Fig6}
\end{figure}

The following theorem goes back to the work of von Staudt \cite{St} in the middle of the last century: 

\begin{thm} 
\label{add}
There are admissible functional arrangements $C_A, C_M$ which {\em  define} 
the functions  
$$
A(z_1, z_2) =  z_1 + z_2~, \quad M(z_1, z_2) =  z_1 \cdot z_2$$
\end{thm}
\proof Consider the functional arrangement 
$C_M$ described on Figure \ref{Fig3}. We omitted from the 
figure several (inessential) {\em lines} and {\em points} of the 
standard triangle $T$, however we still assume that $C_M$ is a 
{\em based arrangement}. A generic projective realization $\psi$ of 
this arrangement is described on Figure \ref{Fig4}. 
Then the point of intersection of the line $M_1:= \psi(m_1)$ and the 
$x$-axis $L_x$ is equal to $ab$. (See \cite[page 45]{Hart}.)
The addition is {\em defined} via the abstract arrangement on the Figure 
\ref{Fig5}. Generic  projective realization of $C_A$ is described 
on the Figure \ref{Fig6}. (See \cite[page 44]{Hart}.) 
We will prove that $C_M$ is a functional arrangement, leaving 
similar case of $C_A$ to the reader. The Axioms {\bf (A1)}, {\bf (A2)} 
are clearly satisfied by the arrangement $C_M$\ , it is also 
easy to see that the morphism $\Pi: BR_0(C_M)\to \A^2$ is 
a bijection of complex points. 

The problem is to prove that $\Pi$ is invertible 
as a morphism. The example that the reader should keep in 
mind is the following. Consider the identity map
$$
id: \{ x= 0: x\in \C\} \lra \{x^2= 0: x\in \C\}
$$
Then $id$ is a  morphism which is bijective on complex points but 
not invertible as a morphism. 

We will prove that $\Pi: BR_0(C_M)\to \A^2$ induces a bijection 
$\Pi_R$ of $\r$-points for any commutative ring $\r$. This will 
imply that $\Pi$ is an isomorphism by \cite[Proposition IV-2]{EH}. 

Let $A_{\nu}$ be a based marked arrangement.  
We first interpret the finiteness condition 
$\psi(v_j)\notin L_{\8}, v_j\in\nu$,  
scheme-theortically for any $\r$-valued point $\psi\in BR_0(A)(\r)$. 
We will construct a subfunctor $U\subset \P^2$ corresponding to the 
affine plane $\P^2- L_{\8}$. Let $\r$ be a commutative ring, 
$\{e_1, e_2, e_3\}$ be the standard basis of $\r^3$ and let 
$M\subset \r^3$ be a submodule defined by $M= \r e_1 + \r e_2$. We define 
$U(\r)\subset \P^2(\r)$ by
$$
U(\r)= \{V \in \P^2(\r) : V \ \ \hbox{is a complement to}\ \   M\}
$$
It is then immediate, see \cite[ \S 1.3.9]{DG}, that $V\in U(\r)$ 
implies that $V$ is free and contains a unique vector 
(necessarily a basis) $u$ of the form 
$u= x e_1 + y e_2 + e_3$. The map $\omega: U(\r)\to \r^2 + e_3$ 
given by $\om(V)= (x, y,1)$ is a natural (with respect to  
$\r$) bijection. Consequently $U$ is represented by the polynomial 
algebra $\Z[X,Y]$. By \cite{DG}, loc. cit., $U$ corresponds 
to an open subvariety, again denoted $U$, of $\P^2$ 
isomorphic to $\A^2$ (by the previous sentence). We observe 
that $U(\C)= P^2(\C) - L_{\8}$. Thus $U$ gives us the 
scheme-theoretic definition of $\P^2- L_{\8}$. 

We now give the scheme-theoretic definition of the space of 
{\em finite} realizations $BR_0(A_{\nu})$ by requiring that 
$\psi(v_j)\in U$, $v_j\in \nu$. In the case 
$\psi\in BR_0(A_{\nu})(\r)$ we have 
$$
\om(\psi(v_j))= (x_j, y_j, 1), v_j \in \nu
$$ 
Moreover, since $\i(v_j, l_x)$ we have: $y_j=0, v_j \in \nu$. Thus we 
have associated a point $(x_1,..., x_n)\in \r^n$ to each 
$\psi\in BR_0(A_{\nu})(\r)$. We will see in what follows that 
the coordinates $(x_1,..., x_n)$ will play the same role as 
complex or real variables in the classical arguments.

We are now ready to prove

\begin{prop}
$\Pi_{\r}: BR_0(C_M)(\r) \lra \r^2$ is a bijection. 
\end{prop}
\proof Let $\psi\in BR_0(C_M)(\r)$. We will prove that $\psi$ is determined 
by $\psi(v_1), \psi(v_2)$ and the incidence relations in $C_M$. If 
$u= (x, y,z)$ we will use $[x: y : z]$ below to denote both the 
element in $\P^2(\r)$ with the basis $u$ and the element in 
$\P^2(\r)^{\vee}$ with the basis $u^{\vee}$; the meaning will be 
clear from the context. 

We now  prove the proposition. Let $\om(\psi(v_1))= (a, 0, 1)$ 
and $\om(\psi(v_2))= (b, 0, 1)$. Now 
$\om(\psi(v_1))\times (0, 1, 0)= (1, 0, -a)$. 
Hence by the cross-product test $[1: 0: -a]$ is the unique 
line joining $\psi(v_1)$ and $[0, 1, 0]$. But 
$\psi(l_1)= \psi(v_1)\bul L_y$, whence 
$\psi(l_1)= [1: 0: -a]$. We continue in this way, at each stage the 
cross-product test applies and we find in order
$$
\psi(v)= \psi(l_1)\bul \psi(l_d)= [a: a: 1]
$$
$$
\psi(l_2)= \psi(v_{11})\bul \psi(v_2)= [1: b-1: -b]
$$
$$
\psi(u)= \psi(l_2)\bul \psi(l_{\8})= [b-1: -1: 0]
$$
$$
\psi(m_1)= \psi(u)\bul \psi(v)= [-1: 1-b: ab]
$$
$$
\psi(w_1)= \psi(m_1)\bul \psi(l_x)= [ab: 0: 1] 
$$
This concludes the proof of Theorem \ref{add} in the case of the 
arrangement $C_M$. The argument in the case of $C_A$ is 
similar and is left to the reader.   $\qed$

\begin{lem}
\label{min}
For the function $D(x, y)= x-y$ there is a functional 
arrangement $C_D$ which defines the function $D(x, y)$.
\end{lem}
\proof Take the arrangement $C_A$ corresponding to the function 
$A(x, y)= x+y$ and reverse the roles of input-output {\em points}  
$v_2, w_1$. $\qed$ 

\begin{rem}
\label{use}
Suppose that $A$ is one of the arrangements $C_A, C_M, X_D$ and \newline 
$\psi\in BR(A,\P_2(\C))$  is a realization such that 
$\Pi(\psi)= (0, 0)$, hence $\psi(w_1)= 0$, $\psi(m_1)= L_y$. 
Then $\psi(A)\subset \phi_T(T)\cup \{(\8, \8)\}$. 
(Recall that $\phi_T$ is the 
standard realization of the standard triangle.) Let's 
verify this for $C_M$. If $\psi(v_1)= 0, \psi(v_2)= 0$ 
then $\psi(l_1)= \psi(l_y), \psi(l_2)= \psi(l_d)$. 
Hence $\psi(u)= (\8,\8)$ and $\psi(v)= (0,0)$.  
\end{rem}

\begin{figure}[tbh]
\leavevmode
\centerline{\epsfxsize=3.5in \epsfbox{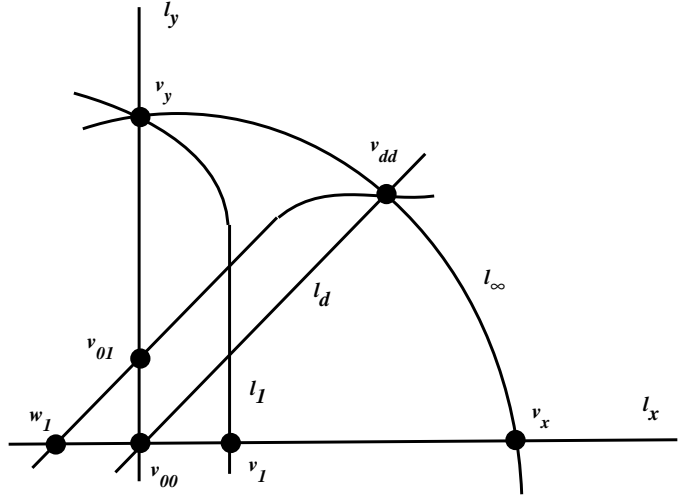}}
\caption{\sl Projective arrangement for the constant function 
$f_-(z)= -1$.}
\label{Fig7}
\end{figure}

\begin{lem}
\label{const}
There exist admissible functional arrangements $C^+$, $C^-$ which define 
the constant functions 
$$
f_+(z) = 1, f_- (z) = -1$$
\end{lem}
We describe configuration only for the function $f_-$, the 
other case is similar. Consider the arrangement on the Figure 
\ref{Fig7}, as usual we omit inessential lines and points 
form the standard triangle. Let $v_1$ be the input and  $w_1$ 
be the output points. 
%I.e. we add to the standard 
%triangle $T$ two extra {\em lines} $l_1, m_1$ and  extra {\em points} 
%$v_1$ (input) and  $w_1$ (output) such that:
%$$
%\i(v_{01}, m_1), \i(v_{dd}, m_1), \i(w_1, m_1), \i(v_1, l_1), 
%\i(v_y, l_1)
%$$
Under each realization $\phi$ of $C^-$ the image $\phi(w_1)$ 
is the point $(-1, 0)$. The image of $l_1$ is any vertical line in $\A^2$. We 
leave the proof of the fact that the arrangement $C^-$ is 
functional to the reader. $\qed$

\begin{cor}
\label{polyn}
For any polynomial $f$ of $n$ variables with integer coefficients 
there exists a functional arrangement $A^f$ which {\em defines} $f$.
\end{cor}
\proof Any such polynomial is a composition of the constant functions $f_{\pm}$, addition and multiplication.  
Thus the asssertion is a straightforward  application of Lemmas \ref{comp},  \ref{add}, \ref{min}, \ref{const} and Corollary \ref{ftof}.   $\qed$

\medskip
Now we  construct arrangements corresponding to polynomial 
vector-functions defined over $\Z$.

\begin{lem}
\label{eq}
Suppose that $f_1,..., f_m \in \Z [x_1,..., x_n]$. Then there 
exists an admissible abstract functional arrangement $A^F$ which {\em defines} the vector-function $F= (f_1,..., f_m)$. 
\end{lem}
\proof By Corollary \ref{polyn} there exist functional 
arrangements $A^{f_1},..., A^{f_m}$ which {\em define} the polynomials 
$f_j$. Let $v_{ij}$ be input point of $A^{f_i}$ corresponding to 
the variable $x_j$, $1\le i\le m$, $1\le j \le n$. Let 
$C= \{v_1,..., v_n\}$ be the arrangement which consists only of 
{\em points} $v_j$ and has no {\em  lines}, $B= C \cup T$. For each 
{\em point}  $v_j$ define $\psi_i : v_j \mapsto v_i$, this gives 
us embeddings $\psi_i: C\hook A^{f_i}$. Using these embeddings 
define  $A= (A^{f_1}* ... * A^{f_m}) \times_{C} B$. In  the 
arrangement $A$ we have $v_{ij}\equiv v_{kj}\equiv v_j$, 
$1\le j, k\le n$, these are the {\em input points} of $A$. 
The output points $w_1,..., w_m$ correspond to the {\em  output 
points} of the functional arrangements  $A^{f_1},..., A^{f_m}$. 
Then the fact that $A= A^F$ is the functional arrangement 
{\em defining} the vector-function $F$ follows by iterated application 
of Corollary \ref{C3.2} and Theorem \ref{ftof} similarly to the 
proof of Lemma \ref{comp}. $\qed$  

\section{Systems of polynomial equations}

Suppose that we have a system of polynomial equations defined 
over $\Q$:
$$
\Si= \left\{ \matrix f_1(x)= 0\cr f_2(x)= 0\\ \vdots\\ f_m(x)= 0 \endmatrix 
\right.
$$
where $x= (x_1,..., x_n)$, $x_j \in \C$. These equations determine an affine 
variety $S\subset \A^n$ defined over $\Q$. In the previous section we had 
constructed a functional arrangement $A= A^F$ which {\em  defines} 
the vector-function $F= (f_1,..., f_m)$.  Recall that we have 
two projection morphisms 
$$
\Pi: BR_0(A, \P^2)\to \A^n\ , \quad\Del: BR_0(A, \P^2)\to \A^m
$$ 
so that the diagram of morphisms 
\[
\begin{array}{ccc}
BR_0(A, \P^2) & \stackrel{\Pi}{\longleftrightarrow} & \C^n\\
\Del {\big\downarrow} & ~ & F{\big\downarrow} \\
\C^m & \= & \C^m\\
\end{array}
\]
is commutative. By the Axiom (A3) the projection $\Pi$ is an 
isomorphism. Let $T\subset A$ be the standard triangle, 
$w= v_{00}$ be its vertex. Define the new abstract arrangement 
$A^{\Si}$ as the join
$$
A^{\Si}=(\ldots ((A*_{w_1\equiv w} T)*_{w_2\equiv w}T)
\ldots *_{w_m\equiv w}T) 
$$
Then $BR_0(A^{\Si})$ is the fiber product 
$\{\psi \in BR_0(A, \P^2): \Del(\psi)= 0\}$ 
and $BR_0(A^{\Si}) \cong \{x \in \A^n : F(x)=0\}$, where the 
isomorphism is given by the restriction $\tau$ of $\Pi$ 
to $BR_0(A)$. Thus we get the following

\begin{thm}
\label{systems}
For any system of polynomial equations with integer coefficients
$$
\Si=\left\{ \matrix f_1(x)= 0\cr f_2(x)= 0\cr \vdots\cr f_m(x)= 0 
\endmatrix \right.
$$
there is an admissible based  arrangement $A= A^{\Si}$ such that 
$\tau: BR_0(A, \P^2)\cong S$ 
is an isomorphism of quasi-projective varieties over $\Q$.  
\end{thm}

\begin{defn}
We call the morphism $geo= \tau^{-1}$ {\bf geometrization}: it allows us to do 
algebra  (i.e. solve the system $\Si$ of algebraic equations) via geometry 
(i.e. by studying the space of projective arrangements).  
\end{defn}

Note that the arrangement $A= A^{\Si}$ is not uniquely determined by the affine 
variety $S$ but also by its affine embedding 
(the system of polynomial equations) and particular formulae used to describe 
these equations. For instance, the equation $x^5=0$ can be described as 
$(x^2 \cdot x^2)\cdot x$ as well as $x\cdot (x^2 \cdot x^2)$ and 
$(x\cdot (x^2))\cdot x^2$, etc. 

\medskip
Suppose that the system of equations $\Si$ is defined over $\Z$ and 
has no constant terms. Then 
the only subarrangements involved in construction of the arrangement 
$A^{\Si}$ are arrangements for multiplication, addition and subtraction 
(described in Lemmas  \ref{add}, \ref{min}) and we do not need 
arrangements for the constant functions $\pm 1$. Let 
$\psi_0\in BR(A, \P^2)$ be the realization corresponding to the point 
$0\in S$ under the isomorphism $\tau$. Take any {\em  line} $l\in A- T$ and a 
{\em  point} $v\in A-T$. Then by using Remark \ref{use} and the fact 
that arrangement for the composition of functions is a join of their 
arrangements (so it has no new points or lines) we conclude 
that the following holds:

\begin{lem}
\label{noco}
$\psi_0(l)$ is one of the {\em lines} $L_x, L_y, L_d$
and $\psi_0(v)$ is one of the {\em points} 
$$(0,0)\ , (0,\8)\ , (\8, 0)\ , (\8, \8)$$ 
for each {\em line} $l\in A- T$  and each {\em point} $v\in A- T$. 
\end{lem} 

Suppose now that $Q$ is an affine variety defined over $\Q$ and $q\in Q$ is 
a rational point. Then we can realize $Q$ as an affine subvariety 
$S \subset \A^N$ defined over $\Z$ so that $q$ goes to zero. Hence we get the following

\begin{cor}
\label{csystems}
For any affine variety $Q$ defined over $\Q$ and  
a rational point $q\in Q$ there exists an abstract admissible arrangement 
$A$ and an isomorphism\footnote{of varieties defined over $\Q$}
$geo: Q \to BR_0(A)$ so that the point $q$ corresponds to a realization 
$\psi$ so that  $\psi(l)$ is one of the {\em lines} $L_x, L_y, L_d$
and $\psi(v)$ is one of the {\em points} 
$$(0,0)\ , (0,\8)\ , (\8, 0)\ , (\8, \8)$$ 
for each {\em line} $l\in A- T$  and each {\em point} $v\in A- T$. 
\end{cor}

\section{Groups corresponding to abstract arrangements}
\label{gaa}

We will define several classes of groups corresponding to 
abstract arrangements. Let $\Ga$ be the bipartite graph corresponding to 
an abstract arrangement $A$. 
We first construct  the Coxeter group $G^c_A$ without assuming 
that $A$ is a based arrangement: we assign the label $2$ to all edges of 
 $\Ga$ and let $G^c_A:= G^c_{\Ga}$.   

From now on we suppose that $A$ is a based arrangement. 
We start by identifying the {\em point} $v_{00}$ with the {\em line} 
$l_{\8}$, the {\em point} $v_x$ with the {\em line} $l_y$ and the {\em point} 
$v_y$ with the {\em line} $l_x$ in the standard triangle $T$. We also 
introduce  the new edges 
$$
[v_{10}, v_{00}], \quad [v_{01}, v_{00}], \quad [v_{11}, v_{00}]
$$
(Where $v_{10}, v_{00}, v_{11}, v_{01}$ are again {\em points} in the standard triangle $T$.) We will use  the notation $\La$ for the resulting graph. 
Put the following labels on the edges of $\La$:  

1) We assign the label $4$ to the edges $[v_{10}, v_{00}], 
[v_{01}, v_{00}]$ and all the edges which contain $v_{11}$ 
as a vertex (with the exception of $[v_{11}, v_{00}]$). We put 
the label $6$ on the edge 
$[v_{11}, v_{00}]$. 

2) We assign the label $2$ to  the rest of the edges. 

\no Let $\Del:= \La- [v_{11}, v_{00}]$. 
Now we have  labelled graphs and we use the procedure 
from the Section \ref{groups} to construct:
 
\medskip
(a) The {\em Artin group} $G_A^a:= G_{\La}^a$\ .

\medskip
(b) We assign the label $3$ to the vertex $v_{11}$ and labels 
$2$ to the rest of the vertices. Then we get the {\em Shephard} 
group $G_A^s:= G_{\La}^s$\ . 

\medskip
\no We will denote generators of the above groups $g_v, g_l$, 
where $v, l$ are elements of $A$ (corresponding to vertices of 
$\La$).

\begin{rem}
Suppose that $A$ is a based arrangement. Then all the group 
$G_A^s$ admit an epimorphism onto free products 
of at least 3 copies of $\Z/2$ and the group $G^a_A$ has an epimorphism 
onto a free group of rank $r\ge 3$, where $r+3$ is the number of {\em lines} 
in $A$. Let's check this for $G= G_A^a$. Construct a new arrangement 
$B$ by removing all the {\em points} in $A$ (and the {\em lines} 
$l_x, l_y, l_{\8}$). Then 
$G_B^s\cong \Z * ...* \Z$ is the $r$-fold free product, 
where $r\ge 3$ is the number of lines in $A- \{l_x, l_y, l_{\8}\}$. 
It's clear that we have an epimorphism $G_A^a \to G_B^a$. Hence 
all the groups   $G_A^a, G_A^s$  are not virtually solvable. 
Suppose that $G^a_A$ is the fundamental group of a smooth 
complex quasi-projective surface $X$. Then it is known %\cite{Bauer} 
that $X$ 
admits a nonconstant holomorphic mapping onto a Riemann surface 
$Y$ of negative Euler characteristic. 
\end{rem}

\begin{figure}[tbh]
\leavevmode
\centerline{\epsfxsize=3.5in \epsfbox{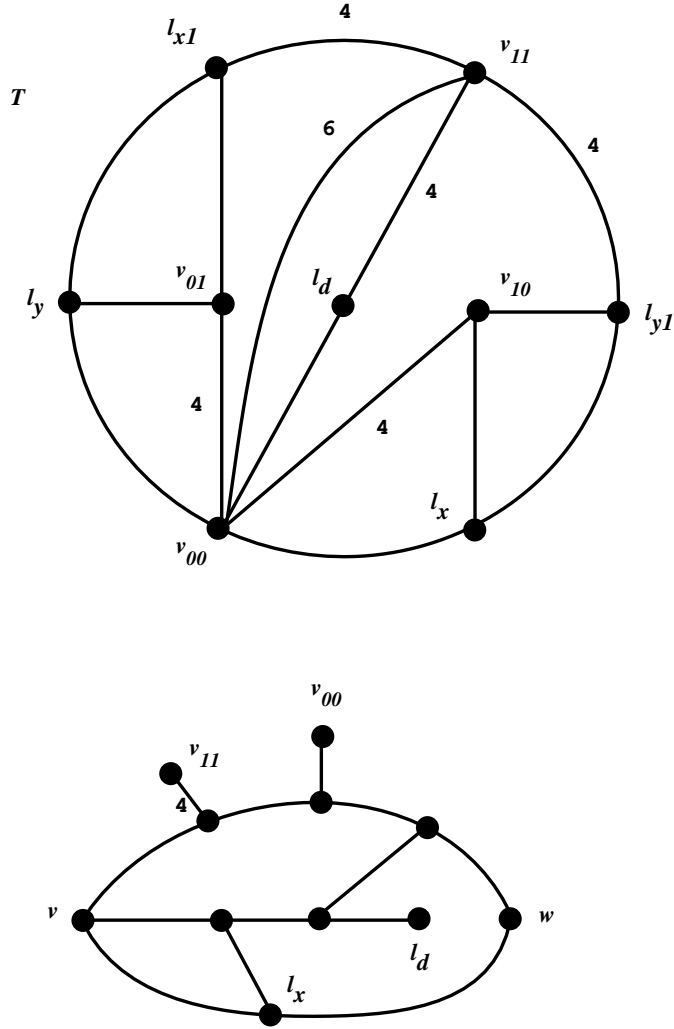}}
\caption{\sl Labelled graph of the functional arrangement for the 
function $x^2$. Identify vertices with the same labels. The {\em point} 
$v$ is the ``input'', the {\em point} $w$ is the ``output''. }
\label{Fig10}
\end{figure}

As an illustration we describe an example of a labelled graph 
corresponding to the based functional arrangement {\em defining} the 
function $x\mapsto x^2$, see the Figure \ref{Fig10}.  

\section{Representations associated with projective arrangements}
\label{S14}

This section is in a sense the heart of the paper. 
We start with an 
outline of the main idea behind it\footnote{Certain versions of this idea 
were used previously in our papers \cite{KM1}, \cite{KM2}, \cite{KM3}.}. 
A projective arrangement $\psi$  is {\it anisotropic} 
of $\psi (v) \in \P^2_0, \psi (l) \in (\P^2_0)^{\vee}$,
all points and lines $v, l\in A$. The anisotropic condition 
defines Zariski open subsets 
of the previous arrangement varieties to be denoted $R (A, \, \P^2_0), 
BR(A, \P^2_0)$ and $BR_0(A, \P^2_0)$ respectively.
 As we already saw in the Section 
\ref{S3}, the correspondence between involutions in $\P^2$ and their 
isolated fixed points is a biregular isomorphism between $PQ_0$ and 
$\P^2_0$: any point $x$ in $\P^2_0$ uniquely determines the  
``Cartan involution'' around this point, i.e. the involution such 
that $x$ is the isolated fixed point. The point $x= (1,1)$ in 
$\A^2\subset \P^2$ 
also determines (uniquely up to inversion) the 
rotation of the order 3 around $x$, so that $x$ is the neutral fixed 
point. (We will choose one of these rotations once and for all.) 
 Similarly any line $L\in (\P^2_0)^{\vee}$  uniquely 
determines the reflection which keeps  $L$ pointwise fixed. Finally we  
encode the incidence relation between points and lines in $\P^2$ 
using algebra: two involutions generate the subgroup 
$\Z/2\times \Z/2$ in $PO(3)$ iff the isolated fixed point 
of one belongs to the fixed line of another, rotations of the 
orders $2$ and $3$  anticommute iff the neutral fixed point 
of the rotation of the order $3$ belongs to the fixed line of the 
involution, etc. Thus, given a geometric object (projective 
arrangement) we can constructs an algebraic object (a representation 
of the associated Shephard group). We call the 
 mapping 
$$
alg: \hbox{~~arrangements~~} \longrightarrow \hbox{~~representations}
$$
{\em algebraization}. This mapping is the key in passing from 
realization spaces of projective arrangements to representation 
varieties. The fact that this correspondence is a homeomorphism 
between sets of $\C$-points of 
appropriate subvarieties will be more or less obvious, 
however we will prove more: the {\em algebraization} is a 
biregular isomorphism of certain (quasi-) projective varieties, 
the latter requires more work.

\subsection{Representations of Shephard groups}
\label{Representations of Shephard groups}

Let $A$ be an abstract based arrangement with the graph $\Ga_A$ 
and $G^s_A$ be the corresponding Shephard group with the graph $\La_A$. 
For all edges $e$ in the graph $\La_A$ 
the edge subgroups $G_e\subset G^s_A$ are finite. Recall that in 
Section \ref{appli} 
we had defined the space $\Hom_f(G, PO(3, \C))$ of nondegenerate 
representations of Shephard groups $G$. Thus as a direct 
corollary of Proposition \ref{P6} we get 

\begin{cor}
\label{zaropen}
$\Hom_f(G^s, PO(3,\C))$ is Zariski open and closed in $\Hom(G^s, PO(3))$.
\end{cor}

The significance of nondegenerate representations of Shephard groups 
is that they correspond to projective arrangements under the ``algebraization'' 
morphism $alg$.

\begin{defn}
Suppose that $\rho\in \Hom_f(G^s, PO(3,\C))$ is a representation, 
then we associate a projective arrangement 
$\psi= \al(\rho)\in R(A, \P^2_0(\C))$ as follows:

(a) If $v\in A$ is a point then we let $\psi(v)$ be the neutral fixed point 
of $\rho(g_v)$.

(b) If $l\in A$ is a line then we let $\psi(l)$ be the fixed line 
of the involution $\rho(g_v)$. 
\end{defn}

Clearly $\psi(v)\in \P^2_0$ and $\psi(l)\in (\P^2_0)^{\vee}$. 
Now we verify that $\psi$ respects the incidence relation. 
Suppose that $e= [v, w]$ is an edge in $\La_A$ where $v, w, e$ 
have the label $2$. Then, 
$\rho\in  \Hom_f(G_e, PO(3,\C))$ implies that $\rho(g_v), \rho(g_w)$ 
{\em anticommute} and hence $\psi(v)\cdot \psi(l)= 0$ (see 
\S \ref{Commuting}). All other edges $e= [v, l]$ in $\La_A$ are labelled 
as: $\del(v)=3, \del(l)=2, \eps(e)= 4$, thus $\rho\in  
\Hom_f(G_e, PO(3,\C))$ again implies that 
$\rho(g_v), \rho(g_l)$ {\em anticommute}:
$$
(\rho(g_v)\rho(g_l))^2 = \1
$$
and hence $\psi(v)\cdot \psi(l)= 0$ (see \S  \ref{Commuting}). 
Note that the mapping 
$$
\al: \Hom_f(G^s, PO(3,\C)) \to R(A, \P^2_0(\C)), \quad 
\rho\mapsto \psi= \al(\rho)
$$ 
is  2-1. Namely, we can modify any representation $\rho \in  
\Hom_f(G^s, PO(3,\C))$ by taking 
$\rho_-(g_{v_{11}}):= \rho(g_{v_{11}})^{-1}$ and 
$\rho_-(g_{w}):= \rho(g_{w})$ for all $w\in V(\Ga)- \{v_{11}\}$, then 
$\al(\rho)= \al(\rho^-)$. We denote the involution $\rho \mapsto \rho_-$ by 
$\om$. The mapping  $\al$ is far from being onto $R(A, \P^2_0(\C))$ 
because of extra identifications and edges in the graph $\La$ 
(comparing to the $\Ga_A$). The mapping $alg$ will be the 
right-inverse to $\al$ if we restrict the targent of $\al$ 
to {\em based realizations}:

\begin{defn}
Suppose that $A$ is an abstract based arrangement,  
$\psi\in BR(A, \P^2_0)$ is a realization.\footnote{In particular 
$\psi(v_{00})= \psi(l_{\8})^{\vee}, \psi(v_x)= 
\psi(l_y)^{\vee}, \psi(v_y)= \psi(l_x)^{\vee}$.} 
We construct a homomorphism 
$$alg(\psi)= \rho_{\psi}: G^s_A\to PO(3,\C)$$ 
as follows. 
If $v$ is a {\em point} in $A- \{v_{11}\}$, let $\rho(g_v)$ be the 
rotation of order $2$ in $\P^2$ with the neutral fixed point 
$\psi(v)$. (Such rotation exists since $\psi(v)$ is anisotropic.)  
If $v=l$ is a {\bf line} in $A$ we let $\rho(g_v)$ be the 
reflection in the line $\psi(l)$ (equivalently this is 
the rotation of order $2$ with the isolated fixed  point  
$\psi(l)^{\vee}\in \P^2_0$. For the vertex 
$v= v_{11}$   we take the rotation of 
the order $3$ around the point $(1,1)$. 
There are two such rotations, so we shall 
choose $\rho(g_v): \psi(v_{00})\mapsto \psi(v_x)$. We call the 
mapping $alg: \psi \mapsto \rho_{\psi}$,  {\bf algebraization}. 
\end{defn}

Below we verify that $\rho$ respects relations in $G^s$ and determines 
 a nondegenerate representation.  For all $v\in \p- \{v_{11}\}\subset A$ 
and each {\em line} $l \in \L \subset A$ we have  
$$
\psi(v)\in \psi(l) \Rightarrow  [\rho(g_l),\rho(g_v)]=\1, \quad 
\rho(g_l)\ne \rho(g_v)
$$
Hence all the relations in $G^s$ corresponding to the 
edges labelled by $2$ are preserved by $\rho$ and $\rho$ is 
{\em nondegenerate} on such edges. Since the spherical distances 
between the points $(0,0)$ and $(1, 0)$, and points  
$(0,0)$ and $(0, 1)$ in $\R\P^2$ equal $\pi/4$ we conclude that 
$$
(g_{v_{00}} g_{v_{10}})^4= \1~~, \quad (g_{v_{00}} 
g_{v_{01}})^4= \1 , (g_{v_{00}} g_{v_{10}})^2\ne \1, (g_{v_{00}} 
\ne g_{v_{01}})^2\ne \1 
$$
This implies that the Artin relations for the edges $[v_{00}, v_{10}], 
[v_{00}, v_{01}]$ are preserved by $\rho$ and $\rho$ is nondegenerate 
on these edges. 

For each {\em line} $l$ incident to the {\em point} 
$v_{11}$ the rotations $\rho(g_l)$ and $\rho(g_{v_{11}})$ 
anticommute by Lemma \ref{anticom} and  we get:
$$
\rho(g_{l}) \rho(g_{v_{11}}) \rho(g_{l}) 
\rho(g_{v_{11}})= \1~,\quad 
\rho(g_{v_{11}}) \rho(g_{l}) \rho(g_{v_{11}}) 
\rho(g_{l}) = \1
$$
Hence the  Artin relations for the edges $[l, v_{11}]\subset 
\La$ are preserved by $\rho$. 
It's  easy to check that the order of the element 
$\rho(g_{v_{00}} g_{v_{11}})$ equals $3$. Recall that the 
edge $e= [v_{00}, v_{11}]$ in $\La$ has the label $6$. 
Thus  the relation
$$
(g_{v_{00}} g_{v_{11}})^3 = (g_{v_{11}} g_{v_{00}})^3
$$
associated with the edge $e$ is preserved by $\rho$. Thus we have proved 

\begin{prop}
The  mapping $alg: BR(A, \P^2_0(\C)) \to \Hom_f(G^s_A, PO(3,\C))$ is 
such that $\al \circ alg = id$. The space $BR(A, \P^2_0(\C))$ lies in 
the image of the mapping $\al$. 
\end{prop}

Let $Res_T: \Hom_f(G^s_A, PO(3,\C)) \to \Hom_f(G^s_T, PO(3,\C))$ 
be the restriction homomorphism. Define the varieties 
$$
\BHom_f(G^s_A, PO(3,\C)):= Res_T^{-1} \{\rho_{\phi_T} , \om( \rho_{\phi_T})\} 
$$
$$
\BHom_f^+(G^s_A, PO(3,\C)):= Res_T^{-1} \{\rho_{\phi_T}\}
$$
of {\em based representations}. Clearly $\BHom_f^+(G^s_A, PO(3,\C))$ is 
the image $alg(BR(A, \P^2_0(\C)))$ (as a set) and the mapping 
$$alg: BR(A, \P^2_0(\C)) \to \BHom_f^+(G^s_A, PO(3,\C))$$
 is a bijection. 

\begin{lem}
\label{fina}
The representation $\rho= \rho_{\phi_T}: G_T^s \to PO(3,\R)\subset 
PO(3,\C)$ corresponding 
to the canonical realization $\phi= \phi_T$ of the standard triangle 
$T$ has finite image.  The centralizer of the group $\rho(G^s_T)$ in 
$PO(3,\C)$ is trivial.  
\end{lem}

\proof It is clear that the group $\rho(G^s_T)$ has invariant  
finite set 
$$\Sigma = \{(1,1), (-1, 1), (-1, -1), (1, -1)\}\subset 
\P^2(\R)\subset \P^2(\C)\quad .$$
 Any three distinct points of $\Sigma$ do not lie on a 
projective line in 
$\P^2(\C)$, thus if $g\in PO(3, \C)$ fixes $\Sigma$ pointwise 
then $g= id$. This implies finiteness of $\rho(G^s_T)$. It is 
easy to check 
that $\rho(G^s_T)$ equals $A_4$, the alternating group of 
the 4-element set $\Sigma$. Any element of $g\in PO(3, \C)$ 
centralizing $\rho(G^s_T)$ must 
fix $\Si$ pointwise, hence $g=\1 $. $\qed$

\begin{prop}
\label{lorigid}
The group $PO(3,\C)$ acts 
simply-transitively by conjugations on the set $\Hom_f(G^s_T, PO(3,\C))$. 
\end{prop}
\proof 
Suppose that $\rho\in \Hom_f(G^s_T, PO(3,\C))$. We note that $G^s= G^s_T$ contains the abelian subgroup $\Z/2 \times \Z/2 \times \Z/2$ generated by 
the involutions $g_{v_{00}}, g_{v_{x}}, g_{v_{y}}$. 
Since $\rho$ is nondegenerate 
we conclude that  the restriction of $\rho$ to this subgroup is 
injective. Thus, by the classification of finite subgroups of 
$PO(3, \C)$, we conclude that $\rho$ can be conjugate to a 
representation (which we again denote by $\rho$) so that 
the projective arrangement 
$\psi= \al(\rho)$ has the property:
$$
\psi(v_{00})= (0,0), \psi(v_{x})= (0, \8), \psi(v_{y})= (\8 , 0)
$$
Thus necessarily: $\psi(l_x)= L_x, \psi(l_y)= L_y, \psi(l_{\8})= L_{\8}$. 
Let $G_{10}$ denote the subgroup of $G^s$ generated by 
$g_{v_{00}}, g_{v_{y}}= g_{l_{x}}, g_{v_{10}}$. The restriction of 
$\rho$ to $G_{10}$ factors through the finite Coxeter group 
$G_{10}^c= G_{10}/\<\< (g_{v_{00}} g_{v_{10}})^2\>\>$.  
There are only two homomorphisms \newline
$\rho^c: G_{10}^c \to PO(3, \C)$ 
such that 
$$
\rho^c(g_{v_{00}})= \rho(g_{v_{00}}), \rho^c(g_{v_{y}})= \rho(g_{v_{y}}), 
\rho^c(g_{v_{10}})\ne  \rho^c(g_{v_{00}}) , \rho^c(g_{v_{10}})\ne  
\rho^c(g_{v_{y}})
$$
Namely, for one of them the isolated fixed point of $\rho^c(g_{v_{10}})$
has the affine coordinates $(1,0)$, for the second it has the affine 
coordinates $(-1,0)$. The reflection in the line $L_y$ conjugates 
one representation to the other and commutes with the elements 
$\rho(g_{v_{00}}), \rho(g_{v_{x}}), \rho(g_{v_{y}})$. Thus, after 
adjusting $\rho$ by this conjugation (if necessary), we conclude that 
\newline 
$\al(\rho)= \psi:  v_{10} \mapsto (1, 0)$. Similar argument works for the vertex $v_{01}$. 
It remains to determine $\psi(l_{x1}), \psi(l_{y1}), \psi(v_{11})$ and $\psi(l_d)$. Since $\rho(g_{l_{x1}})$ commutes with the elements $\rho(g_{v_{01}})$ and $\rho(g_{v_{\8 0}})$ and does not coinside with either, the line $\psi(l_{x1})$ contains the points $(0,1)$ and $(\8, 0)$. Similarly the line $\psi(l_{y1})$ contains $(0,1)$ and $(0, \8)$. 

We next determine $\psi(v_{11})$. Since the edge $[v_{11}, l_{y1}]$ has label 4, the elements $\rho(g_{v_{11}})$ and $\rho(g_{l_{y1}})$ anticommute, so $\psi(v_{11}) \in L_{y1}$ (the line in $\P^2$ joining $(1,0)$ and $(0,\8)$. Similarly $\psi(v_{11})\in L_{x1}$. Finally we determine $\psi(l_d)$. Since $\rho(g_{l_d})$ commutes with $\rho(g_{v_{00}})$ we have $(0,0)\in \psi(l_d)$. Also $\rho(g_{l_d})$ anticommutes with $\rho(g_{v_{11}})$, so $(1,1) \in \psi(l_d)$. We conclude that
$$
\psi= \al(\rho)= \phi_T
$$
and either $\rho= \rho_{\phi_T}$ or  $\rho= \om (\rho_{\phi_T})$. 
It now follows from Lemma \ref{fina} that the action of the group 
$PO(3,\C)$ on $Hom_f(G^s_T, PO(3,\C))$ is free. $\qed$

\begin{cor}
$\Hom_f(G^s_T, PO(3,\C))$ equals the $Ad PO(3,\C)$-orbit of the 
set\newline
  $\{\rho_{\phi_T}, \om (\rho_{\phi_T})\}$. 
In particular, both $\rho_{\phi_T}$ and $\om (\rho_{\phi_T})$ are 
locally rigid. 
\end{cor}

\noindent Note however that this corollary doesn't apriori imply 
that the variety $\Hom_f(G^s_T, PO(3,\C))$ is smooth since it 
could be nonreduced. To prove smoothness we need the following:

\begin{prop}
\label{C11.3}
The representations $\rho_{\phi_T}, \om (\rho_{\phi_T}): 
G_T^s \to PO(3,\C)$ are  
infinitesimally rigid. 
\end{prop}

We will prove the proposition only for $\rho_{\phi_T}$, the 
second case easily follows. Proposition \ref{C11.3} will 
immediately follow from the more general 

\begin{prop}
\label{rigid triangle}
Let $G^a= G_T^a$, define $\t\rho: G^a \to PO(3,\C)$ by composing 
$\rho_{\phi_T}$ with the canonical projection  $G_T^a\to G^s_T$. Then 
the representation $\t\rho$ is infinitesimally rigid.  
\end{prop}
\proof Our proof is based on the results of Section \ref{elem}. The reader 
will notice that the proof follows the lines of the proof of 
Proposition \ref{lorigid}. We first consider the subgroup 
$F$ in $G^a$ generated by $g_{v_{00}}, g_{l_x}, g_{l_y}$, these 
generators mutually commute, hence the subgroup is abelian. 
Let $\si\in Z^1(G^a, ad\circ\t\rho)$. By Lemma \ref{L5.3} the 
restriction of $\si$ 
to each cyclic subgroup $\<g_{v_{00}}\>, \<g_{l_x}\>, \<g_{l_y}\>$ 
is cohomologically trivial, thus $\si|_F$ comes from a cocycle on 
the finite Coxeter group
$$
F/\<\<\ g_{v_{00}}^2\ , \ g_{l_x}^2 \ ,\  g_{l_y}^2\ \>\>  
$$
which implies that $\si|_F$ is a coboundary. By adjusting the 
cocycle $\si$ by a coboundary we may assume that $\si|_F =0$. 
Now we consider the subgroup $F_{10}$ generated by 
$g_{l_x}, g_{v_{10}}$. Because of the Artin relations in 
$G^a_T$ this is again an abelian subgroup whose image under 
$\t\rho$ is dihedral. Hence $\si|_{F_{10}}$ is also a 
coboundary. Let $H_{10}$ denote the subgroup  generated by elements of 
$F_{10}$ and $g_{v_{00}}$, we recall that 
$(g_{v_{00}} g_{v_{10}})^2= (g_{v_{10}}g_{v_{00}} )^2$. 

Since the restriction of $\si$ to each generator of 
$H_{10}$ is exact, the cocycle $\si|_{H_{10}}$ comes 
from the finite Coxeter group 
$H_{10}/\<\< g_{v_{00}}^2, g_{l_x}^2, g_{v_{10}}^2\>\>$. 
This implies that $\si|_{H_{10}}$ is a 
coboundary. Since $\si(g_{l_x})=0, \si( g_{v_{00}})=0$ 
and the fixed points of $\t\rho(g_{l_x} ), \t\rho(g_{v_{00}} )$ 
are distinct we conclude  that $\si|_{H_{10}}=0$. 
The same argument implies that $\si(g_{v_{01}})=0$. We 
repeat our argument for the two abelian subgroups 
generated by $g_{l_x}, g_{v_{10}}, g_{l_{y1}}$ and 
by $g_{l_x}, g_{v_{01}}, g_{l_{x1}}$, it follows 
that $\si( g_{l_{y1}})= 0, \si( g_{l_{x1}})= 0$. 

Then we use Lemma \ref{L5.10}, where $a_1= g_{l_{x1}}, 
a_2= g_{l_{y1}}, a_3= 
g_{v_{00}}, b= g_{v_{11}}$, to conclude that 
$\si(g_{v_{11}})=0$. Finally $\si(g_{l_{d}})=0$ 
since the Shephard subgroup of $G^s$ generated by 
$g_{v_{11}},  g_{l_{d}}, g_{v_{00}}$ is again finite. $\qed$ 

\begin{cor}
\label{cross-sec}
The variety $\Hom_f(G^s_T, PO(3,\C))$ is smooth. The variety 
\newline 
$\BHom_f(G^s_T, PO(3,\C))$ is a (scheme-theoretic) cross-section for 
the action of $PO(3,\C)$ by conjugation  on the variety 
$\Hom_f(G^s_T, PO(3,\C))$. 
\end{cor}
\proof Smoothness of $\Hom_f(G^s_T, PO(3,\C))$ follows from 
infinitesimal rigidity of the representations 
$\rho_{\phi_T}, \om \rho_{\phi_T}$, see \cite[Theorem 2.6]{LM}. 
In Proposition 
\ref{lorigid} we proved that the morphism 
$$
\la: PO(3,\C)\times  \BHom_f(G^s_T, PO(3,\C)) \hook 
\Hom_f(G^s_T, PO(3,\C))
$$ 
given by the action of $PO(3,\C)$ by conjugation, is a bijection. 
Thus $\Hom_f(G^s_T, PO(3,\C))$ is also smooth, which implies 
that the morphism $\la$ is actually an isomorphism of varieties.   
$\qed$

\begin{cor}
\label{cro-section}
The variety $\BHom_f(G^s_A, PO(3,\C))$ is a (scheme-theoretic) 
cross-section for the action 
of $PO(3,\C)$ by conjugation  on the variety $\Hom_f(G^s_A, PO(3,\C))$. 
\end{cor}
\proof Consider the $Ad (PO(3,\C))$-equivariant restriction morphism
$$
Res_T: \Hom_f(G^s_A, PO(3,\C)) \lra \Hom_f(G^s_T, PO(3,\C))
$$
It was proven in Corollary \ref{cross-sec} that the subvariety
 $\BHom_f(G^s_T, PO(3,\C))$ is a cross-section for the action of $PO(3,\C)$ 
by conjugation on $\Hom_f(G^s_T, PO(3,\C))$. Thus the pull-back variety
$$
Res^{-1}_T \Hom_f(G^s_T, PO(3,\C)) \equiv  \BHom_f(G^s_A, PO(3,\C))
$$
is a cross-section as well. $\qed$

\medskip
Our goal is to show that the mapping 
$alg: BR(A, \P^2_0(\C)) \to \BHom_f^+(G^s_A, PO(3,\C))$ is 
an isomorphism of varieties over $\Q$, this will be proven in the next section. 

\subsection{$alg$ is an isomorphism of varieties}

We first establish that $alg: BR(A, \P^2_0(\C)) \to 
\BHom_f^+(G^s_A, PO(3,\C))$ is an isomorphism of varieties in 
two elementary cases. Let $C$ be an arrangement whose graph $\Ga_C$ 
has only one edge $e= [v, l]$, $m$ is the number of isolated vertices $v_j$ 
in $\Ga_C$. Let $A$ be a based arrangement which is the disjoint 
union of the standard triangle and $C$.

\begin{lem}
\label{above}
$alg: BR(A, \P^2_0(\C)) \to \BHom_f^+(G^s_A, PO(3,\C))$ is 
an isomorphism of varieties over $\Q$. 
\end{lem}
\proof We already know that $alg$ is a bijection. 
It is clear that the restriction morphisms
$$
Res: \BHom_f^+(G^s_A, PO(3,\C)) \to \Hom_f(G^s_C, PO(3,\C))$$
$$
res: BR(A, \P^2(\C)) \to BR(C, \P^2(\C)) $$
are isomorphisms of varieties. Let 
$alg: BR(C, \P^2(\C)) \to \Hom_f(G^s_C, PO(3,\C))$ 
denote the mapping induced by the restriction. The group 
$G^s_C$ is the free product 

$$
G^s_e * \underbrace{\Z/2 * ... * \Z/2}_{m\hbox{~times}}
$$
Hence 
$$
\Hom_f(G^s_C, PO(3))\cong \Hom_f(G^s_e, PO(3)) \times 
( \Hom_f(\Z /2, PO(3)) )^m
$$
\noindent Since the edge group $G^s_e$ and the vertex groups $\Z/2$ are  
finite, the variety $\Hom_f(G^s_C, PO(3))$ is smooth 
(Proposition \ref{L5.4}). The quasi-projective variety  
$R(C, \P^2_0)$ again splits as the direct product
$$
{\cal I}_0 \times  \P^2_0 \times ... \times  \P^2_0 \ , 
$$
\no where ${\cal I}_0$ is the anisotropic incidence variety, see 
 Section \ref{Projective arrangements}. The anisotropic incidence variety 
is smooth by  Proposition \ref{is smooth}, hence the product smooth as well. 

 Let $B$ denote the arrangement 
obtained by removing the incidence relation between $v$ and $l$ in $C$ 
and $G_B^s$ be the corresponding Coxeter group. Then 
$$
alg : R(B, \P^2_0) \lra \Hom_f(G_B, PO(3)) 
$$
is an isomorphism of smooth varieties (the left-hand side is 
$(\P^2_0)^{m+2}$ and the right-hand side is $p(R)^{m+2}$, 
see Lemma \ref{fixed point}). Thus 
$alg:  BR(A, \P^2_0)\lra \BHom_f(G, PO(3))$ is an 
isomorphism of smooth varieties. $\qed$ 

\medskip
We now consider a relative version of the above lemma. Suppose that 
a based abstract arrangement $D$ is the fiber sum $T\times_{w\equiv v} C$ 
where $C$ is the arrangement above, $w, v$ are elements of 
$T$ and $C$ respectively. 
 
\begin{lem}
\label{below}
$alg: BR(D, \P^2_0(\C)) \to \BHom_f^+(G^s_D, PO(3,\C))$ is 
an isomorphism of varieties. 
\end{lem}
\proof 
Let $\al = \rho_{\phi_T}(g_v)$,  $G_e\subset G^s_D$ is the  edge subgroup 
corresponding to $e=[v, l]$. 
It is clear that the restriction morphism 
$$
\BHom^+_f(G^s_A, PO(3,\C))\stackrel{Res}{\RA} 
\{\rho\in \Hom_f(G_C^s, PO(3,\C)): \rho(g_v)= \al \}$$
$$ 
\quad \quad \quad =: F_{g_v, \al}(G_C^s, PO(3, \C))$$
is an isomorphism of varieties. The group $G_e$ is finite, hence by 
Proposition \ref{L5.4} the variety $F_{g_v, \al}(G_C^s, PO(3, \C))$ is smooth. 
Similarly let $q:= \phi_T(v)$, the restriction morphism
$$
res: BR(D, \P^2_0(\C)) \to \{\psi \in R(C, \P^2_0(\C)): \psi(v)= q\} =: 
F_{v, q}(C, \P^2_0(\C))
$$
is an isomorphism of varieties. Then $F_{v, q}(C, \P^2_0(\C))$ is 
isomorphic the product of the {\em relative anisotropic incidence variety} 
$$
{\cal I}_0(q)= \{ l \in (\P^2_0)^{\vee} : q \cdot l =0 \}
$$
(see \S \ref{Projective  arrangements}) with $m$ copies of 
$\P^2_0$. It is clear that the mapping 
$alg$ induces a bijection of the sets of complex points 
$$
alg: {\cal I}_0(q) (\C) \lra F_{g_v, \al}(G_e, PO(3, \C))
$$
According to Lemma \ref{is smooth} the variety ${\cal I}(q)$ is 
smooth and we  repeat the arguments from the proof of  
Lemma \ref{above}. $\qed$ 

\bigskip
Now we consider the case when $A$ is a general based arrangement. 
Let $X:= BR(A, \P^2_0)$, $Y:= \BHom_f^+(G_A^s, PO(3))$. 

\begin{thm}
\label{T7.2}
The mapping $alg: X\to Y$ is a 
biregular  isomorphism of quasi-projective varieties. 
\end{thm}
\proof Let $F$ (resp. $R$) be the functor of points of $X$ (resp. $Y$). 
Then $alg$ is an isomorphism of the varieties $X$ and $Y$ if and only if 
$F$ and $R$ are naturally isomorphic (see \cite[Proposition IV-2]{EH}). Let 
$B$ be the based arrangement obtained by removing all edges from the 
graph of $A$ that are not in the graph of $T$ (and retaining all vertices). 
Let $G^s_B$ be the corresponding Shephard group, clearly 
$$
G^s_B\cong G^s_T * \underbrace{\Z/2 * ... * \Z/2}_{m \hbox{~~times}}    
$$ 
where $m$ is the cardinality of ${\cal V}(A- T)$. Thus  
$$BR(B, \P^2_0(\C))\cong \P^2_0(\C)^m, \BHom^+_f(G^s_B, PO(3))\cong 
\Hom_f(\Z/2, PO(3, \C))^m$$
 are smooth varieties and 
$$
\widetilde{alg} : \t{X}:= BR(B)\lra \t{Y}:= \BHom^+_f(G^s_B, PO(3))
$$
is an isomorphism of smooth varieties (see the proof of  Lemma \ref{above}). 
We let $\t{F}$ and 
$\t{R}$ be the  functors of points of $\t{X}$ and $\t{Y}$. The isomorphism 
$\widetilde{alg}: \t{X}\to \t{Y}$ induces a natural isomorphism of functors 
$\eta: \t{F}\to \t{R}$. The functors $F$ and $R$ are subfunctors of 
$\t{F}$ and $\t{R}$. We now make explicit the 
inclusions $F\subset \t{F}, R\subset \t{R}$. 

Let ${\cal E}_0$ be the collection of edges of $\Ga_A$ 
that are not edges of $\Ga_T$.  
Suppose $e_{ij}= [v_i, \ell_j]\in {\cal E}_0$ is such an edge. 
Let $X_{ij}$ be the subvariety of $BR_0(B)$ defined by 
$$
x_i \al_j + y_i \be_j + z_i \ga_j=0
$$
where $[x_i: y_i: z_i]$ and $[\al_j: \be_j :\ga_j]$ are 
homogeneous coordinates on $BR_0(B)$ corresponding to 
$v_i, l_j$ respectively. 
We let $F_{ij}$ be the subfunctor of $\t{F}$ which as the functor of points 
associated to the subvariety $X_{ij}$. Then we have
\begin{equation}
%\label{fuF}
F = \bigcap_{[v_i l_j]\in \E_0} F_{ij}  \subset\t{F}
\end{equation}
Similarly if the edge $e_{ij}$ has the label $\eps_{ij}$ then we let 
$R_{ij}$ be the subfunctor of $\t{R}$ corresponding to the subvariety 
$Y_{ij}$ defined by the Artin relation 
$$
(g_{v_i} g_{\ell_j})^{\eps_{ij}}= (g_{\ell_j}g_{v_i} )^{\eps_{ij}}\ .$$ 
We have  
\begin{equation}
%\label{fuR}
R= \bigcap_{(ij)\in \E} R_{ij} \subset \t{R}
\end{equation}

\begin{lem}
$\eta: \t{F}\to \t{R}$ induces an isomorphism from $F_{ij}$ to $R_{ij}$. 
\end{lem}
\proof Let $X_{ij}$ and $Y_{ij}$ be the subvarieties  of 
$BR(B, \P^2_0), \BHom_f^+(G^s_B, PO(3))$ corresponding to the 
subfunctors $F_{ij},  R_{ij}$. Then Lemmas \ref{above} 
and \ref{below} imply that $alg$ induces 
an isomorphism of smooth varieties $X_{ij}\to Y_{ij}$. 
Hence $\eta$ induces an isomorphism  of the corresponding 
subfunctors. $\qed$

The above lemma and equations (2), (3) immediately imply that   
$\eta$ induces a natural isomorphism from $F$ to $R$. Theorem 
\ref{T7.2} follows.  $\qed$

\medskip
Now we can prove one of the two main results of this paper. Let 
$S \subset \C^n$ be an affine variety defined over 
$\Q$. We will consider $S$ as a quasi-projective variety 
in $\P^n(\C)$.  

\begin{thm}
\label{main}
For any variety $S$ as above there exists a Zariski open 
subset $U\subset S(\C)$ containing all real points and a  
based arrangement 
$A$ so that the corresponding  Shephard  group $G^s_A$ has the property: 

There is a Zariski open subset $W$ in $\Hom(G^s_A, PO(3, \C))\h PO(3, \C)$ 
which is biregular isomorphic to $U$.  
\end{thm}

\begin{rem}
$W$ is never Zariski dense in the character variety $X(G^s_A, PO(3,\C))$. 
\end{rem}

\proof Given the variety $S$ we construct an abstract based 
arrangement $A$ such that $BR_0(A, \P^2)$ is biregular isomorphic 
to $S$ via the isomorphism $geo$ (Theorem \ref{systems} and Corollary \ref{csystems}). Let 
$U:= \tau (BR_0(A, \P^2_0(\C)))$, it is  Zariski open in 
$S(\C)$ and contains all real points   
since the subvariety $BR_0(A, \P^2_0)$ is Zariski open in 
$BR_0(A, \P^2)$ and contains all real points.  Theorem \ref{T7.2} implies that we have 
an isomorphic embedding with 
Zariski open image $\BHom_f^+(G^s_A, PO(3, \C))_0$
$$
alg: BR_0(A, \P^2_0) \hook  \BHom_f^+(G^s_A, PO(3, \C))
$$
Corollary \ref{cro-section} implies that $\BHom_f^+(G^s_A, PO(3))$ 
is a cross-section for the action of $PO(3,\C)$ by conjugation on 
the Zariski component $\Hom_f^+(G^s_A, PO(3))$ of 
$\Hom(G^s_A, PO(3))$. All representations 
in the variety $\Hom_f(G^s_A, PO(3, \C))$ are {\em stable} in the 
sense of \cite[ Theorem 1.1]{JM}, thus we get an open 
monomorphism of varieties
$$
\BHom_f^+(G^s_A, PO(3,\C))_0 \hook X(G^s_A, PO(3,\C))= 
\Hom(G^s_A, PO(3,\C))\h PO(3,\C)
$$
so that the image is a Zariski open subvariety $W$ in the 
character variety. Thus we get the composition 
$$
\theta: U \stackrel{geo}{\lra} BR_0(A, \P^2(\C)) 
\stackrel{alg}{\longrightarrow} 
\BHom_f^+(G^s_A, PO(3,\C))_0 \lra W \subset  X(G^s_A, PO(3,\C)) 
$$
is the required isomorphism onto a Zariski open subvariety $W$ of the 
character variety. $\qed$

\begin{prop}
\label{sfina}
Suppose that  $S$ is an affine variety over $\Q$ and $q\in \Q$ is a rational point. Then there is an abstract arrangement $A$ (as in Theorem \ref{main}) 
so that the representation 
$\rho= alg\circ geo(q)$ has finite image. The centralizer of the subgroup
$\rho(G^s_A)\subset PO(3,\C)$ is trivial. 
\end{prop}
\proof Follows from Corollary \ref{csystems} and Lemma \ref{fina}. $\qed$

\medskip
As an example we consider the configuration space of the 
arrangement\footnote{ 
This is not a based arrangement.} 
$$A= \{v; l_1, l_2: \i(v, l_1), \i(v, l_2)\} $$
Take the corresponding Coxeter group $G= G_A^c$. Then 
$R(A, \P^2_0)\cong \P^2_0 \times \P^1_0 \times \P^1_0$ corresponds to 
the representations of $G$ which are {\em nondegenerate}, 
i.e. the elements $\rho(g_v), \rho(g_{l_1}), \rho(g_{l_2})$, 
$\rho(g_v g_{l_1})$,  $\rho(g_v g_{l_2})$ have order 2. However 
there are some other components of \newline 
$\Hom(G, PO(3,\C))$ which are 
described by assigning which of the elements\newline   
$\rho(g_v), \rho(g_{l_1}), \rho(g_{l_2}), \rho(g_v g_{l_1})$, 
$\rho(g_v g_{l_2})$ are equal to $\1$. If $\rho(g_v)= \1$ then 
any representation factors through the free product 
$\Z/2 * \Z/2$ and the corresponding component of 
$\Hom(G, PO(3,\C))$ is isomorphic to $\P^2_0 \times \P^2_0$. 
The reader will verify that besides 
$$
\Hom_f(G, PO(3,\C))\cong \P^2_0 \times \P^1_0 \times \P^1_0$$ 
and the above component isomorphic to  
$\P^2_0 \times \P^2_0$ there are 4 components isomorphic to 
$\P^2_0\times \P^1_0$, 4 components isomorphic to $\P^2_0$ and 
one component which consists of a single reduced point 
(the last corresponds to the trivial representation).

\subsection{Representations of Artin groups}
\label{repart}

Our next goal is to prove a theorem analogous to Theorem \ref{main} 
for Artin groups. Take a based admissible arrangement $A$, consider 
the germ $(\Hom(G^a_A, PO(3, \C)), \t\rho)$, where $\t\rho= \mu(\rho)$, 
$\rho\in \Hom_f(G^s_A, PO(3,\C))$.  

\begin{thm}
\label{iso}
The morphism $\mu: (\Hom(G^s_A, PO(3,\C)), \rho) 
\to (\Hom(G^a_A, PO(3)), \t\rho)$ is an analytical  
isomorphism of germs. 
\end{thm}
\proof Let $m\ge 1$ be an integer, consider the $m$-th 
order Zariski tangent spaces $T^m_{\rho}$ of the variety 
$\Hom(G^s_A, PO(3,\C))$ at $\rho$ and $T^m_{\t\rho}$ 
of the variety $\Hom(G^a_A, PO(3,\C))$ at $\t\rho$, there 
is a well-defined mapping $\mu_m: T^m_{\rho} \to T^m_{\t\rho}$ 
induced by $\mu$. We will think of elements of $m$-th order 
Zariski tangent spaces as {\em formal curves} or homomorphisms 
from $G^s, G^a$ to 
$PO(3)({\cal A})$ which is the set of  ${\cal A}$-points of $PO(3)$ 
and ${\cal A}$ is a certain Artin local $\C$-algebra, see \S \ref{defof}. 

Our goal is to prove that $\mu_m$ is an isomorphism for all $m$. 

\begin{lem}
$\mu_m$ is injective.
\end{lem}
Suppose that $\al\ne \be\in T^m_{\t\rho}$ are formal 
curves of the order $m$. We think of them 
as representations of $G^s$ into $PO(3)({\cal A})$. 
Then $\mu_m(\al), \mu_m(\be)$ are homomorphisms of 
$G^a$ into $PO(3)({\cal A})$ defined by composing 
with the projection $G^a\to G^s$. Clearly $\mu_m(\al)\ne \mu_m(\be)$. 
$\qed$ 

\begin{lem}
$\mu_m$ is surjective.
\end{lem}
\proof Let $\xi\in T^m_{\t\rho}(G^a, PO(3,\C))$ and $v$ is a 
{\em point} in $A- \{v_{11}\}$.  
To prove that $\xi$ belongs to the image of $\mu_m$ it is enough to check that 
the restriction of $\mu_m$ to each vertex-subgroup of $G^a$ is a {\em trivial deformation} (see Lemma \ref{L7.1}).  We first consider the {\em points} 
$v\in A$ distinct from $v_{11}$. Then $v$ 
is incident to at least one {\em line} $l\in A$ and hence  the image of the 
edge-subgroup $\rho( G_{vl}^s)= \t\rho(G^a_{vl})$ is generated by two 
distinct  commuting involutions (since $\rho$ is ``nondegenerate'').  
Thus by Lemma \ref{L5.3} the restriction $\rho|G_{vl}^a$ is 
infinitesimally rigid, in particular $\xi|G^a_v$ is a 
{\em trivial deformation}. In the case $v=v_{11}$ we use   
Proposition \ref{rigid triangle} to conclude that the restriction of 
$\xi$ to the vertex-subgroup 
$G_v^a$ is a trivial deformation as well. Suppose that 
$l\in A$ is a {\em line}. Then admissibility of the arrangement 
$A$ implies that $l$ is incident to at least one {\em point} 
$w$ in $A- \{v_{11}\}$. Hence 
we repeat the same argument as in the case of {\it points} in $A$. $\qed$

Thus, we established that the morphism $\mu$ induces an 
isomorphism of Zariski tangent spaces of all orders. Hence  
Theorem \ref{iso} follows from Lemma \ref{LL3.1}. $\qed$ 

\begin{cor}
\label{inhe}
For any admissible based arrangement $A$ the morphism  
$$\mu: \Hom_f(G_A^s, PO(3))\to \Hom(G_A^a, PO(3))$$
is open (in the classical topology) 
and is an analytic isomorphism onto its image. 
\end{cor}
 
\begin{cor}
For any admissible based arrangement $A$ the character variety of 
$G_A^a$ inherits all 
singularities of the representation variety of the 
group $G_A^s$ corresponding to points of $BR(A, \P^2_0(\C))$. 
\end{cor}

The image $\Hom_f(G^a_A, PO(3,\C)):= \mu(\Hom_f(G_A^s, PO(3,\C)))$ 
is a constructible set. 

\begin{prop}
$\Hom_f(G^a_A, PO(3,\C))$ is closed in $\Hom(G^a_A, PO(3,\C))$ 
\end{prop}
\proof The image $Z$ of the monomorphism
$$
\Hom(G_A^s, PO(3,\C)) \lra \Hom(G_A^a, PO(3,\C))
$$
is described by the equations
$$
\rho(g_v)^{\del(v)}= \1 , \quad v\in \V(\La_A)
$$
where $\del(v)$ is the label of the vertex $v$. Thus $Z$ is  
closed. On the other hand, by Corollary \ref{zaropen}, 
the space $\Hom_f(G^s_A, PO(3,\C))$ is closed in 
$\Hom(G_A^s, PO(3,\C))$, the proposition follows. $\qed$

\begin{cor}
\label{open-artin}
For any admissible based arrangement the space 
$\Hom_f(G^a_A, PO(3,\C))$ is a union of Zariski connected components of 
$\Hom(G^a_A, PO(3, \C))$.
\end{cor}
\proof The set $\Hom_f(G^a_A, PO(3,\C))$ is Zariski open 
since it is constructible and open in the classical 
topology. $\qed$

\begin{thm}
The morphism $\mu: \Hom_f(G_A^s, PO(3))\to  \Hom_f(G_A^a, PO(3)))$ 
is a biregular isomorphism of varieties. 
\end{thm}
\proof We will first prove that the  reduced  varieties corresponding to 
$$
X:= \Hom_f(G_A^s, PO(3)) \hbox{~~and~~} Y:= \Hom_f(G_A^a, PO(3))$$
 are isomorphic. We construct nonsingular varieties $U, W$ containing $X, Y$ 
and an extension of $\mu$ to an isomorphism $U\to W$ which 
carries $X$ to $Y$ bijectively as sets. Let $C$ be the 
arrangement obtained from $A$ by removing the incidence 
relations everywhere outside of the standard triangle. 
Thus we get two nonsingular varieties $Q= \Hom_f(G^a_C, PO(3))$ 
and $U:= \Hom_f(G^s_C, PO(3))$ with $X\subset U$ and $Y\subset Q$. 

Clearly $\mu: U \to V$ is a biregular isomorphism and it 
bijectively carries $X$ to $Y$. Thus the reduced 
varieties $X^r, Y^r$ are isomorphic via $\mu$. Hence the assertion 
of Theorem follows from combination of Corollary  
\ref{inhe} and Theorem \ref{cbasic}. $\qed$ 

\medskip
As a corollary we get the following generalization of Theorem \ref{main}:  

\begin{thm}
\label{Artin}
Let $S$ be an affine variety defined over $\Q$. Then  
there exists an admissible based arrangement 
$A$ so that the corresponding  Artin group $G^a_A$ has the property: 

There is a Zariski open subset in $S$ (containing 
all real points) which is biregular isomorphic to a Zariski open and 
closed subvariety in $\Hom(G^a_A, PO(3,\C))\h PO(3,\C)$. If 
$q\in \Q$ is a rational point. Then the abstract arrangement 
$A$ can be chosen so that the representation 
$\rho= \mu\circ alg\circ geo(q)$ has finite image. The 
centralizer of the subgroup $\rho(G^a_A)\subset PO(3,\C)$ is trivial. 
\end{thm}

\section{Differential graded algebras and Lie algebras}
\label{dgla}

In this section we discuss differential graded Lie algebras, differential algebras  and their Sullivan's minimal models. These definitions will be used in the following two sections.  
Let $\k$ be the ground field. A {\em graded Lie algebra} over $\k$ is a $\k$-vector space
$$
L^{\bul}= \oplus_{i\ge 0}\ L^i
$$
graded by (nonnegative) integers and a family of bilinear mappings
$$
[\cdot , \cdot]: L^i \otimes L^j \lra L^{i+j}
$$ 
satisfying graded skew-commutativity:
$$
[\al, \be] + (-1)^{ij}[\be, \al] =0
$$
and the graded Jacobi identity:
$$
(-1)^{ki}[\al, [\be, \ga]] + (-1)^{ij}[\be, [\ga, \al]]
+ (-1)^{jk}[\ga, [\al, \be]]=0
$$
where $\al\in L^i , \be\in L^j , \ga\in L^k$. We say that $L$ is 
{\em bigraded}  if $L$ is graded as above and 
$$
L= \oplus_{q\ge 0} L_q= \oplus_{i,q\ge 0} \ L^i_q 
$$
where $ L^i_q= L^i\cap L_q$. We also require:
$$
[\cdot , \cdot]: L_q \otimes L_p \lra L_{p+q}
$$

A {\em differential graded Lie algebra} is a pair $(L, d)$ where 
$d: L \to L$ is a {\em derivation} of the degree $\ell$, i.e.
$$
d: L^i \lra L^{i+\ell}, \quad d\circ d=0, \quad \hbox{and}
$$
$$
d[\al, \be]= [d\al, \be] + (-1)^{i\ell}[\al, d\be]
$$
\begin{rem}
In this paper we will use only the degree $\ell=1$. 
\end{rem}

Suppose that $L= L^{\bul}_{\bul}$ is a bigraded Lie algebra, then the \dgla 
$(L, d)$ is called a  {\em differential bigraded Lie algebra} if 
there is a number $s$ such that $d$ has the bidegree $(\ell, s)$:
$$
d: L_q \lra L_{q+s}
$$
\begin{rem}
In this paper we will use only the bidegree $(\ell,s)=(1,0)$. Similarly one defines {\em trigraded} differential Lie algebras $L^{\bul}_{\bul,\bul}$
\end{rem}

Let $L^{\bul}$ be a differential graded Lie algebra, and suppose that  
${\cal J}\subset L^{\bul}$ is a vector subspace. Then $\J$ 
is an {\em ideal} in $L^{\bul}$ if:
\begin{itemize}
\item $\J$ is {\em graded}, i.e. $\J= \oplus_{n=0}^{\8} 
(\J \cap L^n)$; 
\item $d(\J) \subset \J$; 
\item For each $\ga\in \J, \al\in L$ we have: $[\ga, \al]\in \J$.   
\end{itemize}   

\begin{rem}
If $L$ is bigraded then we require ideals in $L$ to be bigraded 
as well, i.e.
$$
\J= \oplus_{n,q\ge 0}\ (\J \cap L^n_q)
$$ 
\end{rem}

\begin{lem}
\label{big}
If $\J \subset L$ is an ideal in a differential graded (bigraded) 
Lie algebra then the quotient $L/\J$ has the natural structure 
of a differential graded (bigraded) Lie algebra so that the 
natural projection $L\to L/\J$ is a morphism. 
\end{lem}
\proof Proof is straightforward and is left to the reader. $\qed$ 

\medskip
Let $(L^{\bul},d)$ be a \dgla and ${\mathfrak g}$ be a Lie algebra. 
An {\em augmentation} is a homomorphism $\eps: L \to {\mathfrak g}$ 
such that $\oplus_{i\ge 1} L^i \subset \ker(\eps)$ and $\eps \ne 0$. 
An {\em augmentation ideal} of $\eps$ is the kernel of $\eps$. 

We recall the definition of the complete local $\k$-algebra 
$R_{L^{\bullet}}$ (see \cite{M}) associated to a \dgla 
$L^{\bullet}$ over $\k$. We will give a definition which 
conveys the intuitive meaning of $R_{L^{\bullet}}$ but has 
the defect of being non-functorial. We have the sequence  
$$
0\lra L^0 \stackrel{d}{\lra} L^1 \stackrel{d}{\lra} 
L^2 \stackrel{d}{\lra} ... 
$$
Choose a complement $C^1\subset L^1$ to the 1-coboundaries 
$B^1\subset L^1$. 
Let ${\cal F}: C^1 \to L^2$ be the polynomial mapping given by 
${\cal F}(\eta)= d\eta + \half [\eta, \eta]$. Then $R_{L^{\bullet}}$ 
is the completion at $0$ of the coordinate ring of the affine 
subvariety of $C^1$ defined by the equation ${\cal F}(\eta)= 0$. 

Since $L^1$ and $C^1$ are infinite dimensional in many applications, 
further definitions of affine variety and completion are needed in 
infinite dimensional vector spaces. These are given in \cite[\S 1]{BuM}. 
Suppose next that $(X, x)$ is an analytic germ. We will say that 
$L^{\bullet}$ {\em controls} $(X,x)$ if $R_{L^{\bullet}}$ is 
isomorphic to the complete local ring $\widehat{O_{X,x}}$. 

\begin{defn}
Let $A^{\bul}, B^{\bul}$ be cochain complexes and $\rho: A^{\bul}\to
B^{\bul}$ is a morphism. Then $\rho$ is called a {\bf quasi-isomorphism} 
if it induces an isomorphism of all cohomology groups. The morphism is 
called a {\bf weak equivalence} if it induces an isomorphism of 
$\H^0, \H^1$ and a monomorphism of $\H^2$.   
\end{defn}

Weak equivalence induces an equivalence relation on the category of differential graded Lie algebras: algebras $A^{\bul}_1, A^{\bul}_m$ 
are (weakly) equivalent if there is a chain of weak equivalences:
$$
A^{\bul}_1 \rightarrow A^{\bul}_2 \leftarrow A^{\bul}_2 \rightarrow ... 
\leftarrow A^{\bul}_m
$$

In the following two sections we will use the following theorem about controlling \dglas proven in \cite{GM}: 

\begin{thm}
\label{GolMil}
Suppose that $L^{\bul}$ and $N^{\bul}$ are weakly equivalent 
\dglas. Then the germs controlled by the algebras 
$L^{\bul}$ and $N^{\bul}$ are analytically isomorphic. 
\end{thm}

A {\em differential graded algebra} $A^{\bul}$ 
is defined similarly to a \dgla: just instead of the Lie bracket 
$[\cdot , \cdot]$ satisfying graded Jacobi identity we have an 
associative multiplication:
$$
\wedge : A^i \otimes A^j \lra A^{i+j}
$$
satisfying the properties:

\begin{itemize}
\item $d: A^i \to A^{i+1}$, $d\circ d=0$; 
\item $d(\al \wedge \be)= d\al \wedge \be + (-1)^i \al \wedge d\be$, for all 
$\al\in A^i$; 
\item $\al\wedge\be= (-1)^{ij}\be \wedge \al$, for all 
$\al\in A^i$, $\be \in A^j$;
\item $A$ has the unit $1\in A^0$. 
\end{itemize}

To get a \dgla from a differential graded algebra $A^{\bul}$ 
take a Lie algebra ${\mathfrak g}$ and let 
$L^{\bul}= A^{\bul} \otimes {\mathfrak g}$ (see \cite{GM} for details). 

Suppose that $V$ is a vector space over $\k$, $(A^{\bul}, d)$ is a 
differential graded algebra and $f: V \to Z^{2}$ is a linear mapping, where $Z^2$ is the space of 2-cocycles of $A^{\bul}$. 
Then the {\em Hirsch extension} $A^{\bul}\otimes_f V$ is a 
differential graded algebra which (as an algebra) equals to 
$A^{\bul}\otimes \La(V)$ and the restriction of the 
differential on $A^{\bul}\otimes_f V$ to $A^{\bul}$ 
equals $d$ and the restriction to $V$ equals $f$. 

A differential graded algebra $\M^{\bul}$ is called {\em 1-minimal} if:

(a) $\M^0 = \k$;

(b) $\M^{\bul}$ is the increasing union of differential subalgebras:
$$
\k = \M_{[0]} \subset \M_{[1]} \subset \M_{[2]}\subset ... 
$$ 
with each $\M_{[i]} \subset \M_{[i+1]}$ a Hirsch extension;

(c) The differential $d$ on $\M$ is decomposable, i.e. for each $\al\in \M$ we have: 
$$d(\al)= \sum_{j,i} \be_j \wedge\ga_i\ , \ \ 
\be_j\ , \ga_i \in \oplus_{s\ge 1}\M^s$$ 

\begin{defn}
Suppose that $A^{\bul}$ is a differential graded algebra. Then a 
{\bf 1-minimal model} for $A^{\bul}$ is a morphism 
$\rho: \M^{\bul} \to A^{\bul}$ such that:
\begin{itemize}
\item The differential graded algebra $\M^{\bul}$ is 1-minimal; 
\item The morphism $\rho$ is a weak equivalence. 
\end{itemize}
\end{defn}

We refer the reader to \cite{Sul}, \cite{GrM}, \cite{Morgan2} 
for further discussion of the definition, properties and 
construction of 1-minimal models.

\section{Hain's theorem and its applications}
\label{Hain's theorem}

In this section we give an exposition of a work of 
R.\ Hain \cite{Hain} which shows that the singularities 
in representation varieties of fundamental  groups of smooth 
complex algebraic varieties 
are quasi-homogeneous. In fact some assumption on the 
representation $\rho$ which is being deformed is  also 
required. In \cite{Hain} the analogue of our Theorem 
\ref{hain} is proven under the assumption that $\rho$ 
was the monodromy representation of an admissible 
variation of a mixed Hodge structure. One doesn't obtain a 
restriction on weights working in this generality. 

Let $M$ be a smooth  connected  
manifold with the fundamental group $\Ga$. Let 
$\G$ be the Lie group of real points of an algebraic group $G$ with the 
Lie algebra $\gg$ defined over $\R$ and $\rho: \Ga \to \G$ is a 
homomorphism. Let $P$ be the flat bundle over $M$ associated to $\rho$ 
and $adP$ the associated $\gg$-bundle. Then $\a^{\bul}(M, adP)$, the 
complex of smooth $adP$-valued differential forms on $M$ is a \dgla. We 
define an augmentation $\eps: 
\a^{\bul}(M, adP)\to \gg$ by evaluating degree zero forms 
at a base-point $x\in M$ and sending the rest of forms to zero. 
Let $\a^{\bul}(M, adP)_0$ be the kernel of $\eps$. The following theorem 
follows immediately from \cite[Theorem 6.8]{GM}. 

\begin{thm}
\label{con germ}
$\a^{\bul}(M, adP)_0$ controls the germ $(\Hom(\Ga, \G),\rho)$. 
\end{thm}

The point of this section is that if $M$ is smooth  connected  complex 
algebraic variety and $\rho$ has finite image then 
$\a^{\bul}(M, adP)_0$ is quasi-isomorphic to a \dgla which  
has a structure  of a mixed Hodge complex. 
By a theorem of Hain this implies that 
$$
R_{\a^{\bul}(M, adP)_0\otimes \C}
$$
is a \qh ring.  We now give details. 

A {\em real mixed Hodge complex} (abbreviated MHC) is a 
pair of complexes $K_{\R}^{\bul}$ and $K_{\C}^{\bul}$ (real 
and complex respectively), 
together with a quasi-isomorphism $\al: K_{\R}^{\bul} 
\otimes \C \to K_{\C}^{\bul}$ such that  $K_{\R}^{\bul}$ 
is a complex of real vector spaces equipped with an 
increasing filtration $W_{\bul}$ (called the {\em weight} 
filtration) and $K_{\C}^{\bul}$ is equipped with an increasing 
(weight) filtration $W_{\bul}$ and a decreasing filtration $F^{\bul}$ 
( called the {\em Hodge} filtration). The data 
$$
K_{\R}^{\bul}, K_{\C}^{\bul}, \al, W_{\bul}, F^{\bul}
$$
satisfies the axioms described in \cite[Scholie 8.1.5]{D2} (take $A= \R)$. 

By a theorem of Deligne (\cite[ Scholie 8.1.9]{D2}) 
the cohomology of a MHC has a 
mixed Hodge structure, \cite[Scholie 2.3.1]{D}. It 
is important in what follows that the filtrations on 
$\H^ {\bul}(K_{\C})$ induced by $W_{\bul}$ and $F^{\bul}$  
can be canonically split, \cite[Section 1.2.8]{D}.  
Thus the filtration  $W_{\bul}$ induces a canonical grading  
 on $\H^ {\ell}(K_{\C})$ and consequently on $\H^ {\ell}(K_{\C})^{*}$, 
$\ell= 0, 1, 2$. 

If $V$ is a finite-dimensional vector space over 
$\C$ we will let $\C\, [[V]]$ denote the completion 
of the symmetric algebra $\C\, [V]$ at the maximal ideal 
corresponding to $0$. Thus $\C\, [[\H^ {\ell}(K_{\C})^*]]$ has 
a canonical decreasing filtration (as an algebra) induced by the grading
 of $\H^ {\ell}(K_{\C})^{*}$, $\ell= 0, 1, 2$.

We will say that a MHC is a {\em mixed 
Hodge differential graded Lie algebra} if the 
complexes  $K_{\R}^{\bul}$ and $K_{\C}^{\bul}$ 
are \dglas (\cite[\S 1.1]{GM}), such that $\al$ is 
bracket preserving and the filtrations satisfy

(i) $[W_p (K_{\R}^{\bul}), W_q (K_{\R}^{\bul})]\subset W_{p+q} 
(K_{\R}^{\bul})$;

(ii) $[F^p (K_{\C}^{\bul}), F^q (K_{\C}^{\bul})]\subset F^{p+q} 
(K_{\C}^{\bul})$.

\smallskip
In this case we will use $L^{\bul}_{\R}$ and  
$L^{\bul}_{\C}$ in place of $K_{\R}^{\bul}$ 
and $K_{\C}^{\bul}$. Let ${\mathfrak m}$ denote the maximal ideal of 
$\C\, [[\H^ 1( L^{\bul}_{\C})^*]]$. We now have

\begin{thm}
\label{ha}
(Hain's Theorem.) Suppose $L^{\bul}= (L^{\bul}_{\R}, 
L^{\bul}_{\C}, \al, W_{\bul}, F^{\bul})$ is a mixed 
Hodge \dgla with $\H^ 0(L^{\bul}_{\R} )= 0$. Then $R_{L^{\bul}_{\C}}$ 
 is \qh (see Section \ref{singu}). 
Moreover there exists a morphism of graded vector spaces 
$$
\del: \H^ 2(L^{\bul}_{\C})^* \to Gr^W  \C\, [[ \H^ 1(L^{\bul}_{\C})^*]]
$$
with image of $\del$ contained in ${\mathfrak m}^2$  
such that $R_{L^{\bul}_{\C}}$  is the quotient 
of $\C\, [[ \H^ 1(L^{\bul}_{\C}) ^* ]]$ by the graded ideal 
generated by the image of $\del$. 
\end{thm}

\begin{rem}
Hain further proves that $\del\equiv Q^* \ (mod ~~{\mathfrak m}^3)$, 
$Q^*$ is dual to $Q$ where
$$
Q: S^2 \H^ 1(L^{\bul})\to \H^ 2(L^{\bul})
$$
is given by the cup-product $Q(\eta)= [\eta, \eta]$.  
\end{rem}

We will say an element $\eta\in \H^ i(L^{\bul}_{\C} )$ has 
{\em weight} $n$ if 
$\eta\in W_n \H^ i( L^{\bul}_{\C})$ but 
$\eta\notin W_{n-1} \H^ i( L^{\bul}_{\C})$. We will combine Hain's 
theorem with the following theorem to obtain our desired result 
about the singularities in representation varieties of fundamental 
groups of smooth complex algebraic varieties.  Suppose now that 
$M$ is a smooth  connected  complex algebraic variety, a representation 
$\rho: \pi_1(M)\to \G$ with finite image, the bundle $adP$, etc., are 
as above.

\begin{thm}
\label{we}
Under the conditions above there is a filtration $W_{\bul}$ on 
$A^{\bul}(M, ad P)$ and a filtration $F^{\bul}$ on  
$A^{\bul}(M, ad P_{\C})$ 
such that for the canonical map \newline 
$\al: A^{\bul}(M, ad P)\otimes \C \to 
A^{\bul}(M, ad P_{\C})$ the algebra 
$$
L^{\bul} = ( A^{\bul}(M, ad P), A^{\bul}(M, ad P_{\C}), 
\al, W_{\bul}, F^{\bul})$$
is a mixed Hodge \dgla. Moreover the weights of  
$\H^ 1( A^{\bul}(M, ad P_{\C}) )$ are $1$ and $2$
and the weights of $\H^ 2( A^{\bul}(M, ad P_{\C}) )$ 
are $2, 3$ and $4$. 
\end{thm}
\proof Let $\t{M}$ be the finite cover of $M$ corresponding 
to $\ker(\rho)$. Let $\Phi\cong \rho(\Ga)$ be the group of 
covering transformations. 
By \cite{Su} there exists an equivariant completion $\t{N}$ 
of  $\t{M}$. But according to \cite{BM} there is also a 
canonical resolution of singularities $N$ of $\t{N}$ so that the 
complement $N- \t{M}$ is a divisor $D= D_1 \cup ... \cup D_m$ with normal crossings.  
 
Hence the action of $\Phi$ extends to $N$, which is a 
smooth $\Phi$-equivariant completion  of $\t{M}$. 
Hence $\Phi$ acts on the log-complex of $\t{M}$ 
defined using the compactification   $N$. It is a 
basic result of Deligne (\cite[Theorem 1.5]{D}) that 
one may use the log-complex to define:

(a) A subcomplex $A^{\bul}(\t{M})\subset \a^{\bul}(\t{M})$ 
so that the inclusion is a quasi-isomorphism. 

(b) Filtrations $W_{\bul}, F^{\bul}$ on $A^{\bul}(\t{M})$ and 
$A^{\bul}(\t{M})\otimes \C$ which satisfy the axioms 
of MHC. 

By construction these filtrations are 
$\Phi$-invariant.  We tensor with $\gg$ (regarded 
as a mixed Hodge \dgla concentrated in degree zero). 
Since $\Phi$ acts on 
$\gg$ via $ad \rho$, it also acts on the tensor 
product. We obtain the required mixed Hodge 
\dgla $L^{\bul}$ by taking $\Phi$-invariants. 
To derive weight restrictions we use results of
 Morgan \cite{Morgan1}, \cite{Morgan2}, who proved it 
for  $A^{\bul}(\t{M})\otimes \C$, etc. The operations of 
tensoring with $\gg$ and taking $\Phi$-invariants 
will not change  these restrictions on weights. $\qed$

\begin{rem}
Choose a point $m\in X$ and a point $\t{m}\in \t{M}$ over 
$m$. We define augmentations (see \cite[\S 3.1]{GM}) 
$\eps: A^{\bullet}(M, ad P) \to \gg$ and 
$\t\eps: A^{\bullet}(\t{M})\otimes \gg\to \gg$ by 
evaluation at $m$ and $\t{m}$ respectively. We let 
$A^{\bullet}(M, ad P)_0$ and $A^{\bullet}(\t{M})_0\otimes \gg$ 
be the augmentation ideals. Then all statements in Theorem \ref{we} 
hold when   $A^{\bullet}(M, ad P)$ and $A^{\bullet}(M, ad P_{\C})$  
are replaced by $A^{\bullet}(M, ad P)_0$ and 
$A^{\bullet}(M, ad P_{\C})_0$. We abbreviate the corresponding 
mixed Hodge \dglas by $L^{\bullet}_0$. 
\end{rem}

%\begin{rem}
%Note that a group $\Ga$ is the fundamental group of a smooth connected 
%complex quasi-projective variety if and only if it 
%is the fundamental group of a smooth complex quasi-projective surface. 
%Namely, let $X\subset \C\P^k$ be a smooth quasi-projective variety, 
%$\dim(X)\ge 3$.  
%Take the intersection $Y$ of $X$ with a generic hyperplane. Then $Y$ is 
%smooth and the generalized Lefshetz hyperplane section 
%theorem \cite{GoMa} implies that 
%$\pi_1(Y)\cong\pi_1(X)$. By repeating this argument  we reduce the 
%dimension of $X$ to $2$. However it is unclear if the same conclusion holds 
%for smooth open complex algebraic varieties 
%(which are not assumed to be quasi-projective).    
%\end{rem}

Let $\Ga, \rho, G$ be as above. Let $Z$ be 
the representation variety $\Hom(\Ga, G)$. 
%Assume $\H^ 0(\Ga, \gg_{ad})= 0$.  
By combining Theorems \ref{ha}, \ref{we} we obtain

\begin{thm}
\label{hain}
The germ $(Z_{\C}, \rho)$ is analytically equivalent to a 
quasi-homogeneous cone with generators of weights 
$1$ and $2$ and relations of weights $2, 3$ and $4$. 
Suppose that there is a  local 
cross-section $S$ through $\rho$ to $G$-orbits. Then the  conclusion 
is valid  not just for the germ $(Z_{\C}, \rho)$ but 
also for $(\Hom(\Ga, G)\h G, [\rho])$. 
\end{thm}
\proof We apply Theorems \ref{ha} and \ref{we} to deduce 
that the complete local $\C$-algebra $(R_{L^{\bul}_0})_{\C}$ 
has a presentation of the required type. But by  
Theorem \ref{con germ}, $(R_{L^{\bul}_0})_{\C}$ 
 is isomorphic to the 
complete local ring associated to the germ $(Z_{\C}, \rho)$. 
We obtain the corresponding result for 
$(\Hom(\Ga, G_{\C})\h G_{\C}, [\rho])$ by 
replacing $L^{\bul}_0$ by 
$L^{\bul}$ and applying \cite[Theorem 2.4 ]{KM3}.  Note that if $S$ 
exists then $\H^ 0(L^{\bul})=0$. $\qed$

\medskip
There are infinitely many germs $(Y, 0)$, where 
$Y\subset \C^n$ is an affine variety defined over 
$\Z$, so that   $(Y, 0)$ is not quasi-homogeneous 
with the weights of relations between $2$ and $4$.
We can even assume that $0$ is an isolated singular point, see 
the Section \ref{singu}. 
Thus as a consequence of Theorems \ref{hain},  
\ref{main} we obtain the following 

\begin{thm}
\label{T10.3}
Among the Artin groups $G_A^a$ there are infinitely 
many mutually nonisomorphic groups which are not 
isomorphic to fundamental groups of smooth complex 
algebraic varieties. 
\end{thm}
\proof Let $Y$ be an affine variety defined over $\Q$ and $y\in Y$ 
be a rational point. Assume that the analytical germ $(Y, x)$ is 
not \qh (with the weights of variables $1,2$ and weights of generators $2,3,4$). Let $A$ be an affine arrangement corresponding to the pair 
$(Y, y)$ as in Proposition \ref{sfina}, so that the representation 
$\rho_s: G^s_A\to PO(3,\C)$ corresponding to $y$ has finite image and 
the group $\rho_s(G^s_A)$ has trivial centralizer in $PO(3,\C)$. Let 
$\rho= \mu(\rho_s): G^a_A\to PO(3,\C)$, where $G_A^a$ is the {\em Artin 
group} of the arrangement $A$. Recall that we have an open embedding 
$$
\mu \circ alg \circ geo: Y\cong BR_0(A)\hook \Hom(G_A^a, PO(3))\h PO(3) = X(G_A^a, PO(3))
$$
Suppose that $G_A^a$ is 
the fundamental group of a smooth complex algebraic  
variety.  Then Theorem \ref{hain} can be applied to 
the germ $(X(G_A^a, PO(3)), [\rho])$ 
provided we can construct a local cross-section through $\rho$ to the $PO(3,\C)$-orbits. In the definition of local cross-section we take 
$$
U:= \Hom_f(G_A^a, PO(3,\C)) \hbox{~~and~~} 
S:= \mu \BHom_f(G_A^s, PO(3,\C))\ .
$$
 Then $U$ is open by Corollary \ref{open-artin} and $S$ is a 
cross-section because $\BHom_f(G_A^s, PO(3,\C))$ is a 
cross-section for the action of $PO(3,\C)$ on 
$\Hom_f(G_A^s, PO(3,\C))$ and the morphism 
$\mu: \Hom_f(G_A^s, PO(3,\C)) \to U$ is an isomorphism. We get a 
contradiction. To see that there 
are infinitely many nonisomorphic examples we refer to 
the argument at the end of the introduction.  $\qed$ 

\begin{figure}[tbh]
\leavevmode
\centerline{\epsfxsize=5.5in \epsfbox{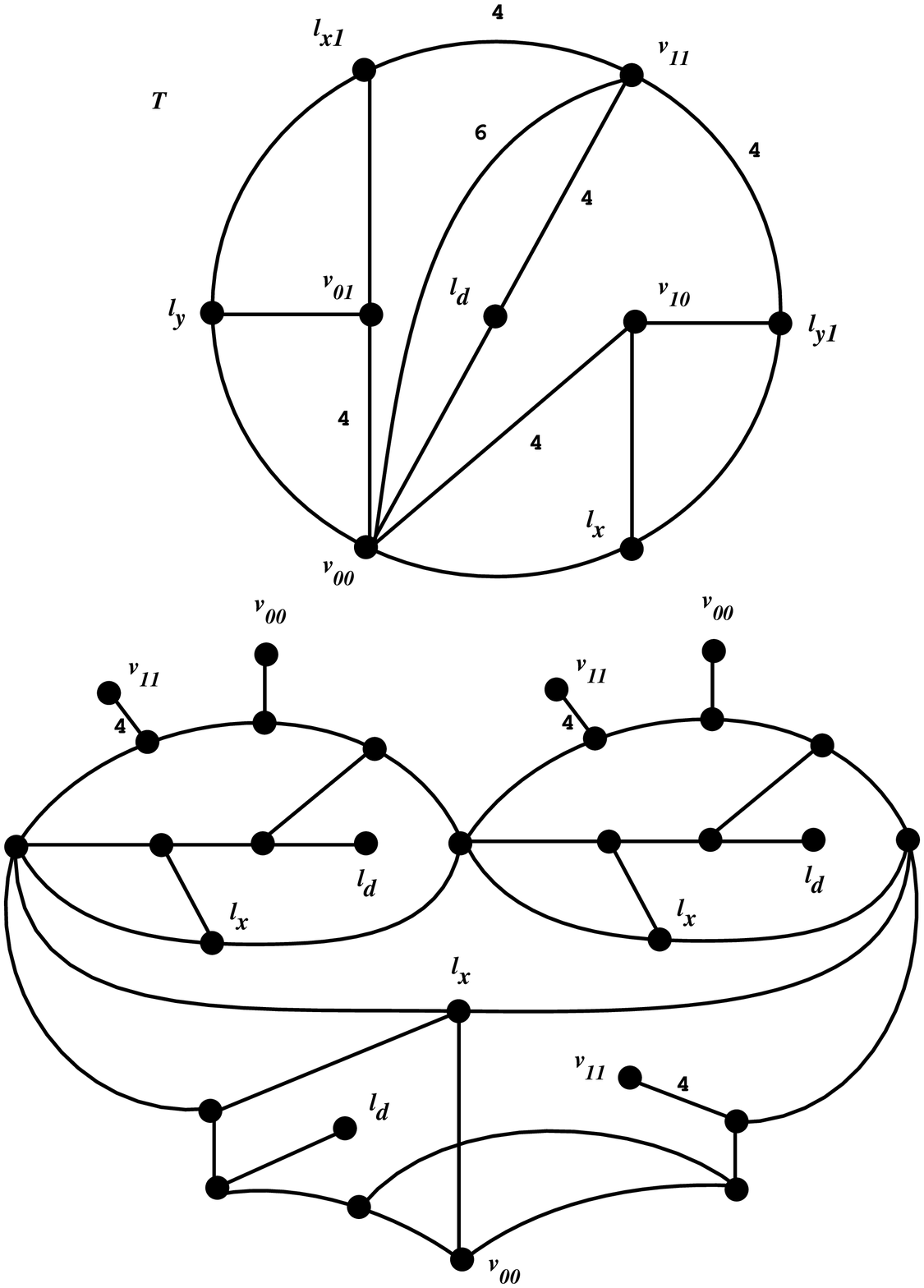}}
\caption{\sl Labelled graph of an Artin group.}
\label{Fig9}
\end{figure}

\medskip
As the simplest example of $(Y,0)$ we can take the germ 
$(\{x^5=0\}, 0)$.   We describe the Coxeter graph of the 
Artin group corresponding to this singularity on the 
Figure \ref{Fig9}. We let $x^5= (x^2)^2 \cdot x$, to get the 
labelled graph\footnote{See \S\S \ref{groups}, \ref{gaa}.} 
$\La$ of $G^a$  from this diagram identify vertices marked 
by the same symbols.  

\section{Sullivan's minimal models and singularities of 
representation varieties}
\label{direct}

The goal of this section is to give a direct proof of Theorem \ref{hai} below
(which is Theorem \ref{t7} of the Introduction). 

\begin{thm}
\label{hai}
Let $M$ be a smooth connected 
complex algebraic variety with the fundamental group $\Gamma$. Let $\G$ be 
the Lie group of real points of an algebraic group $G$ defined over 
${\mathbb R}$~; let $\gg$ be the  Lie algebra of $\G$. 
 Suppose that $\rho: \Ga\to \G$ is a representation 
with finite image. Then the germ  $(\Hom(\Gamma, G), \rho)$
 is analytically isomorphic to a 
quasi-homogeneous cone with generators of weights $1$ and $2$ and 
relations of weights $2, 3$ and $4$. 
In the case there is a local 
cross-section through $\rho$ to $Ad(G)$-orbits, then 
the same conclusion is valid for the quotient germ 
$(X(\Ga, G), [\rho])$ of the character variety.   
\end{thm}
\proof We begin the proof by choosing a smooth $\Phi$-equivariant compactification $\t{N}= \t{M}\cup D$ of $\t{M}$ as in the proof of Theorem \ref{we}.

In what follows we shall use the following simple lemma: 

\begin{lem}
\label{fi}
Suppose that $\Phi$ is a finite group and 
$$
\begin{array}{ccccc}
~ & ~ & H\\
~ & ~ & g {\Big\downarrow}\\
E & \stackrel{f}{\longrightarrow} & F 
\end{array}
$$
is a diagram of morphisms of $\Phi$-modules over $\C$ such that 
$f(E)\subset g(H)$. Then $f$ admits a $\Phi$-equivariant lifting 
$\t{f}: E \to H$. 
\end{lem}
\proof Since $f(E)\subset g(H)$ there exists a linear mapping $h: E \to H$ which lifts $f$. Then we let
$$
\t{f}= Av(h):= |\Phi|^{-1}\sum_{\phi\in\Phi} \phi \circ h \circ \phi^{-1} \quad \qed
$$ 

Recall that Morgan in \cite{Morgan2} defines a {\em mixed Hodge diagram}  
$$
\E(log D) \stackrel{\varphi}{\longleftarrow}  E_{C^{\8}}(\t{M}) \otimes \C \stackrel{\bar\varphi}{\lra} \bar{\E}(log D)
$$
associated with the pair $(\t{N}, D)$. The log-complex $\E(log D)$ is a subcomplex of $\C$-valued differential forms on $\t{M}$ and the complex $\bar{\E}(log D)$ is the log-complex with the {\em opposite} complex structure. 

The mixed Hodge diagram must satisfy 
certain properties described in \cite{Morgan2}. 
In particular, it has a structure 
of a mixed Hodge complex, i.e.  $E_{C^{\8}}(\t{M})$ has an increasing 
filtration $W$ and $\E(log D)$ has a pair of filtrations: an increasing 
weight filtration $W$ and decreasing Hodge filtration $F$; $\varphi, \bar\varphi$ must preserve the weight filtrations and be quasi-isomorphisms. 

\begin{prop}
\label{mhd}
There exists a $\Phi$-invariant mixed Hodge diagram with 
$\Phi$-invariant structure of a mixed Hodge complex. The identity embedding
$$
id: \E(log D) \hook \a(\t{M})\otimes\C 
$$
is a quasi-isomorphism. 
\end{prop}
\proof First we describe  the {\em log-complex} $\E(log D)$ on 
$\t{M}$ associated with the 
compactification $N$. Let $z\in D$ be a point of $p$-fold intersection
$$
z\in D_{i_1} \cap ...\cap D_{i_p}
$$
where each $D_{i_j}$ is locally (near $z$) is given by the equation 
$z_{i_j}=0$. Then elements $\sigma$ of $\E$ are $\C$-valued differential forms 
on $\t{M}$ which can be (locally with respect to the $z_j$-coordinates)
written as
$$
\sum_J \eta_J \frac{dz_{i_1}}{z_{i_1}} \wedge ... \wedge \frac{dz_{i_p}}{z_{i_p}}
$$
where each $\eta_j$ extends to a $C^{\8}$-form in a \nbd of $z$ in $N$. 
Thus near any generic point $z\in D_j$ the form $\sigma$ has at worst a simple 
pole. Since the group $\Phi$ acts holomorphically on $N$ leaving $D$ 
invariant we conclude that this group acts naturally on the log-complex 
$\E(log D)$. The complex $\E(log D)$ has the {\em weight} and {\em Hodge} filtrations which are defined in a canonical way. Recall that 
$W_{\ell}(\E(log D))$ consists of differential forms of the type
$$
\om = \sum_J \om_J \wedge \frac{dz_{j_1}}{z_{j_1}}\wedge ... \wedge \frac{dz_{j_t}}{z_{j_t}} \ ,\ t\le \ell
$$
and $F^p(\E(log D))$ consists of differential forms of the type
$$
\om = \sum_J \om_J \wedge dz_{j_1}\wedge ... dz_{j_s} \wedge  \frac{dz_{j_{s+1}}}{z_{j_{s+1}}}\wedge ... \wedge \frac{dz_{j_t}}{z_{j_t}} \ , \  t\ge p 
$$
where the forms $\om_J$ extend smoothly over the divisor $D$. Thus both filtrations are $\Phi$-invariant. 

Now we describe the second complex  
$E_{C^{\8}}(\t{M})$ associated to $(N, D)$. For each component $D_j$ of $D$ choose a regular \nbd $N_j$ in $\Phi$-invariant way, i.e. if $\phi\in \Phi$ and 
$\phi: D_i \to D_j$ then $\phi: N_i \to N_j$. Let 
$[\om_j]\in \H^2(N_j, \D N_j; \R)$ be the Thom's class. By using Lemma \ref{fi} 
we choose a 2-form $\om_j$ representing the class $[\om_j]$ so that for $\phi\in \Phi$ mapping $D_i$ to $D_j$ we have: $\phi^*\om_j =\om_i$. Then Morgan takes 1-forms $\gamma_j \in \E(log D)$ supported on $N_j$ so that 
$d\ga_j= \om_j$. We can choose $\ga_j$ so that $\phi^*\ga_j =\ga_i$ for each $\phi\in \Phi$ such that $\phi: D_i \to D_j$ by using Lemma \ref{fi} again. 

\medskip
The complex $E_{C^{\8}}(\t{M})$ consists of global 
sections of a certain sheaf ${\cal S}$ that we will describe 
below. Let $U\subset \t{M}$ be an open subset missing all regular 
\nbds $N_j$. Then sections of  ${\cal S}$ over $U$ are real-valued 
infinitely differentiable differential forms on $U$. 

Let $N^p$ denote subset of $\t{N}$ consisting of $p$-fold intersections between the regular \nbds $N_j$. Take a connected open subset $U\subset N^p$ which is disjoint from $N^{p+1}$.  We suppose that $U$ is contained in the $p$-fold intersection $N_{i_1} \cap ... \cap N_{i_p}$. Then sections of 
${\cal S}$ over $U$ are elements of the Hirsch extension:  
$$
\a^{\bul}(U) \otimes_d \La(\tau_{i_1},..., \tau_{i_p})
$$
where $\a^{\bul}(U)$ is the complex of real-valued differential 
forms on $U$ and 
$$
d\tau_{i_j}= \om_{i_j}|_U
$$ 
It is clear that the group $\Phi$ acts on the sheaf ${\cal S}$ and on the 
complex $E_{C^{\8}}(\t{M})$ of its sections as well. The complex 
 $E_{C^{\8}}(\t{M})$ has a canonically defined {\em weight} 
filtration $W_{\bul}$ which is therefore invariant under the action of 
$\Phi$. (This filtration is similar to the weight filtration on 
the log-complex, just use $\tau_j$-s instead of the forms $dz_j/z_j$.) 
Finally Morgan defines a morphism
$$
\al: E_{C^{\8}} (\t{M}) \otimes \C \to \E(log D)
$$
by mapping $\tau_j$ to $\ga_j$ and differential forms supported in $\t{M}$ 
to themselves. Morgan proves that this morphism and the identify embedding 
$$
\E(log D) \hook \a^{\bul}(\t{M})\otimes \C 
$$
induce isomorphisms of cohomology groups. The mixed 
Hodge structure on the diagram is $\Phi$-invariant by the construction. 
This finishes the proof of Proposition \ref{mhd}. $\qed$

\begin{rem}
In what follows we shall use the notation  $A^{\bul}= A^{\bul}(\t{M})$ to denote the complex $E_{C^{\8}}(\t{M})\otimes \C$. 
\end{rem}

Given the weight filtration $W$ on $A^{\bul}$ Morgan defines an increasing 
{\em Dec-weight  filtration} $Dec\, W_{\bul}(A^{\bul})$ as
$$
Dec\, W_{\ell}(A^{k})= \{ x : x\in W_{k-\ell}\, (A^k), dx \in W_{k-\ell+1}\, (A^{k-1})\}
$$
Let $\nu: {\cal N}\to A^{\bul}$ denote a 1-minimal model:  
$$
{\cal N}= (\ 0\lra \N^0= \C \stackrel{0}{\lra} \N^1 \stackrel{d}{\lra} \N^2 \stackrel{d}{\lra} ...\ )
$$
Recall the basic properties of ${\cal N}$ and $\nu$ proven in 
\cite[\S6 and Lemma 7.2]{Morgan2}: 

(a) The $Dec$-weight filtration $Dec\, W_{\bul}(A^{\bul})$ of $A^{\bul}$ 
pulls back to $\N$ to a {\em weight} filtration $W_{\bul}(\N)$ which 
{\em splits} so that $\N$ becomes {\em bigraded} with the differential 
of the bidegree $(1,0)$:

$\N^j= \oplus_{i\ge 0} \N_i^j$, $\N^j_i= \N^j\cap \N_i$,
 where $\N_i$ consists of elements of the {\em weight} $i$,  
$$
d: \N_i \to \N_i, \ \ \wedge : 
\N_i\otimes \N_j \to \N_{i+j}\quad ,$$ 
$\N_0= \N^0$ and each $\N^i_j$ is finite-dimensional. 
 
(b) The weight filtration on $\N$ induces a weight filtration on 
$\H^{\bul}(\N)$ so that the induced weights on $\H^1(\N)$ 
are $1, 2$, the induced weights on $\H^2(\N)$ are $2, 3, 4$. 

(c) The homomorphism $\nu: \N\to A^{\bul}$ is a weak equivalence. 

(d) $\N$ is a 1-minimal differential algebra. Together with (a) it 
implies that the restriction of the differential 
$d$ to $\N^1_1$ is identically zero.

\begin{prop}
\label{lifmin}
Suppose that $A^{\bul}$ is a differential graded commutative algebra 
(over the ground field $\k= \C$ or $\R$) and 
$\Phi$ is a finite group acting on $A^{\bul}$. Then  
the action of $\Phi$ on $A^{\bul}$ lifts 
to its action on a certain 1-minimal model $\N$ for 
$A^{\bul}$. 
\end{prop}
\proof We will prove the proposition by constructing $\N$ in 
$\Phi$-invariant way. Let $\N^0:= \k$. Take $\N_{[1]}^1:= \H^1(A^{\bul})$ 
and $\N_{[1]}$ be the differential graded algebra freely generated by 
$\N^0= \N_{[1]}^0$  and $\N_{[1]}^1$. We need a homomorphism
$$
\nu= \nu_{[1]}: \N_{[1]} \to A^{\bul}
$$
which induces an isomorphism of the 1-st cohomology groups. The group $\Phi$ naturally 
acts on $\H^1(A^{\bul})$. We have the epimorphism of 
$\Phi$-modules 
$$
Z^1(A^{\bul}) \lra \H^1(A^{\bul})
$$
By Lemma \ref{fi} this epimorphism admits a $\Phi$-invariant splitting 
$\nu^1: \H^1(A^{\bul}) \to Z^1(A^{\bul})$. Thus we let $\nu$ be the identity 
embedding of $\N^0=\k$ to $A^0$ and $\nu|\N_{[1]}^1$ be $\nu^1$. 
We continue the construction of $(\N, \nu)$ by induction. Suppose 
that  $(\N_{[i]}, \nu_{[i]})$ 
are constructed and the homomorphism of $\Phi$-modules 
$\nu_{[i]}: \N_{[i]} \to A^{\bul}$ induces an isomorphism 
of $\H^1$ but the induced mapping of $\H^2$ has nonzero kernel. 
This kernel is canonically isomorphic to the relative cohomology group $\H^2(\N_{[i]}^{\bul}, A^{\bul})$. Choose a $\Phi$-equivariant section $\si_{[i]}$ to the projection from the space of relative cocycles 
$Z^2(\N_{[i]}^{\bul}, A^{\bul})$ onto $\H^2(\N_{[i]}, A^{\bul})$. Let
$$
p_1: Z^2(\N_{[i]}^{\bul}, A^{\bul}) \to Z^2(\N_{[i]}^{\bul})$$
and
$$
p_2: Z^2(\N_{[i]}^{\bul}, A^{\bul}) \to A^1$$
be the projections. Both are $\Phi$-equivariant. Define 
$$
d: \H^2(\N_{[i]}^{\bul}, A^{\bul}) \to Z^2(\N_{[i]})
$$
by $d= p_1 \circ \si_{[i]}$ and let
$$
\nu: \H^2(\N_{[i]}^{\bul}, A^{\bul}) \to A^1\ , \quad \nu= p_2 \circ \si_{[i]}
$$
Define
$$
\N_{[i+1]}:= \N_{[i]} \otimes_d \H^2(\N_{[i]}, A^{\bul})
$$
and extend $\nu$ to $\nu_{[i+1]}: \N_{[i+1]} \to A^{\bul}$ multiplicatively. 
The group $\Phi$ acts on $\N_{[i+1]}$ in the natural way and the homomorphism $\nu_{[i+1]}$ is $\Phi$-equivariant. 
$\qed$ 

\medskip
We now assume that we have a mixed Hodge diagram:
$$
{\cal E} \stackrel{\varphi}{\longleftarrow} A_{\C} 
\stackrel{\overline{\varphi}}{\lra} \overline{\cal E}
$$
and a finite group $\Phi$ acting on the diagram compatibly with the \dga structures so that $\varphi$ and $\overline{\varphi}$ are $\Phi$-equivariant. 
In \cite[\S 6]{Morgan2}, Morgan constructs a trigraded 1-minimal model 
\[
\begin{array}{ccccc}
 {\cal E} & \stackrel{\varphi}{\longleftarrow} & A_{\C} & 
\stackrel{\overline{\varphi}}{\lra} & \overline{\cal E} \\
~ &  \nwarrow & \uparrow & \nearrow \\
~ & ~ & \N^{\bul}_{\bul,\bul} & ~ & ~   \\
\end{array}
\]
for the above mixed Hodge diagram (see \cite[Page 270]{Morgan2} for definition). 
We will say that $\N^{\bul}_{\bul,\bul}$ is a $\Phi$-equivariant 
trigraded 1-minimal model for the above mixed Hodge diagram if all three morphisms with the source $\N$ are $\Phi$-equivariant and the trigrading of $\N$ is $\Phi$-invariant.

\begin{rem}
Morgan calls $\N^{\bul}_{\bul,\bul}$ a {\bf bigraded} minimal model. 
\end{rem}

\begin{prop}
\label{lif}
There exists a $\Phi$-equivariant trigraded 
1-minimal model for the mixed Hodge diagram 
$$
{\cal E} \stackrel{\varphi}{\longleftarrow} A_{\C} 
\stackrel{\overline{\varphi}}{\lra} \overline{\cal E}
$$
The filtration $Dec\, W_{\bul}(A)$ pulls back to a filtration $Dec\, W_{\bul}(\N)$ given by
$$
Dec\, W_{q}(\N):= \oplus_{r+s\le q}\ \N_{r,s}
$$
Consequently the filtration $Dec\, W_{\bul}(\N)$ is $\Phi$-equivariantly split. 
\end{prop}
\proof We will check that Morgan's construction can be made $\Phi$-equivariant. To do this we examine Morgan's induction step when he passes from a trigrading on $\N_{[i]}$ to one on $\N_{[i+1]}$. This step is carried on the page 176 of Morgan's paper and involves a study of the diagram
\[
\begin{array}{ccccc}
Z^2(\N_{[i]}, {\cal E}) & \stackrel{j}{\longleftarrow} & Z^2(\N_{[i]}, A_{\C}) & \stackrel{\overline{j}}{\lra} & Z^2(\N_{[i]}, \overline{\cal E}) \\
~ &  \searrow \! \nwarrow s & \downarrow \!\uparrow p &  s'\nearrow\! \swarrow 
\\~ & ~ & \H^2(\N_{[i]}, A_{\C}) & ~ & ~   \\
\end{array}
\]

By induction $\H^2(\N_{[i]}^{\bul})$ has a $\Phi$-equivariant mixed Hodge structure. Since $\H^2(\N_{[i]}^{\bul}, A_{\C}^{\bul})$ is the kernel of the canonical (thus $\Phi$-equivariant) morphism
$$
\H^2(\N_{[i]}^{\bul})\to \H^2(A_{\C}^{\bul})
$$
it inherits a $\Phi$-equivariant mixed Hodge structure. Consequently by Deligne's Theorem (see \cite[Proposition 1.9]{Morgan2}) 
$\H^2(\N_{[i]}^{\bul}, A_{\C}^{\bul})$ has a canonical 
(hence $\Phi$-invariant) bigrading. Thus it suffices to check that the cross-sections $s, p$ and $s'$ as well as the maps $h, h'$ on page 176 of \cite{Morgan2} can be chosen $\Phi$-equivariant. The cross-sections $s, p, s'$ are required to satisfy the linear conditions (1)--(3) of 
\cite[Page 176]{Morgan2}. It is immediate that our averages 
$Av(s), Av(p), Av(s')$ (defined as in Lemma \ref{fi} with respect 
to the action of $\Phi$) again satisfy (1)--(3). Finally the maps
$$
h: \H^2(\N_{[i]}^{\bul}, A_{\C}^{\bul}) \to Dec W_{r+s}({\cal E})
$$
and
$$
h': \H^2(\N_{[i]}^{\bul}, A_{\C}^{\bul}) \to Dec W_{r+s}(\overline{\cal E})
$$
must satisfy certain lifting conditions. By Lemma \ref{fi} we can take $h, h'$ to be $\Phi$-equivariant. 

By definition (see Proposition \ref{lifmin}) we have
$$
\N_{[i+1]}= \N_{[i]}\otimes_d \H^2(\N_{[i]}^{\bul}, A^{\bul}_{\C})
$$
we extend the trigrading from $\N_{[i]}$ and $\H^2(\N_{[i]}^{\bul}, A^{\bul}_{\C})$ to $\N_{[i+1]}$ multiplicatively. We obtain a $\Phi$-equivariant  trigrading on  $\N_{[i+1]}$ and a new diagram
\[
\begin{array}{ccccc}
 {\cal E} & \longleftarrow & A_{\C} & 
\lra & \overline{\cal E} \\
~ &  \nwarrow & \uparrow & \nearrow \\
~ & ~ & \N_{[i+1]} & ~ & ~   \\
\end{array}
\]
of equivariant maps satisfying Morgan's axioms. This completes the 
proof of the proposition. $\qed$ 

\medskip
We will no longer need the trigrading on $\N$ and instead will consider the bigrading:
$$
\N^k_{q}:= \oplus_{r+s=q}\ \N^k_{r,s}
$$
which defines a splitting of the Dec-filtration $Dec\, W_{\bul}(\N)$. 

\medskip

 Let $P$ be the flat bundle over $M$ associated 
to $\rho$ and $adP$ the associated $\gg$-bundle; similarly 
$ad\tilde{P}= \t{M}\times \gg$. We let $adP_{\C}$, $ad\t{P}_{\C}$  
denote the complexifications of these vector bundles. 
Let $\a^{\bul}(M, adP)$, $\a^{\bul}(\t{M}, ad\tilde{P})$ denote the 
complexes of $adP$, $ad\tilde{P}$--valued differential forms. 
According to  \cite[Theorem 2.4]{KM3} we have: 

\begin{thm}
\label{con germ1}
If there is a local cross-section through $\rho$ to the $Ad(G)$-orbits then the 
differential graded Lie algebra $\a^{\bul}(M, adP_{\C})$ {\bf controls} the 
germ $(X(\Gamma, G(\C)),[\rho])$. 
\end{thm} 

We tensor $A^{\bul}$ with the Lie algebra $\gg$ 
(regarded as a differential graded Lie algebra concentrated in the 
degree zero). Since 
$\Phi$ acts on $\gg$ via $ad \rho$, it also 
acts on the tensor product. 
Clearly $\a^{\bul}(M, adP_{\C})\cong 
\a^{\bul}(\t{M}, ad\tilde{P}\otimes\C)^{\Phi}$. 
Since the inclusion $A^{\bul}(\t{M})\hook \a^{\bul}(\t{M}, \C)$
induces an isomorphism of cohomology groups,  the inclusion
$$
A^{\bul}(\t{M})\otimes \gg\hook \a^{\bul}(\t{M}, adP_{\C}) $$
also induces an isomorphism of cohomology groups of the differential 
graded Lie algebras. Hence the inclusion
$$
(A^{\bul}(\t{M})\otimes \gg)^{\Phi}\hook 
(\a^{\bul}(\t{M}, ad\t{P}_{\C})^{\Phi}) \cong 
\a^{\bul}(M, adP_{\C})
$$
will also induce an isomorphism of cohomology groups. We conclude that
 under the conditions of Theorem  \ref{con germ1} the differential 
graded Lie algebra 
$(A^{\bul}(\t{M})\otimes \gg)^{\Phi}$ controls the germ of the 
character variety $(X(\Gamma, G),[\rho])$ (see Theorem \ref{GolMil}).

\medskip
We will need similar results for the representation variety 
$\Hom(\Ga, G)$ itself. Pick a point $m \in \tilde{M}$.  
We define an augmentation 
$\epsilon: \a^{\bul}(M, adP_{\C})\to \gg\otimes\C$ 
by evaluating degree zero forms 
at a base-point $m\in M$ and sending the rest of forms to zero. 
Let $A^{\bul}(M, adP_{\C})_0$ be the kernel of $\epsilon$. 
Recall that by Theorem \ref{con germ}
$\a^{\bul}(M, adP_{\C})_0$ controls the germ $(\Hom(\Gamma, G(\C)),\rho)$. 

\medskip
We lift the augmentation $\epsilon$ to 
$\a^{\bul}(\t{M}, ad\tilde{P}_{\C})$ as follows. 
Let $\t{m}$ be a point in $\t{M}$ which projects to $m$. Then for each 
$\om \in \a^{0}(\t{M}, ad\tilde{P}_{\C})$ let
$$
\t\epsilon(\om):= |\Phi|^{-1}\sum_{\ga\in \Phi} \om(\ga\t{m})
$$
where $|\Phi|$ is the order of the group $\Phi$. We extend $\t\eps$ to 
the rest 
of $\a^{\bul}(M, ad\tilde{P}_{\C})$ by zero. 
It is clear that the restriction of $\t\eps$ to 
$$
\a^{\bul}(\t{M}, ad\t{P}_{\C})^{\Phi}\cong \a^{\bul}(M, adP_{\C})$$
is the same as $\eps$. We let 
$A^{\bul}_0:= A^{\bul}\cap \ker(\t\eps)$. 
Again the inclusion
$$
(A^{\bul}_0\otimes \gg )^{\Phi}\hook 
\a^{\bul}(\t{M}, ad\t{P}_{\C})^{\Phi}_0
$$
is a quasi-isomorphism (i.e. it induces an isomorphism on cohomologies).

The action of the finite group $\Phi$ lifts from $A^{\bul}\otimes\gg$  
to the tensor product $\N \otimes \gg$ (recall that the action on the 
Lie algebra $\gg$ is induced by the adjoint representation $ad\rho$ 
of the group $\pi_1(M)$). 

Let $\M \subset \N \otimes \gg$ denote the subalgebra defined as:
$$
 \M  = (\N)^{\Phi} \hbox{~~is the space of~~} \Phi-\hbox{invariants.} 
$$
Let $\mu$ denote the restriction of $\nu$ to $\M$, the image of 
$\mu$ lies in the algebra of $\Phi$-invariants 
$(A^{\bul}\otimes\gg)^{\Phi}$. Similarly we let  $\L \subset \M$ denote 
the kernel of $\t\eps\circ \mu$.

\begin{lem}
The homomorphisms $\mu: \M \to (A^{\bul}\otimes \gg)^{\Phi}, 
\mu: \L \to (A^{\bul}\otimes \gg)_0^{\Phi}$ induce  
isomorphisms of $\H^0, \H^1$ and monomorphisms of $\H^2$.  
\end{lem}
\proof Standard. $\qed$ 

\begin{rem}
Note that $\H^0(\a^{\bul}(M, adP)_0)\cong \H^0((A^{\bul}\otimes 
\gg)^{\Phi}_0) \cong \H^0(\L)\cong 0$ and, under the assumption 
of Theorem \ref{con germ1}, 
$\H^0(\a^{\bul}(M, adP))\cong \H^0((A^{\bul}\otimes 
\gg)^{\Phi}) \cong \H^0(\M)\cong 0$. 
\end{rem}

\begin{cor}
The differential graded Lie algebra $\L$ controls the germ 
$(\Hom(\Ga, G), \rho)$ and (under the conditions of Theorem 
\ref{con germ1}) the differential graded Lie algebra $\M$ controls 
the germ $(X(\Ga, G), [\rho])$. 
\end{cor}

The split weight filtration on $\N$ defines a split weight filtration 
on $\N \otimes \gg$ (by taking tensor products of the 
components $\N_i$ with $\gg$). Using Proposition
 \ref{lif} restrict the split weight  filtration from $\N$ to  $\M$ 
and $\L$. Clearly these split filtrations of $\M$ and $\L$ satisfy the 
properties (a)-(c), where in (a) instead of $\wedge$ we take the Lie 
bracket. The property (d) can apriori fail, however we 
will need only its weak version:

(d') The restriction of the differential to $\M^1_1$ and $\L^1_1$ 
is identically zero. 

\medskip
\noindent Our further arguments are the same in the cases of $\M$ and $\L$ 
so we will discuss only $\M$.  

\medskip
Note that $\J_5:= \oplus_{i=5}^{\infty} \M_i$  
is an {\em ideal} in $\M$. We let $\J$ be the ideal generated by $\J_5$ and $\M^1_4$:
$$
\J = \J_5 \oplus \M^1_4 \oplus d(\M^1_4)
$$ 
The quotient $\M/\J$ is again a differential graded Lie algebra. 
The property (b) of $\M$ (weight restrictions 
on the cohomology groups) implies that the projection morphism
$\M \to {\cal Q}:= \M/\J$ 
induces an isomorphism of the 1-st and 2-nd cohomology groups. Hence by Theorem \ref{GolMil} the differential graded Lie algebras ${\cal Q}$ and $\M$ control 
germs which are analytically isomorphic. So it is enough to prove that 
${\cal Q}$ controls a quasi-homogeneous germ with the correct weights and we shall consider the \dgla ${\cal Q}$ from now on. 

\begin{rem}
We shall use the notation ${\cal Q}_p$ to denote 
the projection of $\M_p$ ($p\le 4$), elements of ${\cal Q}_p$  
will be denoted $\eta_p$. Lemma \ref{big} implies 
 that  ${\cal Q}= \oplus_p {\cal Q}_p$ so that:
\begin{itemize}
\item $d: {\cal Q}_p \to {\cal Q}_p$ for each $p$; 
\item $[\cdot , \cdot ]: {\cal Q}_p \otimes {\cal Q}_q \to {\cal Q}_{p+q}$\ ;
\item The induced weights on $\H^1({\cal Q}^{\bul})$ are $1, 2$ and the induced weights on $\H^2({\cal Q}^{\bul})$ are $2, 3, 4$; 
\item $d({\cal Q}^1_1) =0$. 
\item ${\cal Q}^0 = {\cal Q}_0= 0$. 
\end{itemize} 
\end{rem}

\noindent We split each vector space ${\cal Q}^k_p$ into the direct sum
${\cal H}^k_p \oplus B^k_p \oplus C_p^k$, where
 
\begin{itemize}
\item The space of coboundaries $B^k_p$ is the image of 
$d_p: {\cal Q}^{k-1}_p \to {\cal Q}^k_p$; 

\item  The space ${\cal H}^k_p$ of ``harmonic forms'' is 
a complement to $B^k_p$ in \newline 
$Z^k_p:= \ker(d_p: {\cal Q}^k_p \to {\cal Q}^{k+1}_p)$;

\item  $C_p^k$ is a complement to $Z^k_p$ in ${\cal Q}^k_p$. 
\end{itemize}

We let $\beta_p: {\cal Q}^k_p \to B^k_p$ denote the projection 
with the kernel ${\cal H}^k_p \oplus C_p^k$ (``coclosed $k$-forms''). We 
let $I_p: B^k_p \to C^{k-1}_p$ denote the inverse to the differential $d_p$. 
This allows us to define the ``co-differential''
$$
\delta_p: {\cal Q}^k_p \to C^{k-1}_p \subset {\cal Q}^{k-1}_p, \quad \delta_p= I_p\circ \beta_p
$$
(whose kernel is ${\cal H}^k_p \oplus C_p^k$). Let 
$\Pi_p: {\cal Q}_p \to {\cal H}_p$ denote the projection with the kernel  
$B_p \oplus C_p$. Clearly the projection
$$
{\cal H}^k_p \to \H^k({\cal Q}_p^{\bul})
$$
is an isomorphism of vector spaces. Notice that 
${\cal Q}^1 = {\cal Q}^1_1 \oplus {\cal Q}^1_2 \oplus {\cal Q}^1_3$.  Consider the variety 
$V\subset {\cal Q}^1$ given by the equation
\begin{equation}
\label{e1}
d\eta + [\eta, \eta]/2 =0 , \quad \eta \in {\cal Q}^1
\end{equation}
The algebra ${\cal Q}$ controls the germ $(V, 0)$ (see \S \ref{dgla}) 
and our goal is to show that this 
germ is quasi-homogeneous with correct weights. 
Since $\eta= \eta_1 + \eta_2 + \eta_3$ the equation (\ref{e1}) 
is equivalent to the system:
$$
d\eta_1 = 0, \ \ d\eta_2 + [\eta_1, \eta_1]/2=0, \ \ 
d\eta_3 + [\eta_1, \eta_2]=0, \ \ 
[\eta_1, \eta_3] + [\eta_2, \eta_2]/2= 0
$$  
Recall that the differential $d$ is identically zero on ${\cal Q}^1_1$ and the  
restriction
$$
d:{\cal Q}^1_3 \to {\cal Q}^2_3
$$
has zero kernel (since $\H^1({\cal Q}^{\bul})$ has no weight 3 elements). Therefore 
the equation $d\eta_3 + [\eta_1, \eta_2]=0$ is equivalent to the system:
$$
[\eta_1, d\eta_2]=0, \quad \Pi_3[\eta_1, \eta_2]=0, \quad \eta_3 + \del_3 [\eta_1, \eta_2]=0
$$
Note however that the equation $d\eta_2 + [\eta_1, \eta_1]=0$ together with 
the graded Lie identity imply that $[\eta_1, d\eta_2]=0$. 
Thus we eliminate the variable $\eta_3$ and the system of equations 
(\ref{e1}) is equivalent to:
\begin{equation}
\label{e2}
 d\eta_2 + [\eta_1, \eta_1]/2=0,\  \ 
  \Pi_3[\eta_1, \eta_2]=0,  \ \  
-[\eta_1, \del_3 [\eta_1, \eta_2]] + [\eta_2, \eta_2]/2= 0
\end{equation}
The mappings $d, \del_3, \Pi_3$ are linear and the bracket 
$[\cdot , \cdot]$ is quadratic. We conclude that the system of 
equations (\ref{e2}) is quasi-homogeneous where 
the weights of the generators (i.e. the components of) $\eta_j$ are $j=1,2$ 
and the weights of the relations are $2, 3$ and $4$. The only problem is that 
the first polynomial equation has nonzero linear term. To resolve this problem 
we let $\eta_2 = \eta'_2 + \eta''_2$, where 
$\eta'_2 \in Z^1_2, \eta''_2 \in C^1_2$. 
Thus (similarly to the case of $\eta_3$) we get:

the equation $d\eta_2 + [\eta_1, \eta_1]/2 =0$ is equivalent to the system:
$$
\Pi_2[\eta_1, \eta_1]=0, \quad \eta_2'' + \del_2 [\eta_1, \eta_1]/2=0
$$
Thus we have $\eta_2= \eta_2' - \del_2 [\eta_1, \eta_1]/2$ and 
instead of the system of equations (\ref{e2}) we get the system
$$
 \Pi_2[\eta_1, \eta_1]=0, 
  \Pi_3[\eta_1, \eta_2' - \del_2 [\eta_1, \eta_1]]=0, $$
$$ 
-[\eta_1, \del_3 [\eta_1, \eta_2' - \del_2 [\eta_1, \eta_1]/2 ]] + 
[\eta_2 ' - \del_2 [\eta_1, \eta_1]/2, \eta_2' - 
\del_2 [\eta_1, \eta_1]/2]/2= 0
$$
which is quasi-homogeneous with the required weights. 
Theorem \ref{hai} follows. $\qed$

\section{Malcev Lie algebras of Artin groups}
\label{completion}

Out discussion of the  material below follows \cite{Kah}. 
Let {\bf k} be 
a field of zero characteristic and $\Ga$ be a group. We define the 
{\bf k}-{\em unipotent completion} (or {\em Malcev completion}) $\Ga\otimes${\bf k} of $\Ga$  by the following universal property:

\begin{itemize}
\item There is a homomorphism $\eta\otimes $ {\bf k} $:\Ga \to \Ga\otimes 
${\bf k}. 
\item For every {\bf k}-unipotent Lie group $U$ and any homomorphism 
$\rho: \Ga\to U$ there is a lift $\t\rho: \Ga\otimes ${\bf k}$ \to U$ 
so that $\t\rho \circ (\eta\otimes ${\bf k}$)= \rho$. 
\item $\Ga\otimes ${\bf k} and $\eta\otimes ${\bf k} are unique up to 
an isomorphism.
\item  $\Ga\otimes\k$ is $\k$-pro-unipotent. 
\end{itemize}

\begin{rem}
Recall that any group $U$ above is torsion-free and nilpotent. 
\end{rem}

\begin{defn}
The group $\Ga\otimes${\bf k} %is {\em pro}-{\bf k}-{\em unipotent} and it 
has a {\em pro}-{\bf k}-{\em nilpotent} Lie algebra ${\cal L}(\Ga, \k)$. 
This algebra is called the {\bf k}-Malcev Lie  algebra of $\Ga$. 
\end{defn}

We will take $\k=\R$ in what follows. Thus we shall denote 
${\cal L}(\Ga):= {\cal L}(\Ga, \R)$, etc. 

\begin{ex}
Suppose that the group $\Ga$ has a generating set consisting of elements of finite order. Then  ${\cal L}(\Ga)=0$ because $\Ga\otimes\R= \{1\}$. In particular, if $\Ga$ is a Shephard group where all vertices have nonzero 
labels then  ${\cal L}(\Ga)=0$.  
\end{ex}

Let $H$ be a finite-dimensional real vector space, $L(H)$ is the free 
Lie algebra spanned by $H$. It can be described as follows. Consider the 
tensor algebra $T(H)$ of tensors of all possible degrees on $H$, define the 
Lie bracket of $T(H)$ by $[u, v] =u\otimes v- v\otimes u$. Then $L(H)$ is  
the Lie subalgebra in $T(H)$ generated by elements of $H$.

Let $F_r$ be a free group of rank $r$ and $H$ be the $r$-dimensional 
real  vector-space, then
$$
L(H)\cong {\cal L}(F_r) \quad .
$$
An element $u\in L(H)$ is said to have the degree $\le d$ if 
$u\in \oplus_{i\le n} H^{\otimes i}$. The {\em degree}  of $u$ equals 
$d$ if $deg(u)\le d$ but $deg(u)$ is not $\le d-1$. I.e. the degree of 
$u$ is the highest degree of monomial in the expansion of $u$ as a 
linear combination of tensor products of elements of $H$. For instance, 
{\em quadratic elements} of $L(H)$ are elements of the degree $2$, i.e. 
they have the form of nonzero linear combinations
$$
\sum_j [u_j, v_j], \quad u_j~, v_j~\in H
$$  

A {\em quadratically presented} Lie algebra is the quotient $L(H)/J$ where 
$J$ is an ideal generated by a (possibly empty) set of quadratic elements.  

\begin{thm}
(P.\ Deligne, P.\ Griffith, J.\ Morgan, D.\ Sullivan, \cite{DGMS}.) 
 Suppose that $M$ is a compact K\"ahler manifold, then the Malcev Lie 
algebra ${\cal L}(\pi_1(M))$ is quadratically presented. 
\end{thm}

\begin{thm}
(J.\ Morgan, \cite{Morgan1}, \cite{Morgan2}. ) {\bf ``Morgan's test.''} 
Suppose that $M$ is a smooth complex algebraic variety. Then the 
Malcev Lie algebra 
${\cal L}(\pi_1(M))$ is the quotient $L(H)/J$, where $L(H)$ is a free 
Lie algebra and the ideal $J$ is generated by elements of degrees 
$2\le d \le 4$.  
\end{thm}

\begin{rem}
Until now Morgan's theorem was the only known restriction on fundamental 
groups of smooth complex algebraic varieties, besides finite presentability. 
Much more restrictions are known in the case of smooth {\bf complete} 
varieties and compact K\"ahler manifolds, see \cite{Kah}.  
\end{rem}

Below we compute Malcev algebras of Artin groups. Suppose that 
$G^a$ is an Artin group. Let $n$ be the number of generators of 
$G^a$. Define a Lie algebra over $\R$ 
$$
{\cal L}:= \< X_1,..., X_n | [X_i, X_j]=0 \hbox{~if~} \eps(i,j)\ne \8 
\hbox{~is even~},  X_i=  X_j \hbox{~if~} \eps(i,j)\ne \8 \hbox{~is odd~}
\>$$
where $[X, Y]$ denotes the Lie algebra commutator. Clearly this Lie 
algebra is quadratically presentable. 

To compute Malcev completions we will need the following two lemmas 

\begin{lem}
\label{LMa}
Suppose that $\rho: G^a\to N$ is a homomorphism to a torsion-free 
 $s$-step nilpotent group. Then for all $x_i, x_j$ such 
that $2q_{ij}= \eps(i,j)\ne \8$ we have  $[\rho(x_i), \rho(x_j)]=\1$. 
 For all $x_i, x_j$ such that $2q_{ij} +1= \eps(i,j)\ne \8$ we 
have  $\rho(x_i)= \rho(x_j)$.
\end{lem} 
\proof The assertion is obvious if $N$ is Abelian. So we assume that 
the assertion is valid for all $(s-1)$-step nilpotent torsion-free 
groups $\bar{N}$. Let $Z(N)$ be the center of $N$, let 
$\bar{N}:= N/Z(N)$ and $p: N\to \bar {N}$ be the projection. 
Then by the induction hypothesis:

\begin{itemize}
\item $[p(\rho(x_i)), p(\rho(x_j))]=\1$, provided that 
$2q_{ij}= \eps(i,j)\ne \8$. 
\item $p(\rho(x_i))= p(\rho(x_j))$, provided that 
$2q_{ij} +1= \eps(i,j)\ne \8$. 
\end{itemize}

(1) Consider the case $2q_{ij}= \eps(i,j)$. Then 
$\rho(x_i) \rho(x_j)= \rho(x_j) \rho(x_i)z $, for some  $z\in  Z(N)$. 
Thus the relation 
$$
(x_i x_j)^{q_{ij}} =   (x_j x_i)^{q_{ij}}
$$
implies that
$$
z^{q_{ij}} \rho[(x_j x_i)^{q_{ij}}]=  \rho[(x_j x_i)^{q_{ij}}] 
$$
Since $N$ is torsion-free we conclude that $\rho( [x_i, x_j])= \1$. 

(2) Another case is when $2q_{ij} +1= \eps(i,j)$. Then 
$\rho(x_i)=z \rho(x_j)$ 
, for some  $z\in  Z(N)$. The relation
$$
(x_i x_j)^{q_{ij}} x_i =   (x_j x_i)^{q_{ij}} x_j 
$$
implies that $z=1$. $\qed$ 

\medskip
\begin{rem}
In our paper we use only Artin groups with even labels. 
\end{rem}

\begin{lem}
\label{Lcom}
Suppose that $U$ is a unipotent group, $a, b\in U$ are commuting 
elements. Then 
$[\log(a), \log(b)]=0$ in the Lie algebra of $U$. 
\end{lem}
\proof Since $U$ is unipotent we can think of $U$ as the subgroup of 
the group of upper-triangular matrices with $1$-s on the diagonal. 
Then for any $g\in U$ we have: $\log(g)= \log(\1 - (\1- g))$, 
$h= \1- g$ is a nilpotent matrix, thus $\log( \1 -h)$ is a polynomial 
of $h$. Since matrices $a, b$ commute, then any polynomial functions 
of them commute as well. Thus $[\log(a), \log(b)]=0$. 
$\qed$

\begin{thm}
\label{quad}
Under the above conditions ${\cal L}\cong {\cal L}(G^a)$ is the Malcev 
Lie algebra of $G^a$. 
\end{thm}
\proof Let $F$ denote the free group on $x_1,..., x_n$ and $\pi: F\to G^a$ 
be the quotient map. Let $F\otimes \R$ be the $\R$-unipotent completion of 
$F$ and $\eta: F\to F\otimes \R$ be the canonical homomorphism. Let 
$\L(F, \R)$ be the Lie algebra of $F\otimes \R$. Put $g_i= \eta(x_i)$, 
$1\le i\le n$ and $X_i:= \log(g_i)$. Let ${\cal I}$ be the ideal in 
$\L(F, \R)$ generated by the commutators $[X_i, X_j]$, $\eps(i, j)\ne \8$ 
is even and the elements 
$X_i-  X_j$ if $\eps(i,j)\ne \8$ is odd. Let ${\cal Q}:= 
\L(F, \R)/{\cal I}$ be the quotient Lie algebra and $Q$ be the 
corresponding pro-unipotent Lie group over $\R$. Let 
$\hat{\pi}: F\otimes \R \to Q$ be the quotient map. Put 
$\bar{g}_i := \hat{\pi}(g_i)$ and $\bar{X}_i:= d\hat{\pi}(X_i)$. 
Then $\bar{g}_i  = \exp(\bar{X}_i)$, $1\le i\le n$. Consequently 
$[\bar{g}_i, \bar{g}_j]= \1$ for all vertices $i,j$ connected by an 
edge with even label and $\bar{g}_i= \bar{g}_j$ for all vertices 
$i,j$ connected by an edge with odd label. Hence we have a 
commutative diagram
\[
\begin{array}{ccccc}
F & \stackrel{\eta}{\longrightarrow} & F\otimes \R & \lra & \L(F, \R) \\
\pi {\Big\downarrow} & ~ & \hat{\pi}{\Big\downarrow} & ~ & {\Big\downarrow}\\
G^a & \stackrel\tau{\longrightarrow} & Q & \lra & {\cal Q}\\
\end{array}
\]
where $\tau(\pi(x_i))= \bar{g}_i$, $1\le i\le n$. We claim that $Q$ 
is the Malcev completion of $G^a$. It is clear that any homomorphism 
$\rho: Q\to U$ from $Q$ to a unipotent group $U$ is determined by 
its pull-back to $G^a$ (because its pull-back to $F\otimes \R$ is 
determined by its further pull-back to $F$). So let $\rho: G^a\to U$ 
be a homomorphism with $U$ a unipotent group over $\R$. The 
homomorphism $\pi^* \rho$ extends to a morphism $\hat{\rho}: 
F \otimes \R \to U$.  Note that $U$ is necessarily nilpotent 
and torsion-free.

\medskip 
According to Lemma \ref{LMa} for each pair of  
$(i,j)$ such that $\eps(i, j)$ is even we have:
$$
\hat\rho ([g_i, g_j])= \pi^* \rho([x_i, x_j])= \rho([x_i, x_j])= \1
$$
Hence $\log(\hat\rho(g_i))$ and  $\log(\hat\rho(g_j))$ commute in 
the Lie algebra $u$ of $U$ (see Lemma \ref{Lcom}). Therefore 
$d\hat\rho(X_i)$ and  
$d\hat\rho(X_j)$ commute in $u$. The case of odd labels  
$\eps(i, j)$ is similar. This implies that $d\hat\rho$ 
descends to ${\cal Q}$ and consequently $\hat\rho$ descends 
to $Q$. $\qed$ 

\begin{cor}
\label{squad}
If $G^a$ is any Artin group then ${\cal L}(G^a)$ is quadratically presented. 
\end{cor}

Thus the Artin groups constructed in Theorem \ref{T10.3} satisfy 
 {\em Morgan's test} of being fundamental groups of smooth 
complex algebraic varieties.

\section{Representation varieties near the trivial representation}
\label{triv}

The second author would like to thank Carlos Simpson for explaining 
Theorem \ref{clarify} in this section.  
 
Let $\Ga$ be a finitely-generated group, $\Ga \otimes \R$ is 
its Malcev completion, $\eta\otimes \R: \Ga \to   
 \Ga \otimes \R$ is the canonical homomorphism. Let 
${\cal L}(\Ga)$ denote the Malcev Lie algebra of $\Ga$. Let 
$\hbox{\bf G}$ be the set of real points of an algebraic 
group $G$ defined over $\R$, $\gg$ be the Lie algebra of $\hbox{\bf G}$. 
The homomorphism $\eta\otimes \R$ induces the pull-back morphism 
$\eta^*: \Hom(\Ga\otimes\R, \hbox{\bf G}) \to \Hom(\Ga, \hbox{\bf G})$. 
Let $\rho_0: \Ga \otimes\R\to \hbox{\bf G}$ be the trivial 
representation .  

\begin{thm}
\label{clarify}
\begin{enumerate}
\item If the Lie algebra ${\cal L}(\Ga)$ is quadratically 
presentable then the variety $\Hom({\cal L}(\Ga), \gg)$ 
is given by homogeneous quadratic equations. 
\item The varieties $\Hom(\Ga\otimes\R, \hbox{\bf G})$ and  
$\Hom({\cal L}(\Ga), \gg)$ are naturally isomorphic. 
\item The morphism $\eta^*$ induces an isomorphism of germs
$$
(\Hom(\Ga\otimes\R, \hbox{\bf G}),\rho_0) \longrightarrow
( \Hom(\Ga  , \hbox{\bf G}), \eta^*(\rho_0))
$$
\end{enumerate}
\end{thm}
\proof The property (1) is obvious. To prove (2) note 
that the group $\Ga\otimes\R$ is $\R$-pro-unipotent, thus 
we have a natural isomorphism between the representation variety 
of $\Ga\otimes\R$ and of its Lie algebra. 

Now consider (3). Let ${\cal A}$ be an Artin local $\R$-algebra. 
We recall  that $G({\cal A})$ is the set of 
${\cal A}$-points of $G$, algebraically the group $G({\cal A})$ is 
the semidirect product $N_{\a} \rtimes \hbox{\bf G}$, where 
$N_{\a}$ is a certain $\R$-unipotent group (kernel of the 
natural projection  
$p_0: G({\cal A})\to \hbox{\bf G}$). Consider the space
$$
\Hom_0(\Ga, G({\cal A})):= \{ \rho: \Ga \to G({\cal A}) |           
p_0(\rho)= \rho_0\}\cong \Hom(\Ga, N_{\a})
$$
Thus, by the definition of $\Ga\otimes\R$, for all Artin local 
$\R$-algebra $\a$ the morphism $\eta^*$ induces a 
 natural bijection between the $\R$-points of the varieties 
$\Hom(\Ga, N_{\a})$ and $\Hom(\Ga\otimes\R, N_{\a})$. Then we have an 
induced isomorphism between functors 
$$
\Hom(\Ga, G) : \a \mapsto set\ ,  \quad 
\Hom(\Ga\otimes\R, G) : \a \mapsto set
$$ 
of $\a$-points. Therefore (by \cite{GM}) $\eta^*$ induces an  
isomorphism of the germs 
$$
(\Hom(\Ga\otimes\R, \hbox{\bf G}),\rho_0~) \longrightarrow
( \Hom(\Ga  , \hbox{\bf G}), \eta^*(\rho_0))
$$
$\qed$ 

\begin{prop}
Let $\Ga$ be a Coxeter group or a Shephard group (where all vertices 
have labels). Then the 
trivial representation $\rho_0: \Ga\to \hbox{\bf G}$ is 
infinitesimally rigid (and hence is an isolated reduced point) 
in $\Hom(\Ga, \hbox{\bf G})$. 
\end{prop}
\proof The group $\Ga$ is generated by elements of finite order. 
Let $\zeta$ be a cocycle in $Z^1(\Ga, \gg)$. Then 
$\zeta|_{\<g_j\>}$ is a coboundary for each generator $g_j$ 
of $\Ga$ (since $g_j$ has finite order). However $\rho_0(g_j)=\1$, 
hence $\zeta(g_j)=0$ for all $j$. We conclude that $\zeta=0$. 
$\qed$

\begin{thm}
 Let $\Ga$ be any  Artin group. Then the representation 
variety 
$\Hom(\Ga, \hbox{\bf G})$ has at worst quadratic singularity at 
the trivial representation. 
\end{thm}
\proof Combine Corollary \ref{squad} and Theorem \ref{clarify}. 
$\qed$ 

%\newpage

\noindent Michael Kapovich: Department of Mathematics, University of 
Utah,  Salt Lake City, UT 84112, USA ; kapovich$@$math.utah.edu

\smallskip
\noindent John J. Millson: Department of Mathematics, University of 
Maryland, College Park, MD 20742, USA ; jjm$@$julia.umd.edu


\begin{thebibliography}{BaBE}
\addcontentsline{toc}{section}{Bibliography}

\bibitem[ABC]{Kah}
J.~Amoros, M.~Burger, K.~Corlette, D.~Kotschick, D.~Toledo, 
``Fundamental groups of K\"ahler manifolds'', AMS Mathematical 
Surveys and Monographs, Vol. {\bf 44}, 1996. 
 
\bibitem[A1]{Arnold1}
V.~I.~Arnold, {\em Normal forms of functions in \nbds of degenerate 
critical points},  In: Lecture Notes of London Math. Soc., Vol. 
{\bf 53}, ``Singularity Theory'', 91--131. 

\bibitem[A2]{Arnold2}
V.~I.~Arnold, {\em Critical points of smooth functions and their normal 
forms,} In: Lecture Notes of London Math. Soc., Vol. {\bf 53}, 
``Singularity Theory'', 132--206. 

\bibitem[AM]{AM}
M.\ Atiyah, I.\ MacDonald, ``Introduction to Commutative Algebra'', 
Addison-Wesley, 1969.  

%\bibitem[Ba]{Bauer}
%?. Bauer, preprint. 

\bibitem[BiM]{BM}
E.\ Bierstone, P.\ Millman, {\em Canonical desingularization in 
characteristic zero by blowing up maximal strata of a local 
invariant}, preprint. 

\bibitem[BuM]{BuM}
R.\ O.\ Buchweitz, J.\ J.\ Millson, {\em CR-geometry  and deformations of 
isolated singularities,} Memoirs of AMS (to appear). 

\bibitem[B]{Bri}
E.\ Brieskorn, {\em Die fundamentalgruppe des raumes der re\-gu\-la\-ren 
or\-bits einer endlichen komplexen spiegelungsgruppe,} Inv. Math., 
Vol. {\bf 12} 
(1971) 57--61. 

\bibitem[C]{Coxeter}
H.\ S.\ M.\ Coxeter, {\em Finite unitary groups generated by reflections,} 
Abhandlunden aus dem Mathematischen Seminar der Universitat Hamburg, 
Bd. {\bf 31} (1967) 125-- 135. 

\bibitem[D1]{D}
P.\ Deligne, {\em Theorie de Hodge II,} Math. Publ.  of IHES, 
Vol. {\bf 40} (1971) 5-- 58. 

\bibitem[D2]{D2}
P.\ Deligne, {\em Theorie de Hodge III,} Math. Publ.  of IHES, 
Vol. {\bf 44} (1974) 5-- 77. 

\bibitem[DGMS]{DGMS}
P.\ Deligne, P.\ Griffith, J.\ Morgan, D.\ Sullivan, 
{\em Real homotopy theory of K\"ahler manifolds,} Inv. Math., Vol. {\bf 29} 
(1975) 245-- 274. 

\bibitem[DG]{DG}
M.\ Demazure, P.\ Gabriel, ``Groupes algebriques: Vol. I. Geometrie 
algebrique, g\'en\'eralit\'es, groupes commutatifs,'' Paris: Masson, 1970. 

\bibitem[Di]{Dimca}
A.\ Dimca, ``Topics on Real and Complex Singularities'', Advanced Lectures in Mathematics, Vieweg, 1987.  

\bibitem[EH]{EH}
D.\ Eisenbud, J.\ Harris, ``Schemes: the language of the modern 
algebraic geometry'', Wadsworth \& Brooks/Cole Math Series, 1992. 

\bibitem[EN]{Neu}
D.~Eisenbud, W.\ Neumann, ``Three-dimensional link theory and invariants 
of plane curve singularities'', Ann. of Math. Stud., Princeton Univ. Press, 
Vol. {\bf 110} (1985). 

\bibitem[GM]{GM}
 W.\ Goldman, J.\ J.\ Millson, 
{\em The deformation theory of representations of
 fundamental groups of compact K\"ahler manifolds}, Math. Publ. of 
IHES, Vol. {\bf  67} (1988) 43-- 96.

\bibitem[GrM]{GrM}
P.\ Griffiths, J.\ Morgan, ``Rational homotopy theory and differential forms,  Birkhauser'', Progress in Mathematics, 1981. 

\bibitem[Hai]{Hain}
R.\ Hain, in preparation. 

\bibitem[H]{Hart}
R.\ Hartshorne, ``Foundations of Projective Geometry'', Benjamin Inc., 
NY, 1967. 

\bibitem[JM]{JM}
D.\ Johnson, J.\ J.\ Millson, {\em Deformation spaces associated to 
compact hyperbolic manifolds,} In: ``Discrete Groups in Geometry 
and Analysis'', Papers in honor of G.\ D.\ Mostow on his 60-th 
birthsday, R.\ Howe (ed.), Progress in Mathematics, Vol. {\bf 67}, 
Birkhauser; (1987) 48-- 106. 

\bibitem[KM1]{KM1}
 M.\ Kapovich, J.\ J.\ Millson,  {\em The relative deformation theory 
of representations and  flat connections and deformations of 
linkages in constant curvature spaces,} Compositio Math., Vol. {\bf 103}, 
N 3 (1996) 287--317.  

\bibitem[KM2]{KM2}
 M.\ Kapovich, J.\ J.\ Millson, {\em On the deformation theory of 
representations of fundamental groups of closed hyperbolic 
3-manifolds,} Topology, Vol. {\bf 35}, N 4 (1996) 1085--1106.

\bibitem[KM3]{KM3}
 M.\ Kapovich, J.\ J.\ Millson, {\em Hodge theory and the 
art of paper folding,} Publications of RIMS, Kyoto, Vol. 
{\bf 33}, N 1 (1997). 

\bibitem[Le]{Le}
H. van der Lek, {\em Extended Artin groups,} Proc. of Symp. in Pure
Math., Vol. {\bf 40} (1983), Part 2, 117--121. 

\bibitem[Lo]{Lo}
E.\ Looijenga, {\em Invariant theory for generalized root systems,} 
Inv. Math., Vol. {\bf 61} (1980) 1-- 32. 

\bibitem[LM]{LM}
A.\ Lubotzky, A.\ Magid, {\em Varieties of representations for finitely 
generated groups}, Memoirs of AMS, Vol. {\bf 336} (1985), N 5. 

\bibitem[Mi]{M} 
J.\ J.\ Millson, {\em Rational homotopy theory and deformation 
problems from algebraic geometry,} Proc. of ICM 1990, Vol. I, 549-- 558.

\bibitem[Mn]{Mnev}
N.\ Mn\"ev, {\em The universality theorems on the classification 
problem of configuration varieties and convex polytopes varieties}, 
Lecture Notes in Math, Vol. {\bf 1346} (1988) 527--543. 

\bibitem[Mo1]{Morgan1}
J.\ Morgan,  {\em Hodge theory for the algebraic topology of smooth 
algebraic varieties,} Proc. Symp. in Pure Math., Vol. {\bf 32} 
(1978), 119--127.  

\bibitem[Mo2]{Morgan2}
J.\ Morgan, {\em The algebraic topology of smooth algebraic varieties,} 
Math. Publ. of IHES, Vol. {\bf 48} (1978), 137-- 204. 

\bibitem[N]{Newstead}
P.\ E.\ Newstead, {\em Introduction to Moduli Problems and Orbit 
Spaces}, Tata Institute Lecture Notes. 

\bibitem[Sh]{Sh}
G.\ Shephard, {\em Regular complex polytopes,} Proc. London Math. Soc., 
Vol.  {\bf 2} (1952) 82--97. 

\bibitem[Si]{Si}
C.~Simpson, {\em Higgs bundles and local systems,} Publications of IHES, 
{\bf 75} (1992) 5--95. 

\bibitem[St]{St}
K.\ G.\ C.\ von Staudt, ``Beitr\"age zur Geometre der Lage'', Heft 2, 1857. 

\bibitem[Sul]{Sul}
D.\ Sullivan, {\em Infinitesimal computations in topology}, Publ. of IHES. 
vol. {\bf 47} (1977) 269--331.

\bibitem[Sum]{Su}
H.\ Sumihiro, {\em Equivariant completion,} J. Math. Kyoto Univ., 
Vol. {\bf 14} (1974), 1-- 28; Vol. {\bf 15} (1975), 573-- 605. 

\end{thebibliography}
\end{document}